\documentclass[11pt]{amsart}

\usepackage{DKstyle}
\usepackage{bm}
\usepackage{float}
\usepackage[toc,page]{appendix}

\newcommand{\col}{\: : \:}
\newcommand{\Mod}[1]{\ (\mathrm{mod}\ #1)}
\newcommand{\vertiii}[1]{{\left\vert\kern-0.25ex\left\vert\kern-0.25ex\left\vert #1
    \right\vert\kern-0.25ex\right\vert\kern-0.25ex\right\vert}}
 \makeatletter
\newcommand*{\rom}[1]{\expandafter\@slowromancap\romannumeral #1@}
\makeatother
\theoremstyle{plain}

\usepackage{caption}
\usepackage[labelfont=rm]{subcaption}

\usepackage[pagebackref=false, colorlinks]{hyperref}

\hypersetup{pdffitwindow=true,linkcolor=blue,citecolor=blue,urlcolor=blue,menucolor=blue}

\usepackage{comment}

\newcommand {\comm}[1]   {\textcolor{black}{#1}}
\def\d{\text{d}}

\subjclass{}%

\date{\today}%
\dedicatory{}%
\commby{}%

\title{Second moment of the light-cone Siegel transform and applications}
\author{Dubi Kelmer}
\address{Department of Mathematics, Boston College, Chestnut Hill MA 02467-3806, USA}
\email{kelmer@bc.edu}

\author{Shucheng Yu}
\address{School of Mathematical Sciences, University of Science and Technology of China (USTC), 230026, Hefei, China}
\email{yusc@ustc.edu.cn}
\thanks{D.K. and S.Y. were partially supported by NSF CAREER grant DMS-1651563. S.Y. was supported by the Knut and Alice Wallenberg Foundation}
\begin{document}
\begin{abstract}
We study the light-cone Siegel transform, transforming functions on the light cone of a rational indefinite quadratic form $Q$ to a function on the homogenous space $\SO^+_Q(\Z)\bk \SO^+_Q(\R)$. In particular, we prove a second moment formula for this transform for forms of signature $(n+1,1)$, and show how it can be used for various applications involving counting integer points on the light cone. In particular, we prove some new results on intrinsic Diophantine approximations on \comm{ellipsoids} as well as on the distribution of  values of random linear and quadratic forms on the light cone.
\end{abstract}
\maketitle

\tableofcontents

\section{Introduction}
The Siegel transform, introduced in \cite{Siegel1945}, transforms a function, $f$ on $\R^n$ of sufficient decay, to a function $Sf$ on $X_n$, the space of unimodular lattices in $\R^n$, by averaging over the lattice.
The space of lattices $X_n$ can be identified with the homogenous space $\SL_n(\Z)\bk \SL_n(\R)$ and is equipped with  a natural $\SL_n(\R)$-invariant probability measure; the Siegel  mean value theorem expresses the average of $Sf$ in terms of the average of $f$.  
Later, generalizing Siegel's result, Rogers \cite{Rogers1955} proved a general $k$-th moment formula for $Sf$ for $k$ up to $n-1$. 
While in general these formulas are quite involved, the second moment formula for $Sf$ is relatively simple, and can be expressed in terms of the first and (twisted) second moment of $f$. 
A masterful application of the second moment formula, was given by Schmidt \cite{Schmidt1960}, who proved very strong asymptotic formulas counting lattice points in an increasing family of sets. 

In recent years there have been some renewed interests in these moment formulas, which were used in a number of new applications for problems in geometry of numbers and Diophantine approximations.
To list just a few examples, in \cite{AthreyaMargulis09} Athreya and Margulis used the second moment formula to prove a random version of Minkowski's theorem;
 while in  \cite{AthreyaGhoshTseng2015}, by using Siegel's mean value theorem, together with Birkhoff ergodic
theorem, Athreya, Ghosh and Tseng proved a quantitative variant of Dirichlet's theorem with approximants restricted in certain fixed directions. 
This was extended in \cite{AlamGhoshYu2021} by  Alam, Ghosh and Yu, who proved a quantitative variant of Khintchine's theorem with a congruence condition. The strategy used in \cite{AlamGhoshYu2021} follows mostly that of \cite{AthreyaGhoshTseng2015}, but with the ergodic theorem replaced by Schmidt's lattice point counting arguments which rely on moment formulas of Siegel transforms. 
 In another direction, Athreya and Margulis \cite{AthreyaMargulis2018} found another novel use of the second moment formula to study value distribution of random quadratic forms at integers points improving on previous results of 
\cite{GhoshKelmer2018,GhoshGorodnikNevo2020} using homogenous dynamics.  The method of Athreya and Margulis was further refined by the authors in \cite{KelmerYu2020}, giving an asymptotic counting formula for the number of integer solution of a generic form in a shrinking interval, as well as treating more general homogenous polynomials of higher degrees. Similar results were later obtained for more general functions in \cite{BandiGhoshHan2021,KleinbockSkenderi2021,GhoshKelmerYu2022} using the second moment formula as a main ingredient. Finally we mention the results of \cite{BjorklundGorodnik2023,BjorklundGorodnik2023b} in which these moment formulas were a crucial ingredient for proving various limiting 
laws for counting lattice points and some related problems in geometry of numbers.
 
Following this renewed interest in moment formulas for the Siegel transform, in  \cite{KelmerYu2021} it was observed that the Siegel transform can be seen as an incomplete Eisenstein series, and one can use the spectral theory of Eisenstein series to give an alternative proof of the second moment formula. The advantage of this viewpoint is that it opens the way to prove similar moment formulas when the homogenous space $\SL_n(\Z)\bk \SL_n(\R)$ is replaced by other homogenous spaces. To showcase this approach  in \cite{KelmerYu2021} the authors proved a new second moment formula for Siegel transforms restricted to the space of symplectic lattices $\Sp_{2n}(\Z)\bk \Sp_{2n}(\R)$. More recently, using this spectral approach, together with Ghosh \cite{GhoshKelmerYu2022} we proved a new second moment formula for Siegel transforms defined on a certain congruence cover of the space of lattices. 
We mention that in all these new moment formulas, the Siegel transforms under consideration are all defined such that they take a function on the Euclidean space to a function on the corresponding homogeneous space. In this paper we apply this spectral approach to study a new type of Siegel transform where the Euclidean space is replaced by the light cone of a quadratic form.
\comm{Below we prove first and second moment formulas for this Siegel transform and use these moment formulas to prove some new results on intrinsic Diophantine approximation on ellipsoids as well on the value distribution of random linear and quadratic forms on the light cone.}

\subsection{Light-cone Siegel transform and its moment formulas}
Fix an integer $n\geq 1$, and let $Q: \R^{n+2}\to \R$ be a rational $\Q$-isotropic quadratic form of signature $(n+1,1)$. \comm{Recall that $Q$ is \textit{$\Q$-isotropic} means that there exists a nonzero $\bm{v}\in \Q^{n+2}$ such that $Q(\bm{v})=0$; note that
in view of Meyer's theorem (see e.g. \cite[p. 43]{Serre1973}}) this assumption is  satisfied as soon as $n\geq 3$.
The light cone of $Q$ is defined as 
$$
\cV_Q:=\left\{\bm{v}\in \R^{n+2}\setminus \{\bm{0}\} : Q(\bm{v})=0\right\}
.
$$
It has two connected components and we let $\cV_Q^+$ be one \comm{of these}  \comm{components}; \comm{see section \ref{sec:cormea} for more details}. Since $Q$ is assumed to be $\Q$-isotropic, the set of integer points $\cV_Q^+(\Z):=\cV_Q^+\cap \Z^{n+2}$ is infinite. We also let $\cL_Q\subseteq \cV_Q^+(\Z)$ be the subset of primitive integer points. 
Let $G=\SO_Q^+(\R)$ be the identity component of the special orthogonal group preserving $Q$ and let $\G=\SO_{Q}^+(\Z):=G\cap \SL_{n+2}(\Z)$ be the subgroup of integer points (which is a non-uniform lattice of $G$ since $Q$ is $\Q$-isotropic).
Note that $G$ acts transitively on $\cV_Q^+$ via right multiplication and that the \comm{induced} action of $\G$ \comm{on $\cV_Q^+$} preserves $\cL_Q$. We can thus parameterize the space of $\cL_Q$-translates  by the homogeneous space $Y_Q=\G\bk G$ via the map $\G g\in \G\bk G\mapsto \cL_Q g$. Generalizing the classical Siegel transform, given a function $f: \cV_Q^+\to \C$ of sufficient decay, its \textit{light-cone Siegel transform}, is a function on $Y_Q$, defined for any $g\in G$ by
\begin{align*}
S_Qf(g):=\sum_{\bm{v}\in \cL_Q }f(\bm{v}g).
\end{align*}
Our main goal of this paper is to prove moment formulas, relating moments of $S_Qf$ to integrals of the original function $f$ defined on $\cV_Q^+$. 

When taking the light cone Siegel transform of the function $f(\bm{v})=\|\bm{v}\|_Q^{-s}$ for $\|\cdot\|_Q$ an appropriate norm (see \eqref{equ:qnorm} below) and $\Re(s)>n$, we get the \textit{light-cone Eisenstein series}, $E_Q(s,g)$, that was studied in  \cite{KelmerYu2022a}. In particular, it was shown there that $E_Q(s,g)$ has an analytic continuation with a simple pole at $s=n$ with constant residue 
$$
\omega_Q:=\Res_{s=n}E_Q(s,g)
$$ 
and at most one more exceptional simple pole in the half plane $\Re(s)\geq \frac{n}{2}$ located at $s_n:=\left \lfloor{\tfrac{n+2}{2}}\right \rfloor$ with residue $\Res_{s=s_n}E_Q(s,g)\in L^2(Y_Q)$ 
(when $n=1,2$ we have $s_n=n$ and there is no exceptional pole); \comm{see section \ref{sec:lightconeeisen} below for more details}. Moreover, also for $n\geq 3$  it is shown that, in many examples, there is no  exceptional pole; see \rmkref{rmk:noexpole} below.

Before we state the moment formulas, we need to fix some measures. Let $\mu_Q$ be the unique $G$-invariant probability measure on $Y_Q$. Similarly, there is also a unique (up to scaling) $G$-invariant measure on $\cV_Q^+$, which we fix and denote by $m_{\cV^+_Q}$ (see \eqref{equ:lebmea} below). One should think of this measure as the counterpart of the Lebesgue measure in the classical Siegel transform setting. \comm{For simplicity of notation, for any function $f$ on $\cV_Q^+$ we will abbreviate the space average $\int_{\cV_Q^+}f(\bm{v})\ \d m_{\cV_Q^+}(\bm{v})$ simply by $m_{\cV_Q^+}(f)$.} With these measures \comm{and notation} we have the following.

\begin{Thm}\label{thm:mainrest}
Let $f: \cV_Q^+\to \C$ be \comm{measurable}, bounded and compactly supported. Then we have
\begin{align}\label{equ:mainthmfirst}
\int_{Y_Q}S_Qf(g)\,\d\mu_Q(g)=\omega_Q m_{\cV_Q^+}(f).
\end{align} 
\comm{Further assume $f$ is smooth, then}
\begin{align}\label{equ:mainthmsecond}
\int_{Y_Q}\left|S_Qf(g)\right|^2\, \d\mu_Q(g)=|\omega_Qm_{\cV_Q^+}(f)|^2+c_Q M_{f,f}\left(s_n\right)+O\left(m_{\cV_Q^+}(|f|^2)\right),
\end{align}
where $s_n=\left \lfloor{\tfrac{n+2}{2}}\right \rfloor$, the term $M_{f,f}(s)$ is a quadratic form on $f$ given in \eqref{def:mffun} and $c_Q\geq 0$ with $c_Q=0$ if and only if $E_Q(s,g)$ has no exceptional pole.
\end{Thm}

\begin{rem}\label{rmk:noexpole}
Since some of our applications are nicer when there is no secondary term, \comm{i.e. when $c_Q=0$}, it is worthwhile to give some explicit examples for which we know this is the case. First, when $n=1,2$ there is no exceptional pole for any form. Next, for the standard forms
\begin{equation}\label{equ:stformqn}
Q_n(\bm{v})=\sum_{j=1}^{n+1}v_j^2-v_{n+2}^2
\end{equation}
it follows from \cite[Theorem 4.4]{KelmerYu2022a}  that $c_{Q_n}> 0$ if and only if $n>1$ with  $n\equiv 1\Mod{8}$. In fact, the same holds for the forms $Q_{n,d}(\bm{v})=\sum_{j=1}^{n+1}v_j^2-d^2v_{n+2}^2$ with $d$ odd and square-free.
For a general rational $\Q$-isotropic form, $Q(\bm{v})=\bm{v}^tJ\bm{v}$, of signature $(n+1,1)$ it follows from  \cite[Theorem 1.8]{KelmerYu2022a} that  $c_Q=0$ when $n\equiv 0\Mod{4}$ and $\sqrt{-\det(J)}$ is irrational or when $n\equiv 2\Mod{4}$.
\end{rem}
 
For our applications to counting problems, we take the function $f$ to be an indicator function of a measurable set $B\subseteq \cV_Q^+$ of finite measure and, \comm{after approximating $f$ by smooth functions}, we get  mean square bound for the discrepancy 
\begin{equation} 
D(\cL_Qg,B)=|\#(\cL_Q g\cap B)-\omega_Qm_{\cV_Q^+}(B)|.
\end{equation} 
When $c_Q=0$, there is no secondary term in our second moment formula and our variance bound is sharp and holds for any \comm{finite-measure} set. When $c_Q> 0$ some work is needed in order to control the secondary term. While this secondary term is given by some explicit formula, it is in general difficult to bound without imposing additional restrictions on the set $B$ (e.g., assuming that it is given by dilation of some fixed set or a generalized sector  (see \defref{def:gesec})). This restriction is responsible for some constraints in our later applications. 
In the following, let 
\begin{align}\label{equ:defbq}
\beta_Q &:=\left\{\begin{array}{ll}   
	1 & c_Q=0,\\
	\frac{2s_n}{n} & c_Q> 0.
	\end{array}\right.
\end{align}
We prove the following variance bound.

\begin{Cor}\label{cor:vabd}
\comm{Let $B\subseteq \cV_Q^+$ be a Borel set of finite measure. If $c_Q>0$ we further assume that $B$ is a generalized sector.  We then have} 
       \begin{align}\label{equ:disbdwiexpo}
			\int_{Y_Q}\left|D(\cL_Q g, B)\right|^2 \, \d\mu_{Q}(g)\ll_Q 
				m_{\cV_Q^+}(B)^{\beta_Q}+m_{\cV_Q^+}(B).
		\end{align}
\end{Cor}
\begin{remark}
We note that for the case when $c_Q> 0$, 
the exponent $\beta_Q$ is optimal. In fact, given any family of dilation\comm{s} $\{tB\}_{t>1}$, their measure grows like  $m_{\cV_Q^+}(tB)=t^n m_{\cV_Q^+}(B)$ and the secondary terms $M_{\chi_{tB},\chi_{tB}}(s)=t^{\frac{2s}{n}}M_{\chi_B,\chi_B}(s)$  for any $s\in (\frac{n}{2}, n)$ (see \rmkref{rmk:scalar}). From this and the moment formulas in \thmref{thm:mainrest} we get that when $c_Q> 0$, for any $t>1$,
\begin{align*}
\int_{Y_Q}D(\cL_Qg,tB)^2 \,\d\mu_Q(g)\asymp_{B} m_{\cV_Q^+}(tB)^{\beta_Q}.
\end{align*}
\end{remark}

\subsection{Applications to counting points on the light cone}
We now discuss several applications of our moment formulas to counting primitive integer points on the light cone. \comm{Fix $\beta_Q$ as in \eqref{equ:defbq}.}
First, we can apply Schmidt's arguments with the variance bounds in \corref{cor:vabd} to get the following result on counting points of generic $\cL_Q$-translates in increasing sets. \comm{Here a family of Borel sets $\{B_t\}_{t>0}$ in $\cV_Q^+$ is called \textit{increasing} if $B_{t_1}\subseteq B_{t_2}$ whenever $t_1<t_2$.}

\begin{Thm}\label{thm:countschd}
Let $\{B_t\}_{t>0}\subseteq \cV_Q^+$ be an increasing family of Borel sets of finite measure. If $c_Q> 0$ we further assume either $B_t$ are all generalized sectors and that the function $t\mapsto m_{\cV_Q^+}(B_t)$ is continuous, or that $B_t=tB$ are dilations of some fixed star-shaped set $B$. Then we have for $\mu_Q$-a.e. $\G g\in Y_Q$ and for all sufficiently large $t$,
\begin{align*}
\#(\cL_Q g\cap B_t)=\omega_Qm_{\cV_Q^+}(B_t)+O_{g,\e}\left(m_{\cV_Q^+}(B_t)^{\comm{\frac{\beta_Q}{2}}+\e}\right).
\end{align*}
\end{Thm}

Next, following the strategy outlined in  \cite{KelmerYu2020} we can also apply the moment formulas to study the value distribution of random functions at the discrete set $\cL_Q$.
As explained in  \cite{KelmerYu2020}, such results follow from a second moment formula together with effective volume estimates for the functions under consideration. 
We will prove such volume estimates for both linear and quadratic forms (see \propref{p:volLinOmega} and \propref{p:volFinI}), and thus proving the desired quantitative results concerning the value distribution of these functions.
\comm{We note that although both results follow from the same strategy, the volume estimates for quadratic forms are much more involved and require extra assumptions; we thus state these results separately.}
We first state our result for a linear form, \comm{or rather several linear forms simultaneously}. 
\begin{Thm}\label{thm:linearformintro}
Assume \comm{$n\geq 2$} and $c_Q=0$. Let  $m<n$ be a positive integer and let $\{\Omega_T\}_{T>1}\subseteq \R^m$ be a decreasing family of bounded measurable sets with $\vol(\Omega_T)=T^{-a}$ for some $0<a<n-m$. Then there is some $\nu>0$ such that for almost every linear map $L:\R^{n+2}\to \R^m$ with $Q|_{\ker{L}}$ indefinite\comm{,} for all sufficiently large $T$
$$
\#\{\bm{v}\in \cL_Q: \|\bm{v}\|_Q\leq T,\; \bm{v}L\in \Omega_T\}=c(Q,L) T^{n-m-a}(1+O(T^{-\nu})),
$$
where $\|\cdot\|_Q$ is a \comm{certain} norm on $\R^{n+2}$ \comm{depending on $Q$} (see \eqref{equ:qnorm} below) and $c(Q,L)$ is some positive constant depending only on $Q$ and $L$.
\end{Thm}

\begin{remark}\label{rmk:almostevery}
Using a classification result of Sargent \cite[Lemma 2.2]{Sargent2014} (see also \lemref{l:classification} and \lemref{l:QuadClas} below), we will make the notion of ``almost every'' in \thmref{thm:linearformintro} \comm{(as well as in  \thmref{thm:valuedisquadform} below)} more precise; see \rmkref{rmk:almostall}. 
\end{remark}
By taking the sets $\Omega_T$ to be shrinking norm balls in $\R^m$ we get the following corollary extending the result of \cite[Theorem 1.3]{GhoshGorodnikNevo2020}. 

\begin{Cor}\label{c:GGN}
Assume \comm{$n\geq 2$} and $c_Q=0$. Let $m<n$ be a positive integer and  let $\kappa>\comm{\frac{m}{n-m}}$. Given any $\bm\xi\in \R^m$, for almost every linear map $L:\R^{n+2}\to \R^m$ with $Q|_{\ker{L}}$ indefinite, for all $\epsilon>0$ sufficiently small  the system of inequalities 
\begin{equation}\label{e:GGN}\|\bm{v} L-\bm\xi\|< \epsilon \quad \text{and}\quad \|\bm{v}\|< \epsilon^{-\kappa} 
\end{equation}
has a solution $\bm{v}\in \cL_Q$. 
\end{Cor}
\begin{rem} 
As explained in \cite{GhoshGorodnikNevo2020} the threshold of $\kappa>\comm{\frac{m}{n-m}}$ is optimal as can be seen from the pigeonhole principle. 
When $n=m+1$ \corref{c:GGN} follows from \cite[Theorem 1.3]{GhoshGorodnikNevo2020}, but for larger $n$ their method gives a worse threshold for $\kappa$. We should note however that the method of \cite{GhoshGorodnikNevo2020} deals with more general quadratic surfaces of the form  $\{\bm{v}\in \R^{n+2}: Q(\bm{v})=k\}$ for any rational $k\in \Q$ and not just $k=0$ and also does not have the restriction $c_Q=0$ on the quadratic form.
\end{rem}
For quadratic forms on the light cone we need to introduce one more definition. Let $F: \R^{n+2}\to \R$ be a real quadratic form of rank $m<n$, that is, it has a null subspace $V_{\rm{null}}(F)\subseteq \R^{n+2}$ (i.e. $F|_{V_{\rm{null}}(F)}=0$) of dimension $n+2-m$. We say that $F$ is \textit{indefinite on $\cV^+_Q$} if both $F$ and $Q|_{V_{\rm{null}}(F)}$ are indefinite. 
\begin{Thm}\label{thm:valuedisquadform}
Assume $n\geq 3$ and $c_Q=0$. Let $m<n$ be a positive integer and let $\{I_T\}_{T>1}\subseteq \R$ be a decreasing family of bounded measurable sets with $|I_T|=T^{-a}$ for some $0<a<n-2$. Then there is some $\nu>0$ such that for almost every rank $m$ quadratic form $F$ that is indefinite on $\cV^+_Q$ and for all sufficiently large $T$
\begin{align*}
\#\left\{\bm{v}\in \cL_Q: \|\bm{v}\|_Q\leq T,\ F(\bm{v})\in I_T\right\}=c(Q,F)T^{n-2-a}(1+O(T^{-\nu})),
\end{align*}
	where $c(Q,F)$ is some positive constant depending only on $Q$ and $F$.
	\end{Thm}
We note that the assumption $c_Q=0$ \comm{in both \thmref{thm:linearformintro} and \thmref{thm:valuedisquadform}} comes from the fact that the underlying sets within which we count are not generalized sectors. 
Finally we mention that we will prove a similar counting result for a certain more general family of homogeneous functions and \thmref{thm:valuedisquadform} will be a special case of this  result when the degree is two; see \thmref{thm:homogen}.

%
\subsection{Application to intrinsic Diophantine approximations on ellipsoids}\label{sec:indioapp}
 In \cite{KleinbockMerrill2015}, Kleinbock and Merril used  dynamical methods to study problems in intrinsic Diophantine approximations on the sphere, that is, the question of how well can one approximate a point on the sphere by rational points on the sphere of a bounded height (where the height of a rational point 
 $\frac{\bm{p}}{q}$ 
is the denominator $q$ when written in lowest terms). Their methods and results were then generalized by Fishman, Kleinbock, Merrill, and Simmons \cite{FishmanKleinbockMerrillSimmons2022} to deal with more general nonsingular rational quadratic hypersurfaces. In what follows we give a refined and quantitative version of many of their results. 

Before stating our results, let us first fix some notation. For the remaining of this section we consider rational $\Q$-isotropic quadratic forms of the form 
\begin{align}\label{equ:qdfspfo}
Q(\bm{x},y)=\cQ(\bm{x})-y^2
\end{align} 
with $\cQ$ a positive definite rational quadratic form in $n+1$ variables. Let
\begin{align}\label{equ:ellipsoid}
S_{\cQ}:=\left\{\bm{x}\in \R^{n+1}: \cQ(\bm{x})=1\right\}
\end{align}
be the ellipsoid associated to $\cQ$. \comm{For instance when $Q=Q_n$ is the standard form as in \eqref{equ:stformqn}, 
$
S_{\cQ}=S^n:=\left\{\bm{x}\in \R^{n+1}: \|\bm{x}\|=1\right\}
$ 
is the unit $n$-sphere in $\R^{n+1}$. On $S^n$ there is a unique rotation-invariant probability measure which we denote by $\sigma_n$. In general, $S_{\cQ}=S^n\tilde{\tau}^{-1}$ with $\tilde{\tau}\in \GL_{n+1}(\R)$ such that $\cQ(\bm{x})=\|\bm{x}\tilde{\tau}\|^2$. We denote by $\sigma_{\cQ}$ the probability measure on $S_{\cQ}$ obtained as the pushforward of $\sigma_n$ under the right $\tilde{\tau}^{-1}$-multiplication map.} We also let $\|\cdot\|_{\cQ}$ be the norm on $\R^{n+1}$ given by $\|\bm{x}\|_{\cQ}:=\sqrt{\cQ(\bm{x})}$. 
For any $\bm{x}\in S_{\cQ}$ and $r>0$ let
\begin{align}\label{equ:opendiscell}
\mathfrak{D}^{\cQ}_r(\bm{x}):=\left\{\bm{x}'\in S_{\cQ}: \|\bm{x}'-\bm{x}\|_{\cQ}<r\right\}
\end{align} 
be the open ball in $S_{\cQ}$ of radius $r$ with center $\bm{x}$, and for any $T>1$ consider the following counting function:
\begin{align}\label{equ:countingfuncell}
\cN_{\cQ}(\bm{x}, r, T):=\#\left\{\tfrac{\bm{p}}{q}\in \mathfrak{D}^{\cQ}_r(\bm{x}): \gcd(\bm{p},q)=1,\  1\leq q< T\right\}.
\end{align}
\comm{Here $\gcd(\bm{p},q)$ is the greatest common divisor of the entries of the vector $\bm{p}\in \Z^{n+1}$ together with the integer $q\in \N$.}
Note that for $r\geq 2$ this counting function just counts the number of rational points in $S_{\cQ}$ with  denominator bounded by $T$; \comm{denote this counting function by $\cN_{\cQ}(T)$}. It was shown by Duke in \cite{Duke2003} that 
$\cN_{\cQ}(T)\sim \kappa_{\cQ}T^{n}$ for some $\kappa_{\cQ}>0$ and moreover that the rational points are equidistributed on $S_{\cQ}$ in the sense that 
$\cN_{\cQ}(\bm{x},r,T)$ is asymptotic to  $\kappa_{\cQ} T^n\sigma_{\cQ}(\mathfrak{D}^{\cQ}_r(\bm{x}))$. While his results were stated for \comm{the special case of $S_{\cQ}=S^2$,} his methods extend to higher dimensions \comm{and general ellipsoids}. 

Going beyond equidistribution, an interesting question (related to intrinsic Diophantine approximation) is how small can one take $r$ in terms of $T$ and still have this asymptotic formula. A natural threshold is suggested by the results of \comm{Fishman, Kleinbock, Merrill, and Simmons \cite[Theorem 5.1]{FishmanKleinbockMerrillSimmons2022} (and earlier by Kleinbock and Merrill \cite[Theorem 4.1]{KleinbockMerrill2015} for the special case of spheres)}, who proved an analogue of Dirichlet's theorem,  showing that there is a constant $C>0$ such that for any $\bm{x}\in S_{\cQ}$, there always exists a rational point $\frac{\bm{p}}{q}\in S_{\cQ}$ satisfying  
	$$
	\left\|\bm{x}-\frac{\bm{p}}{q}\right\|_{\cQ}<\frac{C}{T^{\frac12}q^{\frac12}}\quad \text{and}\quad 1\leq q\leq T,
	$$
implying in particular that $\cN_{\cQ}(\bm{x}, r, T)> 0$ for $r\geq CT^{-1/2}$. We now give the following quantitative result for this counting function.

\begin{Thm}\label{thm:c1digen}
Let $Q(\bm{x},y)=\cQ(\bm{x})-y^2$ be as in \eqref{equ:qdfspfo}. Let $\kappa_{\cQ}=\frac{\omega_Q}{n}$ and $\beta=\beta_Q$. Then there exists a constant \comm{$C_{\cQ}>0$} such that 
for any $r,T>0$ satisfying 
\comm{$C_{\cQ}T^{-\frac{2-\beta}{3-\beta}}\leq r< 1$} and for any $\bm{x}\in S_{\cQ}$, 
\begin{align}\label{equ:aysformell}
\cN_{\cQ}(\bm{x}, r, T)=\kappa_{\cQ}T^{n}\sigma_{\cQ}(\mathfrak{D}^{\cQ}_r(\bm{x}))\left(1+O_{\cQ}\left(r^{-\frac{(3-\beta)n}{2n+3}}T^{-\frac{(2-\beta)n}{2n+3}}\right)\right).
\end{align}
\end{Thm}

\begin{remark}\label{rmk:dircount}
Fixing $r$ and letting $T\to\infty$ the above counting formula shows that 
\begin{align*}
\cN_{\cQ}(\bm{x}, r, T)=\kappa_{\cQ}T^{n}\sigma_{\cQ}(\mathfrak{D}^{\cQ}_r(\bm{x}))\left(1+O_{\cQ,r}\left(T^{-\frac{(2-\beta)n}{2n+3}}\right)\right),
\end{align*}
 giving an effective equidistribution result with power saving estimates for the remainder. Moreover, \comm{as part of the theorem we have $\kappa_{\cQ}=\frac{\omega_Q}{n}$. When $Q=Q_n$ as in \eqref{equ:stformqn}, i.e. when $S_{\cQ}=S^n$}, $\omega_{Q_n}$ was  explicitly computed in terms of special values of zeta and $L$-functions in \cite[Equation (1.7)]{KelmerYu2022a}, we thus get an explicit expression for $\kappa_{\cQ}$ in this case, \comm{generalizing a formula of Duke \cite[Theorem 1]{Duke2003} for $\kappa_{\cQ}$ when $n=2$}.
\end{remark}

\begin{remark}
\comm{
When $c_Q=0$ (i.e. when $\beta_Q=1$) we can get the asymptotic estimate in the full range of $C_{\cQ}T^{-1/2}\leq r<1$ (where the results of 
\cite{FishmanKleinbockMerrillSimmons2022} imply that $\cN_{\cQ}(\bm{x}, r, T)> 0$). We expect the same estimate should hold in general, however, when $c_Q>0$ we can only establish this asymptotic estimate in a slightly smaller range of  $C_{\cQ}T^{-\frac{n-2}{2n-2}}\leq r<1$ when $n$ is even and $C_{\cQ}T^{-\frac{n-1}{2n-1}}\leq r<1$ when $n$ is odd, which is an artifact of the secondary term in the moment formula. 
}
\end{remark}

Since the asymptotic results of \thmref{thm:c1digen} apply also to rational $\bm{x}\in S_{\cQ}$ the exponent of $T^{-1/2}$ can not be improved, (e.g. \comm{for $S_{\cQ}=S^n$ and} $\bm{x}=(\bm{0},1)$ we have $\cN_{\cQ}(\bm{x}, r, T)=1$ for all $r\leq T^{-1/2}$). 
For generic $\bm{x}\in S_{\cQ}$, on the other hand,  it is expected that one could do even better. Indeed, \comm{for $S_{\cQ}=S^n$,} it was shown by Ghosh, Gorodnik and Nevo \cite{GhoshGorodnikNevo2013} for odd $n\geq 5$ for almost all $\bm{x}\in S_{\cQ}$  that  $\cN_{\cQ}(\bm{x}, T^{-\lambda},T)>0$ for $\lambda\in (0, \frac12+\frac{3}{2n})$ (and similarly for $\lambda\in (0,\frac23)$ when $n=3$). 
For this problem we can improve the above range to the full range of $\lambda\in (0,1)$ while still get an asymptotic counting formula with an explicit power saving error term.
When $r=\frac{1}{T}$ the main term in the asymptotic formula no longer grows with $T$ so $\lambda<1$ is a natural limit for such asymptotic estimates.
Nevertheless, our method also works when taking $r= \frac{(\log T)^{\lambda}}{T}$ as long as the exponent of $\log(T)$ is sufficiently large. Explicitly we show the following.

\begin{Thm}\label{thm:countforgenericgen}
Keep the notation and assumptions as in \thmref{thm:c1digen} 
and let $\varkappa_{\comm{\cQ}}=\tfrac{\G\left(\frac{n+3}{2}\right)\omega_Q}{\pi^{\frac12}n(n+1)\G\left(\frac{n+2}{2}\right)}$. 
\begin{enumerate}
\item
Let $r_T=T^{-\lambda}$ for some $0<\lambda< 1$. Then for $\sigma_{\cQ}$-a.e. $\bm{x}\in S_{\cQ}$ and for all $T$ sufficiently large 
\begin{align*}
\cN_{\cQ}(\bm{x},r_T,T)=\varkappa_{\cQ}T^{n-n\lambda}+O_{Q}\left(T^{(n-n\lambda)(1-\frac{2-\beta}{n+4})}\log (T)+T^{n-n\lambda-2\lambda}\right).
\end{align*}
\item 
Let $r_T=\tfrac{(\log T)^{\lambda}}{T}$ for some $\lambda>\frac{1}{n(2-\beta)}$. Then for $\sigma_{\cQ}$-a.e. $\bm{x}\in S_{\cQ}$ and for all $T$ sufficiently large 
\begin{align*}
\cN_{\cQ}\left(\bm{x},r_T,T\right)=\varkappa_{\cQ}\left(\log T\right)^{n\lambda}+O_{Q,\lambda}\left(\left(\log T\right)^{n\lambda -\frac{n\lambda(2-\beta)-1}{n+4}}\log\log T\right).
\end{align*}
\end{enumerate} 
\end{Thm}

The results of \thmref{thm:countforgenericgen} imply, in particular, that for any $0<\lambda<1$ for  $\sigma_{\cQ}$-a.e. $\bm{x}\in S_{\cQ}$ the \comm{inequality} $\|\bm{x}-\frac{\bm{p}}{q}\|_{\cQ}\leq q^{-\lambda}$ has infinitely many solutions in $\frac{\bm p}{q}\in S_{\cQ}\cap \Q^{n+1}$. One can ask the same question for more general approximating functions $\psi:\N\to (0,\infty)$. Borrowing the terminology of \cite{KleinbockMerrill2015}, we say that a point $\bm{x}\in S_{\cQ}$ is \textit{$\psi$-approximable} in $S_{\cQ}$, if there are infinitely many rational points $\frac{\bm p}{q}\in S_{\cQ}\cap \Q^{n+1}$ satisfying $\|\bm{x}-\frac{\bm p}{q}\|_{\cQ}< \psi(q)$.
For any decreasing function $\psi$, if the series $\sum_{q\in\N}q^{n-1}\psi(q)^n$ is convergent, it follows from the easy half of the Borel-Cantelli lemma that $\bm{x}$ is not $\psi$-approximable in $S_{\cQ}$ for  $\sigma_{\cQ}$-a.e. $\bm{x}\in S_{\cQ}$ (see \cite[p. 297]{KleinbockMerrill2015}).
\comm{Using a dynamical argument it was proved in  
\cite[Theorem 6.3]{FishmanKleinbockMerrillSimmons2022}  (and \cite[Theorem 1.3]{KleinbockMerrill2015} for the special case of the sphere) that the converse is also true, under the additional assumption that the function $q\mapsto q\psi(q)$ is decreasing.}

Their result is analogous to Khintchine's theorem on Diophantine approximation in Euclidean spaces which has a quantitative counterpart proved by Schmidt \cite{Schmidt1960a}. We now give a corresponding quantitative refinement of the above result, providing an asymptotic estimate, analogous to the quantitative result of Schmidt. To this end 
we define the following counting function that for $\bm{x}\in S_{\cQ}$ and $T>1$,
\begin{align}\label{equ:khincountfunell}
\cN_{\cQ, \psi}(\bm{x}, T):=\#\left\{\tfrac{\bm{p}}{q}\in S_{\cQ}: \gcd(\bm{p},q)=1,\  \|\bm{x}-\tfrac{\bm{p}}{q}\|_{\cQ}<\psi(q),\ 1\leq q< T\right\}.
\end{align}
Note that the function $\cN_{\cQ,\psi}(\bm{x}, T)$ is increasing in $T$ and that it is unbounded if and only if  $\bm{x}$ is $\psi$-approximable in $S_{\cQ}$.  We also denote by
\begin{align}\label{equ:partialsum}
 J_{\psi}(T):=\sum_{1\leq q< T}q^{n-1}\psi(q)^n\quad \mbox{and }\quad I_{\psi}(T):=\sum_{1\leq q< T}q^{n-1}\psi(q)^{n+2}.
 \end{align}
\begin{Thm}\label{thm:khinttypegen}
Let $Q(\bm{x},y)=\cQ(\bm{x})-y^2$ be as in \thmref{thm:c1digen} and assume $c_Q=0$. For any decreasing $\psi: \N\to (0,\infty)$ satisfying that $\lim\limits_{q\to\infty}\psi(q)=0$ and $\sum_{q\in\N}q^{n-1}\psi(q)^n=\infty$,  for $\sigma_{\cQ}$-a.e. $\bm{x}\in S_{\cQ}$ and for all sufficiently large $T$,
\begin{align}\label{equ:khintchine}
\cN_{\cQ,\psi}(\bm{x}, T)= n\varkappa_{\cQ}J_{\psi}(T)+O_{\cQ, \psi}\left(J_{\psi}(T)^{\frac{n+3}{n+4}}\log(J_{\psi}(T))+I_{\psi}(T)\right),
\end{align}
where $\varkappa_{\cQ}$ is as in \thmref{thm:countforgenericgen}. 
\end{Thm}

If the series $\sum_{ q=1}^\infty q^{n-1}\psi(q)^{n+2}$ converges (and in particular if $q\psi(q)$ is decreasing), then $I_{\psi}(T)$ is bounded and we get a power saving estimate for the remainder.
For a general decreasing function $\psi$, it is not hard to see that the term $I_{\psi}(T)$ is dominated by the main term so that $\cN_{\cQ,\psi}(\bm{x}, T)\sim n\varkappa_{Q} J_{\psi}(T)$ as $T\to \infty$. In particular, we get the following corollary recovering the result of \cite[Theorem 6.3]{FishmanKleinbockMerrillSimmons2022} 
without the additional assumption that $q\psi(q)$ is decreasing.

\begin{Cor}
Keep the assumptions as in \thmref{thm:khinttypegen}. For any decreasing function $\psi:\N\to (0,\infty)$ we have that $\sigma_{\cQ}$-almost every (respectively almost no) $\bm{x}\in S_{\cQ}$ is $\psi$-approximable if and only if the series $\sum_{q\in\N}q^{n-1}\psi(q)^n$ diverges (respectively converges).
\end{Cor}

It is reasonable to expect that the same result holds also when $c_Q>0$.  
We note that for the special case $\psi(q)=c/q$ with $c>0$, one can get around this restriction using dynamics as illustrated by Alam and Ghosh in \cite{AlamGhosh2022}, who proved a similar counting formula for the sphere in all dimensions for this explicit function, and more recent work of Ouaggag \cite{Ouaggag2023} improving on \cite{AlamGhosh2022} by providing a power-saving estimate on the error term.

\subsection{Structure of the paper}
For the readers convenience we give a short outline of the paper describing the general ideas going into the proofs. In section \ref{sec:setrans} we define the light cone Siegel transform,  establish the moment formula, and then use it to prove the mean square bound for the discrepancy.
In section \ref{sec:meantocount} we prove several new results showing how one can obtain counting estimates from a mean square bound on the discrepancy.
 We expect these ideas to be useful in a larger context so we state and prove them in greater generality. While for a family of increasing sets, the classical argument of Schmidt gives very good estimates that hold for almost all translates, we give new arguments that can also handle families of sets that are not increasing, as well as estimates that hold for almost all translates from a subgroup or even for all translates; see \thmref{thm:cfamd}, \thmref{thm:countforall} and \thmref{thm:countnoninc}. For the latter we need to assume our sets satisfy suitable well roundedness properties. With these general arguments at hand, our main results are then proved by translating the problem to a counting problem in families of sets and then establishing the corresponding well roundedness properties in each case. Explicitly, in section \ref{sec:ranform}, we use the classification of linear and quadratic forms on the light cone to translate the problem to a suitable counting problem in a family of sets and establish the corresponding volume estimates for these sets thus proving \thmref{thm:linearformintro} and \thmref{thm:valuedisquadform}. In section \ref{sec:intrdioapp} we translate the problem of counting rational points on the sphere to counting primitive points on certain families of growing sets on the light cone, we then establish corresponding well-roundedness properties and volume estimates for these sets from which we can deduce Theorems \ref{thm:c1digen}--\ref{thm:khinttypegen}.

\subsection{Notation and conventions}
Throughout the paper, the norm notation $\|\cdot\|$ is reserved exclusively for the Euclidean norm in various Euclidean spaces, \comm{and all functions are assumed to be measurable}. \comm{For any linear transformation $g$ acting on an Euclidean space, we denote by $\|g\|_{\rm op}$ its operator norm with respect to $\|\cdot\|$ on this Euclidean space.}  Let $I\subseteq \R$ be an interval (not necessarily bounded). A function $f : I\to \R$ is called \textit{increasing} (respectively\ \textit{decreasing}) if $f(t_1)\leq f(t_2)$ (respectively\ $f(t_1)\geq f(t_2)$) whenever $t_1<t_2$.  Similarly, we say a family of Borel sets $\{B_t\}_{t\in I}$ in some algebraic variety is \textit{increasing} (respectively \textit{decreasing}) if $B_{t_1}\subseteq B_{t_2}$ (respectively $B_{t_1}\supseteq B_{t_2}$) whenever $t_1<t_2$. 
The notation, $\bm{v}\in\R^n$, is always used to denote column vectors. For two positive quantities $A$ and $B$, we will use the notation $A\ll B$  {or $A=O(B)$} to mean that there is a constant $c>0$ such that $A\leq cB$, and we will use subscripts to indicate the dependence of the constant on parameters. We will write $A\asymp B$ for $A\ll B\ll A$. For any two positive functions $f, g$ defined on $\R$ we write $f(t)\sim g(t)$ to mean that $\lim\limits_{t\to\infty}\frac{f(t)}{g(t)}=1$.

\subsection*{Acknowledgements}
\comm{We would like to thank the anonymous referee for various suggestions that improved the exposition of this paper.}

\section{Light-cone Siegel transform and its moment formulas}\label{sec:setrans}
Let $n\in\N$ be a positive integer and let \comm{$Q: \R^{n+2}\to \R$} be a rational $\Q$-isotropic  quadratic form of signature $(n+1,1)$. \comm{Here recall that $Q$ is $\Q$-isotropic means that there exists a nonzero $\bm{v}\in \Q^{n+2}$ such that $Q(\bm{v})=0$.} Let $G=\SO_Q^+(\R)$ be the identity component of the special orthogonal group preserving $Q$ and let $\G=\SO_Q^+(\Z):=G\cap \SL_{n+2}(\Z)$ be the lattice of integer points. Let 
$$
\cV_Q:=\left\{\bm{v}\in \R^{n+2}\setminus \{\bm{0}\} : Q(\bm{v})=0\right\}
$$
be the light cone of $Q$ and let \comm{$\cV_Q^+\subset \cV_Q$ be a connected component of $\cV_Q$ (which will be chosen more explicitly in section \ref{sec:cormea} below).} \comm{Let 
$$
\Z_{\rm pr}^{n+2}:=\left\{\bm{v}\in \Z^{n+2}: \gcd(\bm{v})=1\right\}
$$
be the set of primitive integer points in $\R^{n+2}$
and denote by $\cL_Q:=\cV_Q^+\cap \Z_{\rm pr}^{n+2}$ the subset of primitive integer points in $\cV_Q^+$.} Since $Q$ is $\Q$-isotropic, $\G$ is a non-uniform lattice of $G$ and $\cL_Q$ is an infinite discrete subset of $\cV_Q^+$. Note that the right action of $G$ on $\cV_Q^+$ is transitive, and \comm{since $\G<\SL_{n+2}(\Z)$, the induced action of $\G$ on $\cV_Q^+$} preserves $\cL_Q$, i.e. $\cL_Q\gamma= \cL_Q$ for all $\gamma\in \G$. As described in the introduction the homogeneous space $Y_Q=\G\bk G$ parameterizes the space of $\cL_Q$-translates via $\G g\leftrightarrow \cL_Q g$.
Generalizing the classical Siegel transform on the space of lattices, for any bounded and compactly supported function $f: \cV^+_Q\to \C$ we define its \textit{light-cone Siegel transform} by 
\begin{align}\label{def:lcstrans}
S_Qf(g):=\sum_{\bm{v}\in \cL_Q}f(\bm{v}g),\quad (g\in G).
\end{align} 
Since $f$ is bounded and compactly supported, for any $g\in G$ the series defining $S_Qf$ is a finite sum and hence absolutely converges. Moreover, $S_Qf$ is left $\G$-invariant, and can be viewed as a function on $Y_Q$. The main goal of this section is to prove the moment formulas of the light-cone Siegel transform stated in \thmref{thm:mainrest}.

\subsection{Preliminaries}
We first briefly introduce some necessary backgrounds and notation and we refer the reader to \cite[Section 2]{KelmerYu2022a} for more details.

\subsubsection{Coordinates and measures}\label{sec:cormea}
The main reference of this subsection is \cite[Section 2.2]{KelmerYu2022a}. 
We first set up the notation for the standard $(n+1,1)$ form $Q_n$ given in \eqref{equ:stformqn}
and then for a general form via a translating matrix.

First for $Q=Q_n$, take
\begin{align}\label{equ:lightconestfo}
\cV_{Q_n}^+=\left\{\bm{v}\in \R^{n+2}: Q(\bm{v})=0,\ v_{n+2}>0\right\},
\end{align} 
and fix $\bm{e}_0=(-1,\bm{0},1)\in \cV_{Q_n}^+$. The group \comm{$\SO_{Q_n}^+(\R)$} has an Iwasawa decomposition $\SO_{Q_n}^+(\R)=\comm{U_{Q_n}A_{Q_n}K_{Q_n}}$ with 
\begin{align*}
\comm{U_{Q_n}}:=\left\{u_{\bm{x}}=\left(\begin{smallmatrix} 1-\frac{\|\bm{x}\|^2}{2} & \bm{x} & \frac{\|\bm{x}\|^2}{2}\\ 
-\bm{x}^t & I_n & \bm{x}^t\\ 
-\frac{\|\bm{x}\|^2}{2} & \bm{x} & 1+\frac{\|\bm{x}\|^2}{2}\end{smallmatrix}\right): \bm{x}\in \R^n\right\}
\end{align*}
the maximal unipotent subgroup fixing $\bm{e}_0$,
\begin{align*}
\comm{A_{Q_n}}:=\left\{a_y=\left(\begin{smallmatrix}
\tfrac{y+y^{-1}}{2} &  & \tfrac{y-y^{-1}}{2}\\
 & I_n & \\
\tfrac{y-y^{-1}}{2} &  &\tfrac{y+y^{-1}}{2} \end{smallmatrix}\right): y>0\right\}
\end{align*}
the $\R$-split torus with $a_y$ acting on $\bm{e}_0$ as $\bm{e}_0a_y=y^{-1}\bm{e}_0$ and 
\begin{align}\label{equ:macomptgpk}
\comm{K_{Q_n}}:=\left\{k=\begin{pmatrix} \tilde{k} &  \\  &1\end{pmatrix}: \tilde{k}\in \SO_{n+1}(\R)\right\}
\end{align}
a maximal compact subgroup. Let $L$ be the stabilizer of $\bm{e}_0$ and let $P$ be the parabolic subgroup fixing the line spanned by $\bm{e}_0$. More precisely, $P=UAM$ and $L=UM$ with
$$
\comm{M_{Q_n}}=\left\{m=\left(\begin{smallmatrix} 1 &  &  \\  &\tilde{m}& \\  & &1\end{smallmatrix}\right): \tilde{m}\in \SO_n(\R)\right\}
$$ 
the centralizer of $A$ in $K$. 

Now for a general $Q$, there exists some $\tau\in \GL_{n+2}(\R)$ such that $Q(\bm{v})=Q_n(\bm{v}\tau)$. We can then take the one-sheeted light-cone $\cV_Q^+=\cV_{Q_n}^+\tau^{-1}$ and the base point $\bm{e}_0=(-1,\bm0,1)\tau^{-1}$. Note also that $G=\tau\SO_{Q_n}^+(\R)\tau^{-1}$ and we can take the subgroups of $G$ to be the $\tau$-conjugation of the corresponding subgroups of $\SO_{Q_n}^+(\R)$, that is, $G=U_QA_QK_Q$ with $U_Q=\tau U_{Q_n}\tau^{-1}$, $A_Q=\tau A_{Q_n}\tau^{-1}$ and $K_Q=\tau K_{Q_n}\tau^{-1}$ with $U_{Q_n},A_{Q_n}, K_{Q_n}$ as above. When the form $Q$ is fixed and there is no risk of ambiguity we will omit these subscripts.

\begin{remark}\label{rmk:formoftau}
Note that we can take the matrix $\tau$ above to be of the form $\tau=ka$ with $k\in \text{O}_{n+2}(\R)$ and $a$ a diagonal matrix. Indeed, let $J\in \GL_{n+2}(\R)$ be the symmetric matrix such that $Q(\bm{v})=\bm{v}J\bm{v}^t$. Since $Q$ is of signature $(n+1,1)$, $J$ has $n+1$ positive eigenvalues $t_1\geq t_2\geq \ldots \geq t_{n+1}$ and $1$ negative eigenvalue $t_{n+2}$. Moreover, since $J$ is real symmetric, it can be diagonalized by an orthogonal matrix, that is, there is some $k\in \text{O}_{n+2}(\R)$ such that $k^tJk=\diag(t_1,\ldots, t_{n+2})$. Now take $a=\diag(|t_1|^{\frac12},\ldots, |t_{n+2}|^{\frac12})$ so that $a^{-1}k^tJka^{-1}=\comm{\left(\begin{smallmatrix}
I_{n+1} & \\
 & -1\end{smallmatrix}\right)} 
 $, 
 or equivalently, $J=ka\comm{\left(\begin{smallmatrix}
I_{n+1} & \\
 & -1\end{smallmatrix}\right)}a^tk^t$. This then implies that $Q(\bm{v})=Q_n(\bm{v}\tau)$. Throughout this paper, we will take $\tau=ka$ with $k$ and $a$ as above. Clearly, $a$ is uniquely determined by $J$, but $k$ is not. Indeed we may replace $k$ by any $k'\in \text{O}_{n+2}(\R)$ such that $k^{-1}k'$ commutes with $\diag(t_1,\ldots, t_{n+2})$.
\end{remark}
Any $g\in G$ can be written as $g=u_{\bm{x}}a_yk$ with $k\in K$ and in these coordinates the Haar measure of $G$ is given (up to scaling) by
\begin{align}\label{equ:Haar1}
\d\mu_G(g)=y^{-(n+1)}\,\d\bm{x}\d y\d\mu_K(k),
\end{align}
where $\d\bm{x}$ is the usual Lebesgue measure on $\R^n$ and $\mu_K$ is the probability Haar measure of $K$. 

The subgroup $L$ is unimodular with its Haar measure given by 
\begin{align}\label{equ:Lhaar}
\d\mu_L(u_{\bm{x}}m)=\d\bm{x}\d\mu_M(m),
\end{align}
where $\mu_M$ is the probability Haar measure of $M\cong \SO_n(\R)$. 
Since $L$ is the stabilizer of $\bm{e}_0$ and $G$ acts transitively on $\cV^+_Q$, we can identify $\cV_Q^+$ with the homogeneous space $L\bk G$, which gives a natural 
right $G$-invariant measure on $\cV^+_Q$. Explicitly, further identifying $L\bk G$ with $A\times M\bk K$ gives 
 natural polar coordinates on $\cV^+_Q$: Every $\bm{v}\in \cV^+_Q$ can be written uniquely as $\bm{v}=\bm{e}_0a_yk$ for some $y>0$ and $k\in M\bk K$. In these coordinates the measure 
\begin{align}\label{equ:lebmea}
\d m_{\cV^+_Q}(\bm{e}_0a_yk):=y^{-(n+1)}\,\d y\d\mu_{M\bk K}(k)
\end{align}
is such an invariant measure. Here $\mu_{M\bk K}$ is the unique right $K$-invariant probability measure on the homogeneous space $M\bk K$ which is homeomorphic to the unit sphere $S^n$.
\begin{remark}\label{rmk:spherehomeo}
When $Q=Q_n$, this homeomorphism can be made explicit: The projection map sending $k\in K$ to $\bm{\alpha}\in S^n$ with $\bm{\alpha}$ such that $\bm{e}_0k=(\bm{\alpha},1)$ induces such a homeomorphism. We note that the pushforward of the invariant probability measure $\mu_{M\bk K}$ under this homeomorphism is the rotation-invariant probability measure $\sigma_n$. We also fix a continuous section 
\begin{align}\label{equ:section}
\mathfrak{s}: S^n\to K
\end{align}
of this projection map from $K$ to $S^n$, and for any $\bm{\alpha}\in S^n$, we denote by $k_{\bm{\alpha}}=\mathfrak{s}(\bm{\alpha})$. In particular, \eqref{equ:lebmea} can be rewritten as 
\begin{align}\label{equ:lebmeasper}
\d m_{\cV^+_{Q_n}}(\bm{e}_0a_yk_{\bm{\alpha}})=y^{-(n+1)}\,\d y\d\sigma_n(\bm{\alpha}).
\end{align}
\end{remark}

\subsubsection{Eisenstein series}\label{sec:eisenseries}
In this subsection we recall necessary backgrounds on Eisenstein series. The main reference is \cite[Sections 2.3 and 2.4]{KelmerYu2022a}.
Let $P$ be the stabilizer of the line spanned by $\bm{e}_0$. It is a parabolic subgroup of $G$ and has a Langlands decomposition $P=UAM$ with $U, A, M$ as above. In general, as discussed in \cite[Section 2.3]{KelmerYu2022a} every parabolic subgroup of $G$ is of the form $P'=kPk^{-1}$ for some $k\in K$ and its unipotent radical is given by $kUk^{-1}$. 
The \textit{cusps} of $\G$ are the $\G$-conjugacy classes of parabolic subgroups of $G$ whose unipotent radicals intersect $\G$ nontrivially. Let $\kappa$ be the number of cusps of $\G$ and let $P_1,\ldots, P_{\kappa}$ be a full set of representatives of these classes. Note that for each $P_i$, we can take $k_i\in K$ such that $P_i=k_iPk_i^{-1}$, and thus it has a Langlands decomposition $P_i=U_iA_iM_i$ with $U_i=k_iUk_i^{-1}$, $A_i=k_iAk_i^{-1}$ and $M_i=k_iMk_i^{-1}$. We denote by $\G_{U_i}:=\G\cap U_i$ and \comm{$\G_{P_i}:=\G\cap P_i$}. By definition, $\G_{U_i}$ is non-trivial, and \comm{indeed} it is a finite-index subgroup of $\G_{P_i}$; see \cite[Lemma 2.1]{KelmerYu2022a}.

For each $P_i$, we fix the \textit{scaling matrix} $\tau_i=k_ia_{y_i}$ with $k_i\in K$ as above and $y_i>0$ the unique number such that $\mu_L(\tau_i^{-1}\G_{P_i}\tau_i\bk L)=1$. 
The \textit{(spherical) Eisenstein series} corresponding to the $i$-th cusp is then defined 
for $\Re(s)>n$ and $g\in G$ by the convergent series
\begin{equation}\label{e:EisensteinSph}
	E_i(s, g)=\sum_{\gamma\in \Gamma_{\!P_i}\bk \Gamma} y(\tau_i^{-1}\gamma g)^s,
\end{equation}
where $y(g)$ is \comm{the $y$-parameter in the Iwasawa decomposition of $g$, that is,} $g=u_{\bm{x}}a_{y(g)}k$. For each $1\leq  j\leq \kappa$ the \comm{Eisenstein series $E_i(s, \tau_ju_{\bm x} g)$ (as a function in $\bm{x}\in \R^n$) is invariant under translation by $\tau_j^{-1}\G_{U_j}\tau_j$; it thus has a Fourier expansion with respect to $\tau_j^{-1}\G_{U_j}\tau_j\bk U$ and the \textit{constant term} of this Fourier expansion is given by} 
\begin{equation}\label{e:ConstSph} 
	c_{ij}(s,g):=\frac{1}{\vol(\tau_j^{-1}\G_{U_j}\tau_j\bk U)}\int_{\tau_j^{-1}\G_{U_j}\tau_j\bk U}E_i(s, \tau_ju_{\bm x} g)\,\d\bm{x},
\end{equation}
which is known to be of the form
\begin{equation}\label{e:ConstSph1}
c_{ij}(s,g)=\delta_{ij}y(g)^s+\varphi_{ij}(s)y(g)^{n-s}
\end{equation}
for some holomorphic function $\varphi_{ij}$ defined for $\Re(s)>n$; \comm{see e.g. \cite[Chapter 6, (1.11)]{CohenSarnak1980}}.

The series $E_i(s, g)$ (and hence also $\varphi_{ij}$) has a meromorphic continuation to the whole $s$-plane, which on the half plane $\Re(s)\geq \frac{n}{2}$ is holomorphic except for a simple pole at $s=n$ (called the \textit{trivial pole}) and possibly finitely many simple poles on the interval $(\frac{n}{2}, n)$ (called \textit{exceptional poles}). We denote by $\cC_{\G}\subseteq (\frac{n}{2}, n)$ the finite set of exceptional poles of all Eisenstein series of $\G$. 

The residue of $E_i(s,g)$ at $s=n$ is a constant which is the same for Eisenstein series at all cusps, given by the reciprocal of the measure of the homogeneous space $Y_Q$, that is, for each $1\leq i\leq \kappa$ and $g\in G$,
\begin{align}\label{equ:rescovol}
\omega_{\G}:=\Res_{s=n}E_i(s,g)=\mu_G(Y_Q)^{-1},
\end{align}
(see e.g. \cite[Proposition 2.4]{KelmerYu2022a}). In view of this relation, the unique $G$-invariant probability measure on $Y_Q$, is given by 
\begin{align}\label{equ:haar3}
\d\mu_{Q}(g)=\omega_{\G}\,y^{-(n+1)}\,\d\bm{x}\d y\d \mu_M(m)\d\mu_{M\bk K}(k),
\end{align} 
where $g=u_{\bm{x}}a_ymk$ with $u_{\bm{x}}\in U$, $a_y\in A$, $m\in M$ and $k\in M\bk K$.

These Eisenstein series  satisfy a functional equation relating $s$ and $n-s$, and a consequence of this functional equation is that on the line $\Re(s)=\frac{n}{2}$ the functions $\varphi_{ij}$ ($1\leq i,j\leq \kappa$) are holomorphic and satisfy the bound that 
\begin{align}\label{equ:funcons1}
	|\varphi_{ij}(s)|\leq 1, \quad \forall\ \Re(s)=\tfrac{n}{2}.
	\end{align}
Moreover, all poles of the Eisenstein series come from the poles of the scattering matrix $\varphi_{ij}(s)$ and for any $1\leq i,j\leq \kappa$,
\begin{align}\label{equ:constomeg}
\Res_{s=n}\varphi_{ij}(s)=\Res_{s=n}E_i(s,g)=\omega_{\G},
\end{align}
and for any of the exceptional poles $\sigma \in \cC_\G$, it follows from the Maass-Selberg relations (given in \cite[Equation (7.44)]{Selberg1989} and
\cite[Chapter 6, (1.62)]{CohenSarnak1980}) that 
\begin{equation}\label{eq:scatteringres}
\Res_{s=\sigma}\varphi_{ij}(s)=\left\langle \Res_{s=\sigma}E_i(s,\cdot),\Res_{s=\sigma}E_j(s,\cdot)\right\rangle,
\end{equation}
where the inner product is with respect to $\mu_Q$.
\subsection{Relations to incomplete Eisenstein series}
Recall that $\G$ acts on the discrete set $\cL_Q$. This action is in general not transitive and the orbits of this action are in one-to-one correspondence with cusps of $\G$ (see \cite[Lemma 2.2]{KelmerYu2022a}). This correspondence implies that the light-cone Siegel transform is a sum of incomplete Eisenstein series  as follows. 

Let $P_1,\ldots, P_{\kappa}$ be a full set of \comm{cusp} representatives as before. For any bounded and compactly supported function $f: \cV^+_Q\to \C$, the \textit{incomplete Eisenstein series attached to $f$ at $P_i$} is defined by 
\begin{align*}
E_i(g\,|\, f):=\sum_{\gamma\in \G_{P_i}\bk \G}f(\bm{e}_0\tau_i^{-1}\gamma g), \quad (g\in G)
\end{align*}
where $\G_{P_i}=\G\cap P_i$ and $\tau_i$ is the $i$-th scaling matrix as before.
Since $f$ is bounded and compactly supported, the defining series 
for $E_i(g\,|\, f)$ is absolutely convergent and since $E_i(\comm{g}\,|\, f)$ is left 
$\G$-invariant, it can be viewed as a function on the homogeneous space $Y_Q$. We have the following relation. 

\begin{Lem}[{\cite[Lemma 3.1]{KelmerYu2022a}}]\label{lem:sumincomein}
There exist constants $\lambda_1,\ldots, \lambda_{\kappa}>0$ such that for any  bounded and compactly supported function $f: \cV_Q^+\to \C$ \comm{and for any $g\in G$},
\begin{align*}
S_Qf\comm{(g)}=\sum_{i=1}^{\kappa}E_i(\comm{g}\,|\, f_{\lambda_i}),
\end{align*}
where for any $\lambda>0$, $f_{\lambda}(\bm{v}):=f(\lambda^{-1}\bm{v})$.
\end{Lem}

\begin{remark}\label{rmk:constinstr}
The constants $\lambda_i$ ($1\leq i\leq \kappa$) measure the size of cusps of the lattice $\G$, and are intrinsic to $\G$ (hence also to $Q$); see \cite[Lemmas 2.3 and 3.1]{KelmerYu2022a}.
\end{remark}
In view of \lemref{lem:sumincomein}, we have
\begin{align}\label{equ:firsmompre}
\int_{Y_Q}S_Qf(g)\, \d\mu_{Q}(g)=\sum_{i=1}^{\kappa}\int_{Y_Q}E_i(g\,|\, f_{\lambda_i})\, \d\mu_{Q}(g),
\end{align}
and
\begin{align}\label{equ:secmompre}
\langle S_Qf, S_Qf\rangle&=\sum_{i,j=1}^{\kappa}\langle E_i(\cdot\,|\, f_{\lambda_i}), E_j(\cdot\,|\, f_{\lambda_j})\rangle.
\end{align}
Thus in order to compute the moments of $S_Qf$, it suffices to prove moment formulas of incomplete Eisenstein series which can be computed using the spectral approach in \cite{KelmerMohammadi12,Yu17}. For reader's convenience, we outline this spectral computation in next subsection. 

\subsection{Moment formulas of incomplete Eisenstein series}

Using a standard unfolding argument we can deduce the following first moment formula of $E_i(g\,|\, f)$ whose standard proof we omit.
\begin{Lem}\label{lem:firstmom}
	For any bounded and compactly supported $f:\cV_Q^+\to \C$ and for each $1\leq i\leq \kappa$,
	\begin{align}\label{equ:firstmoment}
		\int_{Y_Q}E_i(g\,|\, f)\,d\mu_{Q}(g)=\omega_{\G}m_{\cV_Q^+}(f),
		\end{align}
	where $\omega_{\G}$ is the constant residue as in \eqref{equ:rescovol} and $m_{\cV_Q^+}$ is the $G$-invariant measure on $\cV_{Q}^+$ as in \eqref{equ:lebmea}.
	\end{Lem}

Next, we prove an integration formula for products of two incomplete Eisenstein series.
We first introduce some relevant notation. 

Let $c_{\tau}$ be the $\tau$-conjugation isomorphism from $K=K_Q$ to $K_{Q_n}$ and let $\iota: K_{Q_n}\to \SO_{n+1}(\R)$ be the natural isomorphism sending $k=\left(\begin{smallmatrix} \tilde{k} &  \\  &1\end{smallmatrix}\right)\in K_{Q_n}$ to $\tilde{k}\in \SO_{n+1}(\R)$. Then the isomorphism $\phi=\iota\circ c_{\tau}$ induces an isomorphism of representations between $L^2(M\bk K)$ (as a $K$-representation) and $L^2(S^n)$ (as a $\SO_{n+1}(\R)$-representation). Here the $L^2$-norms are with respect to $\mu_{M\bk K}$ and $\sigma_n$ respectively.
By the classical \comm{spherical} harmonic analysis, the function space $L^2(S^n)$ decomposes into irreducible $\SO_{n+1}(\R)$-representations as following:
\begin{align*}
L^2(S^n)=\bigoplus_{d\geq 0}L^2(S^n,d),
\end{align*}
where $L^2(S^n, d)$ is the space of degree $d$ harmonic polynomials in $n+1$ variables restricted to $S^n$. This in turn induces the following decomposition of $L^2(M\bk K)$ into irreducible $K$-representations
\begin{align*}
	L^2(M\bk K)=\bigoplus_{d\geq 0} L^2(M\bk K, d),
\end{align*}
where $L^2(M\bk K, d)$ is the pre-image of $L^2(S^n,d)$ under the above isomorphism between $L^2(M\bk K)$ and $L^2(S^n)$.
For each $d\geq 0$, let us fix an orthonormal basis 
$
\{\psi_{d,l}: 1\leq l\leq \dim_{\C}L^2(M\bk K, d)\}
$
for $L^2(M\bk K, d)$. 
For any $f: \cV_Q^+\to \C$ bounded and of compact support, let 
\begin{align}\label{equ:spexpan}
	f_{d,l}(y):=\int_{M\bk K}f(\bm{e}_0a_yk)\overline{\psi_{d,l}(k)}\,\d\mu_{M\bk K}(k).
\end{align}
so that $f$ has a \textit{spherical expansion} 
\begin{align}\label{equ:sphexpan}
f(\bm{e}_0a_yk)=\sum_{d, l}f_{d,l}(y)\psi_{d,l}(k)
\end{align}
\comm{in $L^2$ and also pointwise if $f$ is smooth}.
\comm{For any function $\rho$ on $\R_{>0}$, we denote by $\hat{\rho}(s):=\int_0^{\infty}\rho(y)y^{-(s+1)}\,\d y$ ($s\in \C$) its \textit{Mellin transform}  whenever this defining integral is absolutely convergent.}
Since $\psi_{0,1}=1$ is constant, we have
\begin{align}\label{equ:spaceave}
m_{\cV_Q^+}(f)=\int_0^{\infty}\int_{M\bk K}f(\bm{e}_0a_yk)\, \d\mu_{M\bk K}(k)y^{-(n+1)}\, \d y=\int_0^{\infty}f_{0,1}(y)y^{-(n+1)}\,\d y=\hat{f}_{0,1}(n).
\end{align}
\comm{Using the spherical expansion we define the following bilinear form on the space $C^\infty_c(\cV_Q^+)$: For any $f,f'\in C^\infty_c(\cV_Q^+)$ and any $s\in  (\tfrac{n}{2}, n)$ define
\begin{align}\label{def:mffun}
	M_{f,f'}(s):=\sum_{d,l}P_d(s)\hat{f_{d,l}}(s)\overline{\hat{f'_{d,l}}(s)},
\end{align}
with $P_1(s):=1$ and $P_d(s):=\prod_{i=0}^{d-1}\frac{n-s+i}{s+i}$ if $d\geq 1$.} 
We note that for each $d\geq 0$, $P_d(s)$ satisfies the functional equation $P_d(s)P_d(n-s)=1$. In particular, when $\Re(s)=\frac{n}{2}$ we have
	\begin{align}\label{equ:funcons2}
		|P_d(s)|^2=P_d(s)P_d(\bar{s})=P_d(s)P_d(n-s)=1.
		\end{align}
We now state our integration formula.
\begin{Thm}\label{thm:secmoment}
	For any $1\leq i,j\leq \kappa$, there exists a bounded linear operator $\cT_{ij}: L^2(\cV_Q^+)\to L^2(\cV_Q^+)$ with operator norm $\|\cT_{ij}\|_{\rm op}\leq 1$ such that for any $f,f'\in C^\infty_c(\cV_Q^+)$, 
	\begin{align*}
		\langle E_i(\cdot\,|\, f),E_j(\cdot\,|\, f')\rangle&=\omega_{\G}^2m_{\cV_Q^+}(f)m_{\cV_Q^+}(\overline{f'})+\omega_{\G}\left\langle \delta_{ij}f+\cT_{ij}(f), f'\right\rangle+\omega_{\G}\sum_{s_l\in\cC_{\G}}M_{f,f'}(s_l)\Res_{s=s_l}\varphi_{ij}(s),
	\end{align*}
	\end{Thm}
\begin{proof}
	For each $d\geq 0$, let $\{\psi_{d,l}: 1\leq l\leq \dim_{\C}L^2(M\bk K, d)\}$ be the orthonormal basis for $L^2(M\bk K, d)$ fixed as above. Let $f, f': \cV_Q^+\to \C$ be two \comm{smooth} and compactly supported functions. 
	Doing an unfolding argument, using the left $\G$-invariance of $E_i(g\,|\, f)$, and using the Haar measure description \eqref{equ:haar3} we get
	\begin{align}\label{equ:presecfo}
			\langle E_i(\cdot\,|\, f),E_j(\cdot\,|\, f')\rangle&
			=\omega_{\G}\int_0^{\infty}\int_{M\bk K}\overline{f'(\bm{e}_0a_yk)}\cI_{ij}(f)(\bm{e}_0a_yk)\,y^{-(n+1)}\, \d y\d\mu_{M\bk K}(k),
	\end{align} 
	where
	\begin{align*}
	\cI_{ij}(f)(\bm{e}_0a_yk):=\int_{\tau_j^{-1}\G_{P_j}\tau_j\bk UM}E_i(\tau_j u_{\bm{x}}ma_yk\,|\, f)\, \d\bm{x}\d\mu_{M}(m).
	\end{align*}
Now write $f(\bm{e}_0a_yk)=\sum_{d,l}f_{d,l}(y)\psi_{d,l}(k)$ and define $\cT_{ij}(f)$ by 
\begin{align}\label{equ:defingequtij}
\cT_{ij}(f)(\bm{e}_0a_yk)&:=\cI_{ij}(f)(\bm{e}_0a_yk)-\delta_{ij}f(\bm{e}_0a_yk)-\omega_{\G}\hat{f}_{0,1}(n)\\
&-\sum_{s_l\in\cC_{\G}}\sum_{d,l}P_d(s)\hat{f}_{d,l}(s)y^{n-s}\psi_{d,l}(k)\Res_{s=s_l}\varphi_{ij}(s).\nonumber
\end{align}
This relation equivalently gives a formula for $\cI_{ij}(f)$. We note that the desired second moment formula can be checked to follow by plugging this formula into \eqref{equ:presecfo} and noting that $\hat{f}_{0,1}(n)=m_{\cV_Q^+}(f)$ (see \eqref{equ:spaceave}).
Since $C_c^{\infty}(\cV_Q^+)$ is dense in $L^2(\cV_Q^+)$, any bounded linear operator $\cT$ defined on $C_c^{\infty}(\cV_Q^+)$ uniquely extends to a bounded linear operator on $L^2(\cV_Q^+)$ with the same operator norm. Clearly $\cT_{ij}$ is linear; it thus suffices to show \comm{$\|\cT_{ij}\|_{\rm op}\leq 1$}. 

To show this we will show that for any $f\in C^\infty_c(\cV_Q^+)$
\begin{align}\label{equ:alttijdef}
	\cT_{ij}(f)(\bm{e}_0a_yk)=\sum_{d,l}\frac{1}{2\pi i}\int_{(\tfrac{n}{2})}P_d(s)\varphi_{ij}(s)\hat{f}_{d,l}(s)y^{n-s}\psi_{d,l}(k)\, \d s.
\end{align}
\comm{Here for any $\sigma\in \R$ and any function $F$ on $\C$, $\frac{1}{2\pi i}\int_{(\sigma)}F(s)\, \d s:=\frac{1}{2\pi}\int_{-\infty}^{\infty}F(\sigma+it)\, \d t$.} 
We note that assuming this claim, then applying the Plancherel's theorem and using the estimate $|P_d(s)\varphi_{ij}(s)|\leq 1$ for $\Re(s)=\frac{n}{2}$ (cf. \eqref{equ:funcons1} and \eqref{equ:funcons2}), we see that indeed $\|\cT_{ij}(f)\|_2\leq \|f\|_2$, implying the desired operator norm bound. It thus remains to prove \eqref{equ:alttijdef}.

Now using the fact that $\G_{U_j}$ is a finite-index subgroup of $\G_{P_j}$ (see \cite[Lemma 2.1]{KelmerYu2022a}) we get 
\begin{align*}
	\cI_{ij}(f)(\bm{e}_0a_yk)=\frac{1}{[\G_{P_j}: \G_{U_j}]}\int_{\tau_j^{-1}\G_{U_j}\tau_j\bk U}\int_M E_i(\tau_ju_{\bm{x}}ma_yk\,|\, f)\, \d\bm{x} \d\mu_M(m).
\end{align*}
Next, note that by our normalization, 
$$
1=\mu_L(\tau_j^{-1}\G_{P_j}\tau_j\bk UM)=\frac{\mu_L(\tau_j^{-1}\G_{U_j}\tau_j\bk UM)}{[\G_{P_j}: \G_{U_j}]}
=\frac{\vol(\tau_j^{-1}\G_{U_j}\tau_j\bk U)}{[\G_{P_j}: \G_{U_j}]}.
$$
This implies that $[\G_{P_j}: \G_{U_j}]=\vol(\tau_j^{-1}\G_{U_j}\tau_j\bk U)$. We can then use this identity and change the order of integration to get 
$$\cI_{ij}(f)(\bm{e}_0a_yk)=
	\int_M\comm{\cP_{ij}(f)}(ma_yk)\, \d\mu_M(m),$$
where
\begin{align*}
\comm{\cP_{ij}(f)}(ma_yk):=\frac{1}{\vol(\tau_j^{-1}\G_{U_j}\tau_j\bk U)}\int_{\tau_j^{-1}\comm{\G_{U_j}}\tau_j\bk U}E_i(\tau_ju_{\bm{x}}ma_yk\,|\, f)\, \d\bm{x}.
\end{align*}
First if $f$ is spherical, i.e. $f(\bm{e}_0a_yk)=\rho(y)$ for some $\rho\in C_c^{\infty}(\R_{>0})$, applying Mellin inversion formula $\rho(y)=\frac{1}{2\pi i}\int_{(\sigma)}\hat{\rho}(s)y^{\comm{s}}\, \d s$ ($\sigma\in \R$) \comm{(see e.g. \cite[p.\ 55-56]{Bump1997})} and the constant term formula \eqref{e:ConstSph1} we get that for any $\sigma>n$
\begin{align*}
	\comm{\cP_{ij}(f)}(ma_yk)&=\delta_{ij}\rho(y)+\frac{1}{2\pi i}\int_{(\sigma)}\varphi_{ij}(s)\hat{\rho}(s)y^{n-s}\, \d s.
\end{align*}
Next, for functions of the form $f(\bm{e}_0a_yk)=\rho(y)\psi_{d,l}(k)$ by applying suitable \comm{raising} operators to the above identity (see \cite[p. 460-463]{Yu17}, especially \cite[(3.12)]{Yu17}) we get that 
for any $\sigma>n$,
\begin{align*}
	\comm{\cP_{ij}(f)}(ma_yk)&=\delta_{ij}f(\bm{e}_0a_yk)+\frac{1}{2\pi i}\int_{(\sigma)}P_d(s)\varphi_{ij}(s)\hat{\rho}(s)y^{n-s}\psi_{d,l}(k)\, \d s.
\end{align*}
Finally, for a general function, we can expand $f(\bm{e}_0a_yk)=\sum_{d,l} f_{d,l}(y)\psi_{d,l}(k)$ to get that
\begin{align*}
	\comm{\cP_{ij}(f)}(ma_yk)&=\delta_{ij}f(\bm{e}_0a_yk)+\sum_{d,l}\frac{1}{2\pi i}\int_{(\sigma)}P_d(s)\varphi_{ij}(s)\hat{f}_{d,l}(s)y^{n-s}\psi_{d,l}(k)\, \d s.
\end{align*}
In particular, we see that $\comm{\cP_{ij}(f)}(ma_yk)$ is independent of the variable $m\in M$. Thus $\cI_{ij}(f)=\comm{\cP_{ij}(f)}$, or equivalently, 
$$
\left(\cI_{ij}(f)-\delta_{ij}f\right)(\bm{e}_0a_yk)=
	\sum_{d,l}\frac{1}{2\pi i}\int_{(\sigma)}P_d(s)\varphi_{ij}(s)\hat{f}_{d,l}(s)y^{n-s}\psi_{d,l}(k)\, \d s.
	$$
We can now further simplify the right hand side by shifting the contour of integration to the critical line $\Re(s)=\frac{n}{2}$ and picking up the contribution of the exceptional poles $s_{l}\in \cC_{\G}$ and the trivial pole $s=n$ (and noting that $P_d(n)=\delta_{d0}$ and $\Res_{s=n}\varphi_{ij}(s)=\omega_{\G}$ (see \eqref{equ:constomeg})) to get 
\begin{align*}
	\left(\cI_{ij}(f)-\delta_{ij}f\right)(\bm{e}_0a_yk)&=\omega_{\G}\hat{f}_{0,1}(n)+\sum_{s_l\in\cC_{\G}}\sum_{d,l}P_d(s)\hat{f}_{d,l}(s)y^{n-s}\psi_{d,l}(k)\Res_{s=s_l}\varphi_{ij}(s)\\
	&+\sum_{d,l}\frac{1}{2\pi i}\int_{(\tfrac{n}{2})}P_d(s)\varphi_{ij}(s)\hat{f}_{d,l}(s)y^{n-s}\psi_{d,l}(k)\,\d s.\nonumber
\end{align*}
The formula \eqref{equ:alttijdef} then follows immediately by comparing \eqref{equ:defingequtij} and the above equation. This finishes \comm{the proof}.
	\end{proof}
\begin{remark}\label{rmk:scalar}
	For later purpose, we note that \comm{for any $f\in C_c^{\infty}(\cV_Q^+)$} if we define $f_{\lambda}(\bm{v}):=f(\lambda^{-1}\bm{v})$ ($\lambda>0$), then the relation $\bm{e}_0a_y=y^{-1}\bm{e}_0$ ($y>0$), implies that $(f_{\lambda})_{d,l}(y)=f_{d,l}(\lambda y)$. Hence, for any $\lambda_1,\lambda_2>0$ and any $s\in (\frac{n}{2}, n)$,
	\begin{align}\label{equ:idmf}
		M_{f_{\lambda_1}, f'_{\lambda_2}}(s)=\lambda_1^s\lambda_2^sM_{f,f'}(s).
		\end{align}
	\end{remark}

\subsection{Proof of Theorem \ref*{thm:mainrest}}\label{sec:lightconeeisen}
 Before giving the proof, let us give a quick review of the light-cone Eisenstein series studied in \cite{KelmerYu2022a}. 
 Let $Q(\bm{v})=Q_n(\bm{v}\tau)$ with $\tau=ka$ as in \rmkref{rmk:formoftau}, and define the $K$-invariant norm $\|\cdot\|_Q$ on $\R^{n+2}$ by 
\begin{align}\label{equ:qnorm}
\|\bm{v}\|_Q:=\|\bm{v}\tau\|.
\end{align}
Observe that this norm does not depend on the choice of orthogonal matrix $k$. Recall from \cite{KelmerYu2022a} that the \textit{light-cone Eisenstein series} of $Q$ is defined for $\Re(s)>n$ by the series
\begin{align}\label{equ:ligconeeise}
E_Q(s,g)&:=\|\bm{e}_0\|_Q^s\sum_{\bm{v}\in \cL_Q}\|\bm{v}g\|_Q^{-s}.
\end{align}

Similar to the light-cone Siegel transform, $E_Q(s,g)$ can be written as a sum of Eisenstein series of $\G$ at all cusps; see \cite[Corollary 3.2]{KelmerYu2022a}. Explicitly, we have for any $\Re(s)>n$,
\begin{align}\label{equy:relaeisen}
E_Q(s,g)&=\sum_{i=1}^{\kappa}\lambda_i^sE_i(s,g),
\end{align}
where $\lambda_1,\ldots, \lambda_{\kappa}$ are as in \lemref{lem:sumincomein}.
In particular, we see from this relation that $E_Q(s,g)$ has a meromorphic continuation to the whole $s$-plane and is holomorphic on the half plane $\Re(s)\geq \frac{n}{2}$ except for a simple pole at $s=n$ with constant residue, $\omega_Q:=\Res_{s=n}E_Q(s,g)$, and at most finitely many simple poles in $(\tfrac{n}{2}, n)$ contained in the set $\cC_{\G}$. Moreover, it follows from \cite[Theorem 1.8]{KelmerYu2022a} that  in fact there is at most one such pole located at $s_n=\left \lfloor{\tfrac{n+2}{2}}\right \rfloor$ if $n\geq 3$ and no poles for $n<3$.
The residue $\omega_Q$ then satisfies the relation
\begin{align}\label{equ:omegaq}
\omega_Q=\omega_{\G}\sum_{i=1}^{\kappa}\lambda_i^n.
\end{align}

\begin{proof}[Proof of \thmref{thm:mainrest}]
The first moment formula follows immediately by combining \lemref{lem:firstmom} and the relations \eqref{equ:firsmompre} and \eqref{equ:omegaq}.
	
For the second moment formula, \comm{further assume $f$ is smooth}. For each $1\leq i,j\leq \kappa$ we use
	\thmref{thm:secmoment} together with the observation that for any $\lambda>0$ we have  $m_{\cV_Q^+}(f_{\lambda})=\lambda^nm_{\cV_Q^+}(f)$ and the relation \eqref{equ:idmf} to get that 
\begin{align*}
     \langle E_i(\cdot\,|\, f_{\lambda_i}),E_j(\cdot\,|\, f_{\lambda_j})\rangle&=
     \omega_{\G}^2\lambda_i^n\lambda_j^n|m_{\cV_Q^+}(f)|^2+\omega_{\G}\sum_{s_l\in\cC_{\G}}\lambda_i^{s_l}\lambda_j^{s_l}M_{f,f}(s_l)\Res_{s=s_l}\varphi_{ij}(s)+O_{Q}\left(m_{\cV_Q^+}(|f|^2)\right).
\end{align*}
Next, for each pole $s_l\in \cC_\G$,  expanding $E_Q(s,g)=\sum_{j=1}^{\kappa} \lambda_j^s E_j(s,g)$ and taking residues we get that 
$$\Res_{s=s_l}E_Q(s,g)=\sum_{j=1}^{\kappa} \lambda_j^{s_l} \Res_{s=s_l}E_j(s,g).$$
Taking an inner product and using \eqref{eq:scatteringres} we get that
\begin{align*}
	\|\Res_{s=s_l}E_Q(s,g)\|^2&=\sum_{i,j=1}^{\kappa}\lambda_i^{s_l}\lambda_j^{s_l}\langle \Res_{s=s_l}E_i(s,\cdot),\Res_{s=s_l}E_j(s,\cdot)\rangle\\
	&=\sum_{i,j=1}^{\kappa}\lambda_i^{s_l}\lambda_j^{s_l}\Res_{s=s_l}\varphi_{ij}(s).
\end{align*}
Using the relations \eqref{equ:secmompre}, \eqref{equ:omegaq} we get that 
\begin{align*}
     	\langle S_Qf, S_Qf\rangle&=
     	\left|\omega_Qm_{\cV_Q^+}(f)\right|^2+\omega_\G \sum_{s_l\in\cC_{\G}} \|\Res_{s=s_l}E_Q(s,g) \|^2 M_{f,f}(s_l)+O_{Q}\left(m_{\cV_Q^+}(|f|^2)\right).
\end{align*}
Since $\Res_{s=s_l}E_Q(s,g)=0$ unless $s_l$ is a pole of $E_Q(s,g)$ we get that the only secondary term comes from the possible pole at $s=s_n$ and setting
$c_Q=\omega_\G\|\Res_{s=s_n}E_Q(s,g) \|^2$ concludes the proof.

	\end{proof}

	\subsection{Applications to discrepancies}\label{sec:varbd}
	{In this subsection we prove the discrepancy bound in \corref{cor:vabd}
	
using the moment formulas proved above  \comm{together with some approximating arguments}. 
As mentioned in the introduction, the case when $c_Q> 0$ is more tricky due to the existence of an exceptional pole. In this case, \comm{we need to bound the secondary term which} we can only prove for certain Borel sets defined below. }
\begin{Def}\label{def:gesec}
We say a subset $B\subseteq \cV_{Q}^+$ is a \textit{generalized sector} if $B$ is \comm{of finite measure} and its indicator function is of the form $\chi_B(\bm{e}_0a_yk)=\rho(y)\phi(k)$, where $\rho$ and $\phi$ are indicator functions of some Borel sets in $\R_{>0}$ and $M\bk K$ respectively.  
\end{Def}

\begin{remark}\label{rmk:diffculty}
\comm{The main benefit of working with generalized sectors is that their indicator functions naturally \textit{separate the variables}. In practice, we can take $f\in C_c^{\infty}(\cV_Q^+)$ of the form $f(\bm{e}_0a_yk)=\rho(y)\phi(k)$ to approximate the indicator function of a generalized sector. For such $f$, the secondary term $M_{f,f}(s_n)$ has a simpler expression (see \eqref{equ:expseofvar} below) from which one can invoke some bounds from spherical harmonic analysis and use the fact that $\rho$ is close to some indicator function to bound $M_{f,f}(s_n)$ in terms of the $L^1$ and $L^2$ norms of $f$. For $f\in C_c^{\infty}(\cV_Q^+)$ approximating a general finite-measure Borel set, such a simpler expression is not available and it is not clear to us how to bound $M_{f,f}(s_n)$.}
\end{remark}

		\begin{proof}[Proof of \corref{cor:vabd}]
	 First assume that $c_Q=0$. In this case, using the second moment formula for smooth compactly supported functions, we can deduce the same formula for measurable, bounded and compactly supported functions. Indeed, for any measurable, bounded and compactly supported $f: \cV_Q^+: \to \C$, we can find a sequence $\{f_j\}_{j\in\N}\subset C^\infty_c(\cV_Q^+)$ with $f_j$ converging to $f$ in $L^1$ and in $L^2$ (and hence, after taking a subsequence  also pointwise almost everywhere). We then have that $S_Q f_j$ converges pointwise almost everywhere to $S_Qf$. 
	 To show that this convergence is also in $L^2(Y_Q)$ note that the sequence $\{f_j\}_{j\in\N}$ is a Cauchy sequence in $ L^1(\cV_Q^+)\cap L^2(\cV_Q^+)$ and we can use the second moment formula to bound for any $i\neq j$,
	 $$\|S_Qf_i-S_Qf_j\|^2=\langle S_Q(f_i-f_j),S_Q(f_i-f_j)\rangle\ll_Q \|f_i-f_j\|^2_1+\| f_i-f_j\|^2_2.$$
	Hence the sequence $\{S_Qf_j\}_{j\in\N}$ is a Cauchy sequence in $L^2(Y_Q)$ and thus converges to $S_Qf$ in $L^2(Y_Q)$. We can now write 
	 $\langle S_Qf, S_Qf\rangle=\lim_{j\to \infty} \langle S_Qf_j, S_Qf_j\rangle$ and deduce the second moment formula for $S_Qf$ from the second moment formula for $S_Qf_j$. Moreover, if $f=\chi_B$ is the indicator function of some measurable set $B$ of finite measure, then we can approximate $f$ pointwise from below by a monotone sequence $\{\chi_{B_j}\}_{j\in\N}$ with $\chi_{B_j}$ measurable and of compact support. Hence $S_Q\chi_{B_j}\to S_Qf$ pointwise monotonously and by the monotone convergence theorem we conclude that the same \comm{first and} second moment formula holds for $S_Qf$.
	 
	Now let  $f=\chi_B$  so that $S_Qf(g)=\#(\cL_Qg\cap B)$. 
	 Expanding the square $$D(\cL_Q g,B)^2=\left|S_Qf(g)-\omega_Qm_{\cV_Q^+}(B)\right|^2=|S_Qf(g)|^2+\omega^2_Qm_{\cV_Q^+}(B)^2-2S_Qf(g)\omega_Qm_{\cV_Q^+}(B),$$ 
	 and applying \eqref{equ:mainthmfirst} we get
		\begin{align*}
		\int_{Y_Q}\left|D(\cL_Q g, B)\right|^2 \,\d\mu_{Q}(g)&=\int_{Y_Q}\left(S_Qf(g)\right)^2\, \d\mu_{Q}(g)-\omega_Q^2m_{\cV_Q^+}(B)^2.
		\end{align*}
		Recall we assume $c_Q=0$ and hence \eqref{equ:disbdwiexpo} immediately follows from  \eqref{equ:mainthmsecond}.

		Next when $c_Q>0$ ({in particular, we can assume $n\geq 3$}), we need to further assume that $B$ is a generalized sector. In this case we can write $f=\chi_B$ as  $f(\bm{e}_0a_yk)=\rho(y)\phi(k)$ for some indicator functions $\rho$ and $\phi$ and note that $m_{\cV_Q^+}(B)=\hat{\rho}(n)\mu_{M\bk K}(\phi)$ and $\mu_{M\bk K}(\phi)=\|\phi\|_1=\|\phi\|_2^2$. 
		
We first assume $f$ is of compact support. 
 While we cannot apply the second moment formula directly to $f$, we can find a sequence of smooth non-negative compactly supported functions $f_j(\bm{e}_0a_yk)=\rho_j(y)\phi_j(k)$ with $\rho_j\to \rho$ and $\phi_j\to \phi$ pointwise almost everywhere. Then also $S_Qf_j\to S_Qf$ pointwise almost everywhere. Since $f$ is of compact support we can find some non-negative, smooth compactly supported function dominating $f$ and all $f_j$  and then by dominated  convergence we have $\langle S_Qf, S_Qf\rangle=\lim_{j\to\infty} \langle S_Qf_j, S_Qf_j\rangle$.  We can now use the second moment formula for $S_Qf_j$ and take a limit. However, to get our result we need an estimate on the secondary term $M_{f_j,f_j}(s_n)$.  Using the fact that $f_j(\bm{e}_0a_yk)=\rho_j(y)\phi_j(k)$ factors it is not hard to see that for any $s\in (\frac{n}{2}, n)$,
		\begin{align}\label{equ:expseofvar}
			M_{f_j,f_j}(s)=\left|\hat{\rho}_j(s)\right|^2\sum_{d\geq 0}P_d(s)\|\phi_{j,d}\|_2^2,
		\end{align}
		where for any $d\geq 0$, $\phi_{j,d}=\text{pr}_d(\phi_j)$ with $\text{pr}_d: L^2(M\bk K)\to L^2(M\bk K, d)$ the projection operator from $L^2(M\bk K)$ to $L^2(M\bk K,d)$.
		It was shown in \cite[p. 468-470]{Yu17} that for any function $\phi \in L^2(M\bk K)$
				$$
		\sum_{d\geq 0}P_d(s)\|\phi_d\|_2^2\ll_{s,n} \|\phi\|^{2(\frac{2s}{n}-1)}_1\|\phi\|^{4(1-\frac{s}{n})}_2.
		$$
		Hence 
			$$
			M_{f_j,f_j}(s)\ll_{s,n}
			\left|\hat{\rho_j}(s)\right|^2\|\phi_j\|^{2(\frac{2s}{n}-1)}_1\|\phi_j\|^{4(1-\frac{s}{n})}_2,
			$$ 
			and we can take the limit and apply the dominated convergence theorem to the right side of the above inequality to get
			$$\limsup_{j\to\infty} 	M_{f_j,f_j}(s)\ll_{s,n} \left|\hat{\rho}(s)\right|^2\|\phi\|^{2(\frac{2s}{n}-1)}_1\|\phi\|^{4(1-\frac{s}{n})}_2=\left|\hat{\rho}(s)\right|^2\mu_{M\bk K}(\phi)^{\frac{2s}{n}}.$$
			Moreover, using the fact that $\rho$ is an indicator function and applying H\"{o}lder's inequality one can show, 
		$	\hat{\rho}(s)\leq 2\hat{\rho}(n)^{\frac{s}{n}}$; 
		see \cite[Lemma 4.3.1]{Yu2018} and recall that we assume here $n\geq 3$. We can thus bound 
		$$\limsup_{j\to\infty} 	M_{f_j,f_j}(s)\ll_{s,n}  \hat{\rho}(n)^{\frac{2s}{n}}\mu_{M\bk K}(\phi)^{\frac{2s}{n}}=m_{\cV_Q^+}(B)^{\frac{2s}{n}}.$$
		Taking the limit in the second moment formula for $S_Qf_j$ and applying the above bound we can conclude the proof for a precompact generalized sector. Finally for a general generalized sector $B$, we can similarly approximate $\chi_B$ from below by a monotone sequence $\{\chi_{B_j}\}_{j\in \N}$  of indicator functions of precompact generalized sectors $B_j$ and use the monotone convergence theorem to deduce the second moment formula
\begin{align*}
\int_{Y_Q}\left|S_Q\chi_B(g)\right|^2\, \d\mu_Q(g)=\omega_Q^2m_{\cV_Q^+}(B)^2+O_Q\left(m_{\cV_Q^+}(B)^{\beta_Q}+m_{\cV_Q^+}(B)\right),
\end{align*} 
and the desired discrepancy bound.	
		\end{proof}

	\begin{remark}\label{rmk:diffbd}
	\comm{For applications we will need to apply the discrepancy bound \eqref{equ:disbdwiexpo} for difference sets of two nested generalized sectors which may no longer be a generalized sector. Nevertheless, these difference sets can be expressed as a disjoint union of at most two generalized sectors.
	 To handle these sets, for any $k\in \N$} using the definition \eqref{def:mffun}, the fact that the map $f\mapsto \hat{f}_{d,l}$ is linear in $f$ for any $d\geq 0, 1\leq l\leq \dim_{\C} L^2(S^n,d)$ and the estimate $\left|\sum_{i=1}^kz_i\right|^2\leq k\sum_{i=1}^m\left|z_i\right|^2$ ($z_i\in\C$) 
	we have \comm{for any $s\in (\frac{n}{2}, n)$ and} for any $k$ \comm{smooth} compactly supported functions $f_i: \cV_Q^+\to \C$, and $f:=\sum_i f_i$
	\begin{align*}
	M_{\comm{f, f}}(s)\leq k\left(M_{f_1,\comm{f_1}}(s)+\cdots \comm{+}M_{f_k,\comm{f_k}}(s)\right).
	\end{align*}
	In particular, if $B=\bigsqcup_{i=1}^k B_i$ is a disjoint union of $k$ generalized sectors $B_1,\ldots, B_k$, then \comm{using the above estimate together with an approximating argument we get}
	\begin{align*}
	\int_{Y_Q}\left|D(\cL_Q g, B)\right|^2 \, \d\mu_{Q}(g)&
	\ll_{Q, k}m_{\cV_Q^+}(B)+ \sum_{i=1}^km_{\cV_Q^+}(B_i)^{\beta_Q}\leq m_{\cV_Q^+}(B)+m_{\cV_Q^+}(B)^{\beta_Q}.
	\end{align*}
	\end{remark}
\subsection{Conjugation}\label{rmk:expcord}\label{sec:conj}
Let $\tau=ka\in \GL_{n+2}(\R)$ be as in \rmkref{rmk:formoftau}.
Later for our applications, we will need to work with explicit coordinates, hence, instead of working on $Y_Q$, it will be more convenient to work on homogeneous spaces of $\SO_{Q_n}^+(\R)$ for which we can use the explicit coordinates introduced in section \ref{sec:cormea}. For this, we will translate the above discrepancy bound on $Y_Q$ to equivalent bounds on a homogeneous space of $\SO_{Q_n}^+(\R)$ via the translating matrix $\tau$. 

Recall that $\SO_Q^+(\R)=\tau \SO_{Q_n}^+(\R)\tau^{-1}$, and note that the $\tau$-conjugation map $c_{\tau}: \SO_Q^+(\R)\to \SO_{Q_n}^+(\R)$ sending $h\in \SO_Q^+(\R)$ to $\tau ^{-1}h\tau\in \SO_{Q_n}^+(\R)$ induces a homeomorphism between the homogenous spaces $Y_Q$ and $\tilde{Y}_Q:=\tau^{-1}\SO_Q^+(\Z)\tau\bk \SO_{Q_n}^+(\R)$ with the pushforward of $\mu_Q$ being the probability $\SO_{Q_n}^+(\R)$-invariant measure on $\tilde{Y}_Q$, which we denote by $\tilde{\mu}_{Q}$.
Moreover, recall that $\cV_{Q_n}^+=\cV_{Q}^+\tau$ and note that the pushforward of $m_{\cV_Q^+}$ under this right $\tau$-multiplication map 
is exactly $m_{\cV_{Q_n}^+}$, and \comm{one can also deduce} that a Borel set $B\subseteq \cV_Q^+$ is a generalized sector if and only if $B\tau\subseteq \cV_{Q_n}^+$ is a generalized sector. Denote by $\tilde \cL_Q=\cL_Q \tau \subseteq \cV_{Q_n}^+$, then  the discrepancy bound proved in \corref{cor:vabd} is equivalent to the following discrepancy bound on $\tilde{Y}_Q$: For any finite-measure Borel subset $B\subseteq \cV_{Q_n}^+$ \comm{(which is further assumed to be a generalized sector if $c_Q>0$)}, 
       \begin{align}\label{equ:disbdwiexpotran}
			\int_{\tilde{Y}_Q}\left|\#(\tilde{\cL}_Q g\cap B)-\omega_{Q}m_{\cV^+_{Q_n}}(B)\right|^2 \, \d\tilde{\mu}_{Q}(g)\ll_Q 	m_{\cV^+_{Q_n}}(B)^{\beta_Q}+m_{\cV^+_{Q_n}}(B).
		\end{align}

\section{From mean square discrepancy to counting}\label{sec:meantocount}
In this section we show how one can use a mean square bound on the discrepancy in order to obtain effective counting estimates. These arguments are quite general and can be applicable in other settings so we present them in greater generality. For this we assume that $\cV\subseteq \R^n$ is some algebraic variety with a transitive action of a real unimodular algebraic group $G$. Let $\G\leq G$ be a non-uniform lattice and $\cL\subseteq \cV$ an infinite discrete set that is stable under the right action of $\G$, i.e. $\cL \gamma=\cL$ for any $\g\in \G$. Let $\mu_{\G}$ be the probability $G$-invariant measure on $\G\bk G$. With slight abuse of notation, we also denote by $\mu_{\G}$ the Haar measure of $G$ that locally agrees with $\mu_{\G}$. We further assume that there is an appropriately normalized $G$-invariant measure $m_\cV$ on $\cV$ such that for any Borel set $B\subseteq \cV$ of finite measure
$$\int_{\G\bk G} \#(\cL g\cap B)\, \d\mu_{\G}(g)=m_\cV(B).$$
For any $g\in G$ and $B\subseteq \cV$ of finite measure, define the discrepancy $\comm{D(\cL g, B)}=|\#(\cL g\cap B)-m_\cV(B)|$. Assuming a mean square bound for the discrepancy that holds for all sets is sometimes too strong, but in many cases we can have a bound on a large family of sets. To make this precise,
for $\beta\in [1,2)$ and $c>0$ let $\cA_{\beta,c}$ be the family of all finite-measure Borel subsets $B\subseteq \cV$ satisfying that 
\begin{equation}\label{e:DB}
\int_{\G\bk G}|D(\cL g, B)|^2\, \d\mu_{\G}(g)\leq c\left(m_\cV(B)^\beta+m_\cV(B)\right).
\end{equation}
In the rest of this section we show how one can use this assumption to get  various estimates on the counting function $\#(\cL g\cap B_T)$ for families of sets in the family $\cA_{\beta,c}$ under various regularity assumptions.

\begin{remark}\label{rmk:disbd}
Let $Q(\bm{v})=Q_n(\bm{v}\tau)$ be as before and let $\beta_Q$ \comm{be} as in \eqref{equ:defbq}. As discussed in section \ref{sec:conj}, for the case of $\cV=\cV_{Q_n}^+$ and $\cL=\cL_Q\tau$, we can take $G=\SO_{Q_n}^+(\R)$, $\G=\tau^{-1}\SO_Q^+(\Z)\tau$, $\mu_{\G}=\tilde\mu_Q$ and $m_{\cV}=\omega_Q m_{\cV_Q^+}$. In this case, by \eqref{equ:disbdwiexpotran}, 
there exists some $c>0$ depending only on $Q$ such that when $c_Q=0$ (so that $\beta_Q=1$), $\cA_{1,c}$ contains all the finite-measure Borel subsets in $\cV_{Q_n}^+$, and when $c_Q>0$, $\cA_{\beta_Q,c}$ contains any Borel set which is a disjoint union of at most two generalized sectors; see also \rmkref{rmk:diffbd}. Then clearly $\cA_{\beta_Q,c}$ contains all the generalized sectors. In addition, if $B_1, B_2$ are two generalized sectors with $B_1\subseteq B_2$, it is then not difficult to see, using the polar coordinates on $\cV_{Q_n}^+$, that the difference set $B_2\setminus B_1$ is a disjoint union of at most two generalized sectors, thus also belongs to $\cA_{\beta_Q,c}$. 
\end{remark}

\subsection{Counting results for almost all translates}\label{sec:congenr}
Using a classical counting argument by Schmidt \cite{Schmidt1960} one can use the mean square bound on the discrepancy to get very good estimates on $\#(\cL g\cap B_T)$ that hold for almost all $\G g\in \G\bk G$. Following his argument one can show the following.
\begin{Thm}\label{t: Schmidt}
For any increasing family of sets $\{B_T\}_{T>0}$ in $\cA_{\beta,c}$ such that $T\mapsto m_{\cV}(B_T)$ is continuous and the differences $B_{T}\setminus B_{T'} \in \cA_{\beta,c}$ for any $T'<T$ we have for $\mu_{\G}$-a.e. $\G g\in \G\bk G$  
$$\#(\cL g\cap B_T)=m_{\cV}(B_T)+O_{g,\e}(m_{\cV}(B_T)^{\comm{\frac{\beta}{2}}+\e}),\quad\text{as $T\to\infty$}.
$$
\end{Thm}

The proof of this result is essentially identical to the proof of Schmidt's result \cite{Schmidt1960}, and we omit the details; \comm{see also \cite[Theorem 6.1]{KelmerYu2021}.} 
\begin{rem}
We note that when $\cA_{\beta,c}$ contains all finite-measure Borel sets in $\cV$ the condition on the continuity of the function $T\mapsto m_{\cV}(B_T)$ is not needed, as one can artificially add sets to the collection to make it so (and the second condition that differences are in $\cA_{\beta,c}$ holds automatically). In particular in this case the result holds for any increasing family of sets in $\cV$.
\end{rem}

For some applications, one is interested in a family of sets that is not \comm{necessarily} increasing.  \comm{For example, for applications to Diophantine approximation, one is naturally interested in sets that increase in one direction and decrease in another.} In order to handle such cases we give here the following general argument that is more flexible\comm{,} but as a tradeoff we get a weaker bound for the remainder.
\begin{Lem}\label{l:sequences}
Let $\Omega\subseteq G$ be a subset endowed with some finite measure $\mu_{\Omega}$ (which could be singular with respect to $\mu_{\G}$).
Let $\cB$ be a family of Borel sets in $\cV$ satisfying that $m_{\cV}(B)>1$ for any $B\in \cB$ and there exists $0<C<1$ such that for any $B\in \cB$ and for any $0<X\leq C m_{\cV}(B)$    
\begin{equation}\label{e:mDlarge}
\mu_{\Omega}\left(\left\{g\in \Omega: D(\cL g,B)\geq X\right\}\right)\ll_\Omega \left(\tfrac{m_{\cV}(B)}{X}\right)^a m_{\cV}(B)^{-b}
\end{equation}
for some $a>1$ and $0<b\leq 1$. 
\begin{enumerate}
\item Then for any sequence $\{B_k\}_{k\in\N}$ of sets from $\cB$ with 
\begin{equation}\label{e:summability}
\sum_k m_{\cV}(B_k)^{a\nu-b} \log\left(m_{\cV}(B_k)\right)^{-a}<\infty
\end{equation}
for some $0<\nu<\frac{b}{a}$, we have for $\mu_{\Omega}$-a.e. $g\in \Omega$ and for all sufficiently large $k$
 $$D(\cL g,B_k)\leq m_{\cV}(B_k)^{1-\nu} \log(m_\cV(B_k)).$$
 \item  Assume we can find two sequences $\{B_k^\pm\}_{k\in\N}$ in $\cB$ each satisfying \eqref{e:summability}  
such that any $B\in \cB$ with $m_{\cV}(B)$ sufficiently large satisfies $B_k^-\subseteq B\subseteq B_k^+$ for some $k$, and that $m_{\cV}(B_k^+\setminus B_k^-)\ll m_{\cV}(B_k^-)^{1-\nu}\log(m_\cV(B_k^-))$
for all sufficiently large $k$. Then for $\mu_{\Omega}$-a.e. $g\in \Omega$, for all $B\in \cB$ with $m_{\cV}(B)$ sufficiently large 
$$
D(\cL g,B)\ll m_{\cV}(B)^{1-\nu}\log(m_{\cV}(B)).
$$
\end{enumerate}
\end{Lem}
\begin{proof}
For the first part let  $\cM_k=\{g\in \Omega: D(\cL g,B_k)> m_{\cV}(B_k)^{1-\nu} \log(m_\cV(B_k))\}$. We then need to show that the set 
$$\cM_\infty=\bigcap_{m\geq 1}\bigcup_{k\geq m} \cM_k$$
is a null set (with respect to $\mu_{\Omega}$).
Note that the summability assumption \eqref{e:summability} implies $\lim\limits_{k\to\infty}m_{\cV}(B_k)=\infty$. Thus we can apply \eqref{e:mDlarge} to get for all $k$ sufficiently large
$$\mu_{\Omega}(\cM_k)\ll_{\Omega} m_\cV(B_k)^{a\nu-b}\log(m_\cV(B_k))^{-a}.$$ 
\comm{Then by the Borel-Cantelli lemma we get from this estimate and the summability assumption \eqref{e:summability}} 
that $\mu_{\Omega}(\cM_{\infty})=0$.
The second part follows from the first together with the observation that 
$$D(\cL g,B)\leq  \max\left\{D(\cL g,B^-_k),D(\cL g,B_k^+)\right\}+m_{\cV}(B_k^+\setminus B_k^-),$$
whenever $B_k^+\subseteq B\subseteq B_k^+$.
\end{proof}

In particular for an increasing family we get the following.
\begin{Thm}\label{thm:incrcounting}
Let $\cB$ be an increasing family of sets satisfying \eqref{e:mDlarge} for some $a>1$ and $0<b\leq 1$ and that $\{m_{\cV}(B): B\in \cB\}\supseteq [V, \infty)$ for some $V>0$. Then for $\mu_{\Omega}$-a.e. $g\in \Omega$ 
and for all $B\in \cB$ with $m_\cV(B)$ sufficiently large 
$$\#(\cL g\cap B)=m_{\cV}(B)+O\left(m_\cV(B)^{1-\frac{b}{a+1}}\log\left(m_\cV(B)\right)\right).$$
\end{Thm}
\begin{proof}
Let $\nu=\frac{b}{a+1}$ and $\alpha=\frac{1}{b-a\nu}=\frac{1}{\nu}$. Let $\{B_k\}_{k\in\N}\subseteq \cB$ be a sequence with $m_{\cV}(B_k)=\max\{V,k^\alpha \log(k)\}$.
Since our family is increasing then any $B\in \cB$ with $m_{\cV}(B)>V$ satisfies $B_k\subseteq B\subseteq B_{k+1}$ for some $k$ and we can estimate for all $k$ sufficiently large
\begin{align*}
 m_{\cV}(B_{k+1}\setminus B_k)
 \ll k^{\alpha-1}\log(k)\ll k^{\alpha-1}\log(k)^{2-\nu}
 \asymp_{a,b} m_{\cV}(B_k)^{1-\nu}\log(m_{\cV}(B_k)).
\end{align*} 
 We can also estimate 
 $$\sum_{k} m_{\cV}(B_k)^{a\nu-b}\log(m_\cV(B_k))^{-a}\ll \sum_k k^{-1}\log(k)^{-a}<\infty,$$
 and the result now follows from \lemref{l:sequences}.
\end{proof}

In order to compare this to \thmref{t: Schmidt}, note that a standard application of the Chebyshev inequality implies that any set $B\in \cA_{\beta,c}$ satisfies \eqref{e:mDlarge} with $a=2$ and $b=2-\beta$.
Hence we can take $\Omega=\G\bk G$ (or more precisely, a fundamental domain for $\G\bk G$) and $\mu_{\Omega}=\mu_{\G}$, and note that $\frac{a+1-b}{a+1}=\frac{1+\beta}{3}$ to get $D(\cL g,B)\leq m_{\cV}(B)^{\frac{1+\beta}{3}}\log(m_{\cV}(B))$ for $\mu_{\G}$-a.e. $\G g\in \G\bk G$.  While this exponent is not as good as the one in \thmref{t: Schmidt}, this argument also works for families that are not increasing as shown in the following.

\begin{Thm}\label{thm:nonestgen}
Let $\cB$ be a family of sets satisfying \eqref{e:mDlarge} for some $a>1$ and $0<b\leq 1$. Let $\{B_{t,r}\}_{t>1,0<r<1}\subseteq \cB$ be increasing in both parameters and assume we have the measure estimates $m_{\cV}(B_{t,r})\asymp t r$  and
\begin{equation}\label{r:mestimate}
 m_{\cV}\left(B_{t_1,r_1}\setminus B_{t_2,r_2}\right)\ll t_1(r_1-r_2)+(t_1-t_2)r_1,\; \forall\ 1< t_2< t_1,\ 0<r_2< r_1< 1.
\end{equation}
We then have the following.
\begin{enumerate}
\item  For any $0<\lambda<1$ let $r_t=t^{-\lambda}$. Then for $\mu_{\Omega}$-a.e. $g\in \Omega$ and all sufficiently large $t$ 
$$\#(\cL g\cap B_{t,r_t})=m_{\cV}(B_{t,r_t})+
O\left(m_{\cV}(B_{t,r_t})^{1-\frac{b}{a+1}}\log\left( m_{\cV}(B_{t,r_t})\right)\right).$$
\item For any $\lambda>\frac{1}{b}$, let $r_t=\tfrac{(\log t)^{\lambda}}{t}$. Then for $\mu_{\Omega}$-a.e. $g\in \Omega$ and all sufficiently large $t$ 
$$\#(\cL g\cap B_{t,r_t})=m_{\cV}(B_{t,r_t})+
O\left(m_{\cV}(B_{t,r_t})^{1-\frac{\lambda b-1}{\lambda(a+1)}}\log( m_{\cV}(B_{t,r_t}))\right).$$
\end{enumerate}
\end{Thm} 
\begin{proof}
For the first part, let $\nu=\frac{b}{a+1},\; \alpha=\frac{1}{(1-\lambda)\nu}$, $t_k=\max\{k^\alpha\log(k),2\}$ and consider the sequences $\{B_k^\pm\}$ with 
$B_k^+=B_{t_{k+1},r_{t_k}}$ and $B_k^-=B_{t_k, r_{t_{k+1}}}$. Then for any $t> t_1$ we have $B_{t,r_t}$ satisfies $B^-_{k}\subseteq B_{t,r_t}\subseteq B_k^+$ (with $k$ such that $t_k\leq t<t_{k+1}$) and we can estimate 
$$m_\cV(B_k^\pm)\gg t_{k}t_{k+1}^{-\lambda}\gg k^{\alpha(1-\lambda)}\log(k)^{1-\lambda}\gg k^{\alpha(1-\lambda)},$$
so that 
$$\sum_k m_{\cV}(B_k^\pm)^{a\nu-b}\log(m_\cV(B_k^\pm))^{-a} \ll \sum_k k^{-1}\log(k)^{-a}<\infty.$$
Moreover we can apply \eqref{r:mestimate} to get 
$m_{\cV}(B_k^+\setminus B_k^-)\ll m_{\cV}(B_k)^{\frac{a+1-b}{a+1}}\log(m_{\cV}(B_k))$,
so the result follows from \lemref{l:sequences}.

For the second part let $\nu=\frac{\lambda b-1}{\lambda(1+a)}$ and $\alpha=\frac{1+a}{\lambda b+a}$. First note that since $\nu=\frac{\lambda b-1}{\lambda(1+a)}$ and $\lambda>1/b$ we have $\alpha=\frac{1+a}{\lambda b+a}\in (0,1)$. Now we again take $B_k^+=B_{t_{k+1},r_{t_k}}$ and $B_k^-=B_{t_k, r_{t_{k+1}}}$ but with the sequence $t_k=e^{k^\alpha}$. 
Since $\{B_{t,r}\}$ is increasing in both parameters and $r_t=\frac{(\log t)^{\lambda}}{t}$ is eventually strictly decreasing in $t$, we have for all $t$ sufficiently large, $B_k^-\subseteq B_{t, r_t}\subseteq B_k^+$ with $k$ such that $t_k\leq t<t_{k+1}$. With this choice we have 
$$
m_{\cV}(B_k^-)\asymp e^{k^{\alpha}-(k+1)^{\alpha}}(k+1)^{\lambda\alpha}\asymp k^{\lambda \alpha}.
$$ 
Here we used that $0<(x+1)^{\alpha}-x^{\alpha}<1$ for any $x>0$ and $0<\alpha<1$. Similarly, we also have $m_{\cV}(B_k^+)\asymp k^{\lambda\alpha}$. Thus
$$\sum_k  m_{\cV}(B_k^\pm)^{a\nu-b} \log(m_{\cV}(B_k^\pm))^{-a}\ll \sum_k  k^{-1} \log(k)^{-a}<\infty.$$
For the difference by \eqref{r:mestimate} 
we have for all $k$ sufficiently large
\begin{align*}
m_{\cV}(B_k^+\setminus B_{k}^-)&\ll e^{(k+1)^{\alpha}}\left(\tfrac{k^{\lambda\alpha}}{e^{k^{\alpha}}}-\tfrac{(k+1)^{\lambda\alpha}}{e^{(k+1)^{\alpha}}}\right)+\left(e^{(k+1)^{\alpha}}-e^{k^{\alpha}}\right)\tfrac{k^{\lambda\alpha}}{e^{k^{\alpha}}}\\
&=2e^{(k+1)^\alpha-k^\alpha}k^{\lambda\alpha}-(k+1)^{\lambda\alpha}-k^{\lambda\alpha}\\
&<2\left(1+O_{\alpha}(k^{\alpha-1})\right)k^{\lambda\alpha}-2k^{\lambda\alpha}\\
&\ll k^{\alpha-1+\lambda\alpha}\ll m_{\cV}(B_k^-)^{1-\nu}.
\end{align*}
Here for the last estimate we used the relation $\alpha-1+\lambda\alpha= \lambda\alpha(1-\nu)$ which is equivalent to $\alpha(1+\lambda\nu)=1$. The latter can be checked using $\nu=\frac{\lambda b-1}{\lambda(1+a)}$ and $\alpha=\frac{1+a}{\lambda b+a}$.
The result then follows from \lemref{l:sequences}.
\end{proof}

\subsection{Counting result for all translates}
In this subsection we give a different (and more direct) approach to get a counting result which holds for \textit{all} $\cL$-translates.
\comm{Such results will then be used in the proof of \thmref{thm:c1digen}}.
In order to deduce such a result we need to impose some regularity conditions on our sets.
Let $\cO=\{\cO_{\e}\}_{0<\e\leq 1}\subseteq G$ denote some fixed increasing family of compact identity neighborhoods. We will assume $\cO$ is invariant under inversion, that is, $\cO_{\e}^{-1}=\cO_{\e}$ for all $0<\e\leq 1$. We say that a family of sets $\cB$ is \textit{strongly $\cA_{\beta,c}$-well rounded with respect to $\cO$} 
if there exist some $C_0\geq 0$ and $\e_0\in (0,1)$ such that for any $\e\in (0,\e_0)$ and $B\in \cB$  there are sets $\underline{B}_{\e},\overline{B}_{\e}\in \cA_{\beta,c}$ satisfying    
\begin{align}\label{equ:regassump}
	\underline{B}_{\e}\subseteq \bigcap_{h\in\cO_{\e}}B h\subseteq \bigcup_{h\in\cO_{\e}}Bh\subseteq \overline{B}_{\e},\quad  m_{\cV}(\overline{B}_{\e}\setminus \underline{B}_{\e})\leq {C_0(\e m_{\cV}(B)+1)}\quad \mbox{and }\quad \overline{B}_{\e}\setminus \underline{B}_{\e}\in \cA_{\beta,c}.
\end{align}

\begin{Thm}\label{thm:cfamd}
Let $\cO=\{\cO_\e\}_{0<\e\leq 1}$ be an  increasing family of compact inversion-invariant identity neighborhoods of $G$ satisfying
\begin{equation}\label{e:OeReg}
\mu_{\G}(\cO_\e) \geq \delta \e^{d}\quad \mbox{ and }\quad \cO_{\e_1}\cO_{\e_2}\subseteq \cO_{\eta(\e_1+\e_2)}
\end{equation}
for some $\delta, d,\eta>0$. Let $\cB$ be a family that is strongly $\cA_{\beta,c}$-well rounded with respect to $\cO$. Then for all $g\in G$ and all $B\in \cB$ with $m_{\cV}(B)\geq \max\{\delta^{-\frac{1}{2-\beta}}, \delta^{\frac{1}{d+\beta}}\}$
\begin{align*}
 			\#(\cL g\cap B)=m_{\cV}(B)+ O_{d, \eta,g}\left( \delta^{-\frac{1}{d+2}}m_{\cV}(B)^{1-\frac{2-\beta}{d+2}}\right).
 		\end{align*}
\end{Thm}
\begin{rem}
Note that the condition $m_{\cV}(B)\geq \max\{\delta^{-\frac{1}{2-\beta}}, \delta^{\frac{1}{d+\beta}}\}$ implies that $\delta^{-\frac{1}{d+2}}m_{\cV}(B)^{1-\frac{2-\beta}{d+2}}\leq m_{\cV}(B)$, and without this condition our bound on the remainder may be larger than the main term. 
\end{rem}
\begin{proof}[Proof of \thmref{thm:cfamd}]
Fix $g\in G$. First note that since $\cO_1$ is compact, $g\cO_{1}$ (hence also $g\cO_{\e}$ for any $0<\e<1$) 
is contained \comm{in} some finite number (which may depend on $g$) of fundamental domains for $\G\bk G$. For any $B\in \cB$ let $\epsilon=\epsilon(B,\delta)$ to be chosen later and $\epsilon'=\frac{\epsilon}{{\max\{1, 2\eta\}}}$ with $\eta$ as in \eqref{e:OeReg}. We then bound 
\begin{align*}
\int_{\cO_{\e'}}|D(\cL g, (\overline{B}_{\epsilon}\setminus \underline{B}_{\epsilon})h)|^2\, \d\mu_{\G}(h)
&= \int_{\cO_{\e'}}|D(\cL gh, (\overline{B}_{\epsilon}\setminus \underline{B}_{\epsilon}))|^2\, \d\mu_{\G}(h)\\
&=  \int_{g\cO_{\e'}}|D(\cL h, (\overline{B}_{\epsilon}\setminus \underline{B}_{\epsilon}))|^2\, \d\mu_{\G}(h)\\
&\ll_g \int_{\G\bk G}|D(\cL h, (\overline{B}_{\epsilon}\setminus \underline{B}_{\epsilon}))|^2\, \d\mu_{\G}(h)\\
&\ll m_{\cV}(\overline{B}_{\epsilon}\setminus \underline{B}_{\epsilon})^\beta+m_{\cV}(\overline{B}_{\epsilon}\setminus \underline{B}_{\epsilon}) \ll_{C_0}(\e m_{\cV}(B))^{\beta},
\end{align*}
where in the first line we used that $D(\cL g, B h)=D(\cL gh^{-1},B)$ and made a change of variable $h\mapsto h^{-1}$ (using also that both the Haar measure and $\cO_{\e'}$ are inversion-invariant), 
in the second we made a change of variables $h\mapsto g^{-1}h$, and in the last estimate we applied the second condition in \eqref{equ:regassump} and restricted our choice of $\e$ to satisfy $\e m_{\cV}(B)\geq 1$.

From this inequality and the lower bound on $\mu_{\G}(\cO_{\e'})$ in \eqref{e:OeReg} we can deduce that  there is some $h_0\in \cO_{\e'}$ with 
$|D(\cL g, (\overline{B}_{\epsilon}\setminus \underline{B}_{\epsilon})h_0)|\ll_{\eta,g} m_{\cV}(B)^{\frac{\beta}{2}} \delta^{-\frac12}\epsilon^{\frac{\beta-d}{2}}$, and hence 
\begin{align}\label{equ:estondiff}
\#(\cL g\cap (\overline{B}_{\epsilon}\setminus \underline{B}_{\epsilon})h_0)\ll_{\eta,g}
 m_{\cV}(B)^{\frac{\beta}{2}} \delta^{-\frac12}\epsilon^{\frac{\beta-d}{2}}+\e m_{\cV}(B) \ll \e m_{\cV}(B).
 \end{align}
For the last inequality we restrict our choice of $\epsilon$ to satisfy 
\begin{align}\label{e:EpsRes1}
\epsilon^{\beta-d-{2}}\leq  \delta  m_{\cV}(B)^{2-\beta}.
\end{align}

Next note that, since any $h\in \cO_{\e'}$ can be written as $h=h_1h_0$ with $h_1\in \cO_\e$ (cf. \eqref{e:OeReg})\comm{,} we have
$\underline{B}_{\e}h_0\subseteq B h \subseteq \overline{B}_{\e}h_0$ (cf. \eqref{equ:regassump}). {On the other hand, since $\e'\leq \e$ (so that $h_0^{-1}\in \cO_{\e'}\subseteq \cO_{\e}$) we also have $\underline{B}_{\e}h_0\subseteq B\subseteq \overline{B}_{\e}h_0$.}  Noting also that $m_{\cV}(Bh)=m_{\cV}(B)$ we can bound
\begin{align*}
\left|D(\cL g, Bh)-D(\cL g, B)\right|&\leq |(\#\cL g\cap Bh)-(\#\cL g\cap B)|\\
&{\leq \#(\cL g\cap Bh\Delta B)}\leq \#(\cL g\cap (\overline{B}_{\e}\setminus\underline{B}_{\e})h_0).
\end{align*}
This, together with \eqref{equ:estondiff} implies that for any $h\in \cO_{\e'}$, 
$$D(\cL g,B h)=D(\cL g,B )+O_{d,\eta,g}\left(\epsilon m_{\cV}(B)\right).$$
Hence there is some constant $C_1>0$ (which may depend on $d, \eta$ and $g$) such that 
$$
D(\cL g,B h)\geq D(\cL g,B )-C_1\e m_{\cV}(B)
$$ 
for all $h\in \cO_{\e'}$. Now we take $\e=m_{\cV}(B)^{-\frac{2-\beta}{d+2}}\delta^{-\frac{1}{d+2}}$ so that $\e m_{\cV}(B)=m_{\cV}(B)^{\frac{\beta}{2}}\delta^{-\frac12}\e^{-\frac{d}{2}}=m_{\cV}(B)^{\frac{d+\beta}{d+2}}\delta^{-\frac{1}{d+2}}$ and one easily checks that 
the conditions $\e m_{\cV}(B)\geq 1$ and \eqref{e:EpsRes1} are satisfied whenever $m_{\cV}(B)\geq \max\{\delta^{-\frac{1}{2-\beta}}, \delta^{\frac{1}{d+\beta}}\}$. 
If $D(\cL g, B)-C_1\e m_{\cV}(B)\leq 0$, then we have 
$$
D(\cL g, B)\leq C_1 m_{\cV}(B)^{\frac{d+\beta}{d+2}}\delta^{-\frac{1}{d+2}}.
$$
Otherwise, $D(\cL g,B )-C_1\e m_{\cV}(B)>0$ and thus we have
\begin{align*}
\int_{\cO_{\e'}}|D(\cL g,B h)|^2\, \d\mu_{\G}(h)\geq |D(\cL g,B )-C_1 \e m_{\cV}(B) |^2  \mu_{\G}(\cO_{\e'})\gg_{\eta, d}  \left|D(\cL g,B )-C_1\e m_{\cV}(B)\right|^2  \delta\e^{d}.
\end{align*}
Since we also have
\begin{align*}
\int_{\cO_{\e'}}|D(\cL g,B h)|^2\, \d\mu_{\G}(h)\ll_g \int_{\G\bk G}|D(\cL h,B )|^2\, \d\mu_{\G}(h) \ll m_{\cV}(B)^\beta,
\end{align*} 
we get that there is some $C_2>0$ (which only depends on $d, \eta$ and $g$) with 
$$ \left|D(\cL g,B )-C_1\e m_{\cV}(B)  \right|^2\leq  C_2 m_{\cV}(B)^{\beta}\delta^{-1} \e^{-d}.$$
Since $D(\cL g,B)-C_1\e m_{\cV}(B)>0$, this implies
$$D(\cL g,B )\leq C_1\e m_{\cV}(B) +\sqrt{C_2}m_{\cV}(B)^{\frac{\beta}{2}}\delta^{-\frac12} \e^{-\frac{d}{2}}= (C_1+\sqrt{C_2})m_{\cV}(B)^{\frac{d+\beta}{d+2}}\delta^{-\frac{1}{d+2}}.$$
To conclude, in both cases $D(\cL g, B)\ll_{d, \eta, g} m_{\cV}(B)^{\frac{d+\beta}{d+2}}\delta^{-\frac{1}{d+2}}$ whenever $m_{\cV}(B)\geq \max\{\delta^{-\frac{1}{2-\beta}}, \delta^{\frac{1}{d+\beta}}\}$. 
\end{proof} 
\subsection{Subgroup translates}
In this subsection, we apply the general counting results obtained in section \ref{sec:congenr} to the special case concerning counting for generic subgroup translates.
Let $H$ be a closed subgroup of $G$ and fix a continuous section, $\iota_H: H\bk G\to G$, of the natural projection from $G$ to $H\bk G$.
We can then write any $g\in G$ of the form $g=h\iota(x)$ for a unique $h\in H$ and $x\in H\bk G$. For some applications one is interested in counting results in translates $\cL h\iota_H(x)$ that hold for \textit{all} $h\in H$ and almost all $x\in H\bk G$. To be more precise we need to specify what measure we are using on the quotient $H\bk G$. When $H$ is unimodular the natural choice is the unique (up to scaling) $G$-invariant measure on the quotient $H\bk G$.
When $H$ is not unimodular, there is no such measure and instead we will use a  \textit{quasi-invariant measure} on $H\bk G$ as introduced in \cite{ReiterStegeman2000}. Explicitly, let $\Delta_H$ be the modular function of $H$ and it is shown in 
\cite[Proposition 8.1.3 and Theorem 8.1.13]{ReiterStegeman2000}  that there exists a strictly positive continuous function $\lambda_H: G\to \R_{>0}$ satisfying the functional equation
\begin{align}\label{equ:funeqqame}
\lambda_H(hg)=\Delta_H(h)\lambda_H(g),\quad \forall\ g\in G,\ h\in H,
\end{align}
and a measure $\mu_{H\bk G}$ on $H\bk G$ satisfying 
\begin{align}\label{equ:Haardec}
\int_G f(g)\lambda_H(g)\, \d\mu_{\G}(g)=\int_{H\bk G}\left(\int_Hf(hg)\, \d\mu_H(h)\right)\, \d\mu_{H\bk G}(Hg),\quad \forall\ f\in C_c(G),
\end{align}
where $\mu_{\G}$ is the Haar measure of $G$ which locally agrees with the probability $G$-invariant measure on $\G\bk G$ {and $\mu_H$ is a right Haar measure of $H$}.
\begin{remark}\label{rmk:parbcase}
When $H\leq G$ is unimodular, then $\lambda_H$ can be taken to be a constant function and $\mu_{H\bk G}$ is a $G$-invariant measure mentioned above. 
For our later applications, we will consider the case when $G=\SO_{Q_n}^+(\R)$ and $H=P$ the parabolic subgroup fixing the line spanned by $\bm{e}_0=(-1,\bm0, 1)\in \cV_{Q_n}^+$. Note that in this case $P\bk G\cong S^n$ and we can take $\iota_P=\mathfrak{s}$ (so that $\iota_P(\bm{\alpha})=k_{\bm{\alpha}}$ for all $\bm{\alpha}\in S^n$). 
Moreover, it is not difficult to see that under the coordinates $g=pk_{\bm{\alpha}}\in G$ with $p=u_{\bm{x}}a_ym\in UAM=P$ and $\bm{\alpha}\in S^n$, the modular function is given by $\Delta_P(p)=y^n$ and we can take $\d\mu_P(p)=y^{-1}\, \d\bm{x}\d y\d\mu_M(m)$, $\lambda_H(g)={\omega_{\G}}y^n$ and $\mu_{P\bk G}=\sigma_n$ (cf. \eqref{equ:haar3}).
\end{remark}

Let $\cO=\{\cO_{\e}\}_{0<\e<1}\subseteq H$ be an increasing family of inversion-invariant identity neighborhoods of $H$. We say a family of finite-measure Borel subsets $\mathcal{B}\subseteq \cV$ is \textit{$\cA_{\beta,c}$-well rounded} with respect to $\cO$ if there exist $T_0>1$, $C_0>0$ and $\e_0\in (0,1)$ such that for any $0<\e\leq \e_0$ \comm{and $B\in \cB$ with $m_{\cV}(B)\geq T_0$}, there exist sets $\overline{B}_{\e},\underline{B}_{\e}\in \cA_{\beta,c}$ satisfying that 
\begin{align}\label{equ:regassumpsubgp}
	\underline{B}_{\e}\subseteq \bigcap_{h\in\cO_{\e}}B h\subseteq \bigcup_{h\in\cO_{\e}}Bh\subseteq \overline{B}_{\e}\quad \text{and}\quad m_{\cV}(\overline{B}_{\e}\setminus \underline{B}_{\e})\leq {C_0(\e m_{\cV}(B)+1)}.
\end{align}
\begin{remark}
	The main difference between the well roundedness \comm{property} defined here and the strong well roundedness \comm{property} defined in \eqref{equ:regassump} is that here we do not require the differences $\overline{B}_{\e}\setminus\underline{B}_{\e}$ to also belong to $\cA_{\beta,c}$.
\end{remark}

We will apply \lemref{l:sequences} to get counting results for subgroup translates. The key ingredient to our arguments is the following proposition which verifies \eqref{e:mDlarge} for well bounded families. 

\begin{Prop}\label{prop:lardiscbdsubg}
Let $H\leq G$ be a closed subgroup and let $\cO=\{\cO_{\e}\}_{0<\e\leq 1}$ be an increasing family of compact inversion-invariant identity neighborhoods of $H$ satisfying that
\begin{align}
\mu_H(\cO_{\e})\gg \e^d,\quad \forall\ 0<\e\leq 1
\end{align}\label{equ:coverpre1}
for some $d>0$. 
Fix $h_0\in H$ and $\Omega$ a compact subset of $H\bk G$. 
Let $\mathcal{B}$ be a family of Borel subsets of $\cV$ which is $\cA_{\beta,c}$-well rounded with respect to $\cO$ (with parameters $T_0, C_0, \e_0$). Then for any $B\in \cB$ with $m_{\cV}(B)>\max\left\{\frac{2}{\e_0},T_0\right\}$ and for any $4C_0<X<2C_0\e_0m_{\cV}(B)$ we have
\begin{align}\label{equ:meestsubg}
\mu_{H\bk G}\left(\left\{x\in \Omega: D(\cL \iota_H(x)^{-1}h_0^{-1}, B)>X\right\}\right)\ll_{h_0,\Omega}\frac{m_{\cV}(B)^{\beta+d}}{X^{2+d}},
\end{align}
where the implied constant may also depend on $\cO$ and the parameters $c, C_0$ and $d$.
\end{Prop}

\begin{proof}
Let $B\in\cB$ be as in this proposition and for any $4C_0<X<2C_0\e_0m_{\cV}(B)$ let us \comm{define}
\begin{align*}
\cM_{B,X}=\cM_{B,X}^{h_0,\Omega}:=\left\{x\in \Omega: D(\cL\iota_H(x)^{-1}h_0^{-1}, B)>X\right\}.
\end{align*}
Let us also denote by 
\begin{align*}
\cM_{B,X}^{\pm}:=\left\{x\in \Omega: \pm\left(\#(\cL\iota_H(x)^{-1}h_0^{-1}\cap B)-m_{\cV}(B)\right)>X\right\}
\end{align*}
so $\cM_{B,X}=\cM_{B,X}^+\bigcup \cM_{B,X}^-$. It thus suffices to show $\mu_{H\bk G}(\cM_{B,X}^{\pm})\ll_{h_0,\Omega} \frac{m_{\cV}(B)^{\beta+d}}{X^{2+d}}
$.


Now take $\e=\frac{X-2C_0}{2C_0m_{\cV}(B)}$ (so $X-{C_0(\e m_{\cV}(B)+1)}= \frac{X}{2}$). Note that since $4C_0<X<2C_0\e_0m_{\cV}(B)$ we have $\e\in (0,\e_0)$. Observe that if $x\in \cM_{B,X}^+$, i.e. $\#(\cL\iota_H(x)^{-1}h_0^{-1}\cap B)-m_{\cV}(B)>X$, then, by the assumption \eqref{equ:regassumpsubgp} for any $h\in \cO_{\e}$ 
\begin{align*}
\#(\cL\iota_H(x)^{-1}h_0^{-1}\cap \overline{B}_{\e}h)&\geq \#(\cL\iota_H(x)^{-1}h_0^{-1}\cap B)>m_{\cV}(B)+X\\
&\geq m_{\cV}(\overline{B}_{\e})-C_0(\e m_{\cV}(B)+1)+X=m_{\cV}(\overline{B}_{\e})+\frac{X}{2},
\end{align*}
implying that
\begin{align*}
D(\cL\iota_H(x)^{-1}h_0^{-1}h^{-1}, \overline{B}_{\e})\geq \#(\cL\iota_H(x)^{-1}h_0^{-1}h^{-1}\cap \overline{B}_{\e})-m_{\cV}(\overline{B}_{\e})>\frac{X}{2}.
\end{align*}
Hence 
\begin{align*}
\frac{X^2}{4}\mu_{H\bk G}(\cM_{B,X}^+)\mu_H(\cO_{\e})<\int_{\Omega}\int_{\cO_{\e}}D(\cL\iota_H(x)^{-1}h_0^{-1}h^{-1}, \overline{B}_{\e})^2\, \d\mu_H(h)\, \d\mu_{H\bk G}(x)=:I.
\end{align*}
On the other hand, using \eqref{equ:Haardec} and making a change of variable $hh_0\mapsto h$ we get
\begin{align*}
I&=\int_{\cS_{\e}} D(\cL\iota_H(x)^{-1}h^{-1}, \overline{B}_{\e})^2\lambda_H(h\iota_H(x))\, \d\mu_{\G} (h\iota_H(x))\\
&=\int_{\cS_{\e}} D(\cL g^{-1}, \overline{B}_{\e})^2\lambda_H(g)\, \d\mu_{\G} (g)=\int_{\cS_{\e}^{-1}} D(\cL g, \overline{B}_{\e})^2\lambda_H(g^{-1})\, \d\mu_{\G} (g).
\end{align*}
Here
$\cS_{\e}=\cS_{\e,h_0,\Omega}:=\left\{h\iota_H(x): h\in \cO_{\e}h_0,\ x\in\Omega\right\}\subseteq G$.
Since $\cO_{\e}h_0\subseteq \cO_1h_0$ and both $\cO_1h_0\subseteq H$ and $\Omega\subseteq H\bk G$ are compact, $\cS_{\e}^{-1}$ can be covered by finitely many (which may depend on $\cO_1h_0$ and $\Omega$) fundamental domains for $\G\bk G$. Moreover, since $\lambda_H$ is continuous and strictly positive on $G$,   we have $\lambda_H(g^{-1})\asymp_{h_0,\cO_1,\Omega} 1$ for all $g\in \cS_{\e}^{-1}$.
Thus
\begin{align*}
I\ll_{h_0,\cO_1,\Omega}\int_{\cS_{\e}^{-1}}D(\cL g, \overline{B}_{\e})^2\, \d\mu_{\G} (g)\ll_{\cO_1h_0,\Omega}\int_{\G\bk G}D(\cL g, \overline{B}_{\e})^2\, \d\mu_{\G} (g)\ll_{c,C_0} m_{\cV}(\overline{B}_{\e})^{\beta}\ll m_{\cV}(B)^{\beta}.
\end{align*}
Here for the last estimate we used that $m_{\cV}(\overline{B}_{\e})<(1+\e)m_{\cV}(B)<2m_{\cV}(B)$. Thus we have
\begin{align*}
\mu_{H\bk G}(\cM_{B,X}^+)\ll_{c, C_0, \cO_1, h_0,\Omega}\tfrac{m_{\cV}(B)^{\beta}}{X^2}\mu_H(\cO_{\e})^{-1}\ll \tfrac{m_{\cV}(B)^{\beta}}{X^2}\left(\tfrac{X-2C_0}{2C_0m_{\cV}(B)}\right)^{-d}\asymp_{C_0,d} \tfrac{m_{\cV}(B)^{\beta+d}}{X^{2+d}}.
\end{align*}
Similarly we also have $\mu_{H\bk G}(\cM_{B,X}^-)\ll
\frac{m_{\cV}(B)^{\beta+d}}{X^{2+d}}$, finishing the proof.
\end{proof}

\comm{Note that \eqref{equ:meestsubg} is a special case of \eqref{e:mDlarge} with $a=d+2$ and $b=2-\beta$. By repeating the proof of \thmref{thm:incrcounting} and \thmref{thm:nonestgen},  replacing \eqref{e:mDlarge}) by \eqref{equ:meestsubg}, we obtain the following two results (which will then be used in the proof of  \thmref{thm:khinttypegen} and \thmref{thm:countforgenericgen} respectively).}

\begin{Thm}\label{thm:countforall}
	Let $H$ and $\cO$ be as in \propref{prop:lardiscbdsubg}.
	Let $\cB$ be an increasing family of Borel subsets of $\cV$ which is $\cA_{\beta,c}$-well rounded with respect to $\cO$ and satisfies $\{m_{\cV}(B): B\in \cB\}\supseteq (V,\infty)$ for some $V>0$. 
Then for all $h\in H$, for $\mu_{H\bk G}$-a.e. $x\in H\bk G$ and for all $B\in \cB$ with $m_{\cV}(B)$ sufficiently large
	\begin{align*}
		\#(\cL \iota_H(x)^{-1}h^{-1}\cap B)=m_{\cV}(B)+O\left(m_{\cV}(B)^{1-\frac{2-\beta}{d+3}}\log (m_{\cV}(B))\right).
	\end{align*}
\end{Thm}
\begin{Thm}\label{thm:countnoninc}
Let $H$ and $\cO$ be as in \propref{prop:lardiscbdsubg}. Let $\cB$ be a family of sets which is $\cA_{\beta,c}$-well rounded with respect to $\cO$. Let $\{B_{t,r}\}_{t>1,r<1}\subseteq \cB$ be increasing in both parameters and satisfy the same measure assumptions as in \thmref{thm:nonestgen}.
\begin{enumerate}
\item Let $r_t=t^{-\lambda}$ for some $0<\lambda<1$. Then for any $h\in H$, for $\mu_{H\bk G}$-a.e. $x\in H\bk G$ and for all $t$ sufficiently large
\begin{align*}
\#(\cL \iota_H(x)^{-1}h^{-1}\cap B_{t, r_t})=m_{\cV}(B_{t, r_t})+O\left(m_{\cV}(B_{t, r_t})^{1-\frac{2-\beta}{d+3}}\log \left(m_{\cV}(B_{t,r_t})\right)\right).
\end{align*}
\item Let $r_t=\frac{(\log t)^{\lambda}}{t}$ for some $\lambda>\frac{1}{2-\beta}$. Then for any $h\in H$, for $\mu_{H\bk G}$-a.e. $x\in H\bk G$ and for all $t$ sufficiently large
\begin{align*}
\#(\cL \iota_H(x)^{-1}h^{-1}\cap B_{t, r_t})=m_{\cV}(B_{t, r_t})+O\left(m_{\cV}(B_{t, r_t})^{1-\frac{\lambda(2-\beta)-1}{\lambda(d+3)}}\log \left(m_{\cV}(B_{t,r_t})\right)\right).
\end{align*}
\end{enumerate}
\end{Thm}
\subsection{Distribution of values of generic functions}
For any fixed function $F:\cV\to \R^m$ and $g\in G$ we let $F_g(\bm{v}):=F(\bm{v}g)$ and consider the family of function $\{F_g(\bm{v}): g\in G\}$. We are interested in the distribution of values $F_g(\cL)$ for generic translates $g\in G$. Explicitly given a family of increasing sets $\cB=\{B_t\}_{t\geq 1}$ and shrinking targets $\{\Omega_t\}_{t>1}\subseteq \R^m$ we want an estimate for 
$$\#\{\bm{v}\in \cL: F_g(\bm{v})\in \Omega_t,\ \bm{v}\in B_t\}.$$
We assume here that all finite-measure Borel sets on $\cV$ are in  $\cA_{1,c}$ for some $c>0$. With this assumption together with some well roundedness assumption on the family of sets $\cB$ and volume estimates for the sets $\{\bm{v}\in B_t: F_g(\bm{v})\in \Omega\}$ we can get an asymptotic estimates that holds for almost all $g$. \comm{The following will be used in the proof of \thmref{thm:linearformintro} and \thmref{thm:valuedisquadform}.}

\begin{Thm}\label{t:GenericFunctions}
Let $\cO=\{ \cO_\e\}_{0<\e\leq 1}$ be a family of identity neighborhoods in $G$ satisfying \eqref{e:OeReg} and let $\cB=\{B_t\}_{t\geq 1}\subseteq \R^n$ be a family of increasing sets satisfying that 
\begin{align}\label{equ:stabsmper}
B_{t(1-\e)}\subseteq \bigcap_{h\in \cO_\e} B_t h\subseteq \bigcup_{h\in \cO_\e} B_th\subseteq B_{t(1+\e)}.
\end{align}
Let $F:\cV\to \R^m$ be a fixed function and assume there is a family, $\Upsilon$,  of finite-volume Borel sets in $\R^m$, satisfying the following measure estimates:  for any $\Omega\in \Upsilon$ and $g\in G$
\begin{equation}\label{e:VolEst}
m_\cV(F_g^{-1}(\Omega)\cap B_t)=c(g)\vol(\Omega)t^a(1+O_g(t^{-\delta})),
\end{equation}
for some uniform $a, \delta>0$, with $c(g)$ and the implied constant are uniformly bounded on compact sets.  Let $\{\Omega_t\}_{t>1}\subseteq \Upsilon$ be a collection of shrinking targets with $\vol(\Omega_t)=t^{-b}$ for some $0< b<a$. Then for any $0<\nu<\min\{\delta, \frac{{a-b}}{d+3}\}$, for $\mu_{\G}$-a.e. $g\in G$ there is $t_0=t_0(g)$ such that for all $t\geq t_0$
\begin{align}\label{equ:valdis}
	\#\{\bm{v}\in \cL: F_g(\bm{v})\in \Omega_t,\ \bm{v}\in B_t\}=c(g)\vol(\Omega_t) t^a(1+O(t^{-\nu})).
\end{align}
\end{Thm}
\begin{proof}
We fix a compact set $\cK\subseteq G$. We show that the result holds for $\mu_{\G}$-a.e. $g\in \cK$. (Since $\cK$ is arbitrary this will conclude the proof.)
We note that since $\cO$ satisfies \eqref{e:OeReg} then for any $0<\e<1$ there is a finite set $\cI_\e\subseteq \cK$ of order $\#\cI_\e\asymp_{\cK, \eta} \e^{-d}$ with 
$\cK \subseteq \bigcup_{h\in \cI_\e} \cO_\e h$ (see \cite[Lemma 2.1]{KelmerYu2020}).  

For any finite-measure Borel set $B\subseteq \cV$ and $X\geq 1$ we denote as before by 
$$\cM_{B,X}=\cM^{(\cK)}_{B,X}=\left\{g\in \cK: D( \cL g, B)\geq X\right\},$$
and note that for any measurable set $B\subseteq \cV$ we have 
$$\mu_{\G}(\cM_{B,X})\ll_\cK \frac{m_{\cV}(B)}{X^2}.$$

Fix $0<\nu<\min\{\delta, \frac{a-b}{d+3}\}$ and take $\nu<\nu_0<\min\{\delta, \frac{a-b}{d+3}\}$.  Let $t_k =k^\alpha$  for some $0<\alpha< \nu_0^{-1}$ to be specified later and let $\e_k=t_k^{-\nu_0}=k^{-\alpha\nu_0}$. Let $\cI_k=\cI_{\e_k}\subseteq \cK$ denote the corresponding covering set.
For any $g\in \cK$ we can write $g=hg_i$ for some $h\in O_{\e_k}$ and $g_i\in \cI_k$ implying that for any $t_k\leq t<t_{k+1}$
$$B_{t_k(1-\e_k)}g_i \subseteq B_{t}g\subseteq B_{t_{k+1}(1+\e_k)}g_i.$$
Consequently, we have
$$F^{-1}(\Omega_{t_{k+1}})\cap B_{t_k(1-\e_k)}g_i
\subseteq F^{-1}(\Omega_{t})\cap B_{t}g
\subseteq  F^{-1}(\Omega_{t_k})\cap  B_{t_{k+1}(1+\e_{k})}g_i,$$
for some $g_i\in \cI_k$. Let $A_{k,i}^\pm$ be given by $A_{k,i}^+= F^{-1}(\Omega_{t_k})\cap  B_{t_{k+1}(1+\e_{k})}g_i$ and $A_{k,i}^-=F^{-1}(\Omega_{t_{k+1}})\cap B_{t_k(1-\e_k)}g_i$.
Then for any $g\in \cK$ and $t_k\leq t<t_{k+1}$,
\begin{align*}
D( F^{-1}(\Omega_{t})\cap B_{t}g,\cL g)&\leq \max\left\{D( \cL g, A_{k,i}^-),D( \cL g, A_{k,i}^+)\right\}+m_{\cV}(A_{k,i}^+\setminus A_{k,i}^-),
\end{align*}
for some $g_i\in \cI_k$.
Using the measure estimate \eqref{e:VolEst}, with our choice of $t_k=k^\alpha$ and $\e_k=k^{-\alpha\nu_0}$ (and noting that $\alpha\nu_0< \min\{1, \alpha\delta\}$) we get that 
$$m_{\cV}(A_{k,i}^{\pm})=c(g_i)k^{\alpha(a-b)}+O_\cK(k^{\alpha(a-b-\nu_0)}),$$
implying that  
\begin{eqnarray*}
m_{\cV}(A_{k,i}^+\setminus A_{k,i}^-)\leq C k^{\alpha(a-b-\nu_0)} \end{eqnarray*}
for some $C$ depending only on $\cK$. Hence for any $X_k\geq 2Ck^{\alpha(a-b-\nu_0)}$ we have
$$
\bigcup_{t_k\leq t<t_{k+1}}\left\{ g\in \cK: D(\cL g, F^{-1}(\Omega_{t})\cap B_{t}g)\geq X_k\right\}\subseteq \bigcup_{g_i\in \cI_k} \left(\cM_{A_{k,i}^+, \frac{X_k}{2}}\cup \cM_{A_{k,i}^-,\frac{X_k}{2}}\right).$$
Hence taking $X_k= k^{\alpha(a-b-\nu)}$ we get that for all sufficiently large $k$ (noting that $\nu<\nu_0$)
\begin{align*}
a_k:&=\mu_{\G}\left(\bigcup_{t_k\leq t< t_{k+1}}\left\{ g\in \cK: D(\cL g, F^{-1}(\Omega_{t})\cap B_{t}g)\geq k^{\alpha(a-b-\nu)}\right\} \right)\\
&\ll  \sum_{g_i\in \cI_{k}} \tfrac{m_\cV(A_{k,i}^+)+m_\cV(A_{k,i}^-)}{X_k^2}\ll   k^{-\alpha(a-b-2\nu-d\nu_0)} < k^{-\alpha(a-b-(d+2)\nu_0)}.
\end{align*}
We choose the parameter $\alpha$ so that $\alpha (a-b-(d+2)\nu_0)>1$. Then $\alpha$ needs to satisfy two conditions that 
$\nu_0^{-1}> \alpha >\frac{1}{ a-b-(d+2)\nu_0}$, which can be satisfied since $\nu_0^{-1} >\frac{1}{ a-b-(d+2)\nu_0}$ (the latter is true since $\nu_0<\frac{{a-b}}{d+3}$). With this choice of $\alpha$ the series $\sum_ka_k$
converges, implying that for $\mu_{\G}$-a.e. $g\in \cK$ there is some $t_0=t_0(g)$ such that for all $t\geq t_0$ we have $D(\cL g, F^{-1}(\Omega_{t})\cap B_{t}g)< t^{a-b-\nu}$, 
which is the same as \eqref{equ:valdis}.
\end{proof}

\section{Random forms on the light cone}\label{sec:ranform}
In this section we prove \thmref{thm:linearformintro} and \thmref{thm:valuedisquadform} by applying \thmref{t:GenericFunctions} in the previous section. In fact for quadratic forms, we will prove a more general theorem concerning the value distribution of a certain family of homogeneous functions and \thmref{thm:valuedisquadform} is a special case of this theorem when the degree is $2$; see \thmref{thm:homogen} below. 

We follow the notation in section \ref{sec:conj}. Assume $n\geq \comm{2}$ 
and let $Q(\bm{v})=Q_n(\bm{v}\tau)$ be a rational $\Q$-isotropic quadratic form of signature $(n+1,1)$ with $\tau\in \GL_{n+2}(\R)$ as in \rmkref{rmk:formoftau}. Let $G=\SO_{Q_n}^+(\R)$ and $\G=\tau^{-1}\SO_Q^+(\Z)\tau< G$ and let $\tilde{\mu}_{Q}$ be the $G$-invariant probability measure on $\G\bk G$. Let $\tilde\cL_Q=\cL_Q\tau\subseteq \cV_{Q_n}^+$. Throughout this section we assume $c_Q=0$ so that $\beta_Q=1$ and the variance bound \eqref{equ:disbdwiexpotran} holds for any finite-measure set $B\subseteq \cV_{Q_n}^+$. 
We abbreviate $\omega_Qm_{\cV_{Q_n}^+}(B)$ by $|B|$. This variance bound will be one of the two main ingredients to study these value distribution problems. In the following we obtain the other main ingredient, namely, effective volume estimates for the functions under consideration. 

We start our discussion by a classification result due to Sargent \cite{Sargent2014}. %

\subsection{Classification of linear and quadratic forms}
We fix $m<n$ a positive integer and denote the standard projection by 
\begin{align}\label{equ:stanproj}
L_0: \R^{n+2}\to \R^m,\qquad \bm{v}L_0=(v_1,\ldots, v_m).
\end{align}

The following lemma is a special case of \cite[Lemma 2.2]{Sargent2014}.
\begin{Lem}\label{l:classification}
	For any linear map $L:\R^{n+2}\to \R^m$ of rank $m$ with $Q_{L}:=Q|_{\ker(L)}$ indefinite, there is $g_L\in \SO_Q^+(\R)$ and $h_L\in \GL_m(\R)$ such that 
	$L=g_L{\tau}L_0h_L$ with $L_0$ the projection as above.
\end{Lem}
\begin{remark}
	In fact, \cite[Lemma 2.2]{Sargent2014} implies that there exist $g_L\in \text{O}_{Q}(\R)$ and $h_L\in \GL_m(\R)$ such that $L=g_L \tau L_0h_L$. Replacing $g_L$ by $g_Lg_0$ with $g_0\in \text{O}_Q(\R)$ of the form  $g_0=\tau\left(\begin{smallmatrix} I_n & &\\ &\pm 1 &\\ &&\pm1\end{smallmatrix}\right)\tau^{-1}$ with the signs \comm{of the last two diagonal entries} properly chosen and noting that $g_0\tau L_0=\tau L_0$ (since $m<n$) we may take $g_L$ from $\SO_Q^+(\R)$ as stated above.
\end{remark}

Moving from a linear form to a quadratic form $F: \R^{n+2}\to \R$ of rank $m$ we can use the following classification. First note since \comm{$F$} is of rank $m$, it has a null subspace $V_{\rm{null}}(F)\subseteq \R^{n+2}$ of dimension $n+2-m$. Also recall that $F$ is called \textit{indefinite on $\cV^+_Q$} if both $F$ and $Q|_{V_{\rm{null}}(F)}$ are indefinite. Combining the classification of linear forms on $\cV^+_Q$ with the classification of quadratic forms on $\R^m$ we get the following as a direct consequence.
\begin{Lem}\label{l:QuadClas}
	\comm{Assume $m\geq 2$.} Let $F$ be a quadratic form of rank $m$ on $\R^{n+2}$ that is indefinite on $\cV^+_Q$. Then there exist integers $p\geq q\geq 1$ with $p+q=m$ and $g\in \SO^+_Q(\R)$ and $h\in \GL_m(\R)$ satisfying 
	$F(\bm{v})=F_{p,q}(\bm{v}g{\tau} L_0h)$, where $L_0$ is the projection as \comm{in \eqref{equ:stanproj}} and $F_{p,q}(\bm{w})=\sum_{j=1}^p w_j^2-\sum_{j=p+1}^m w_j^2$.
\end{Lem}
We now describe the family of homogenous functions that we will be working with. Fix a positive real $1< d<m$ and let $\mathscr{F}_{d,Q}$ be the family consisting of functions of the form
$F:\R^{n+2}\to \R$ with $F(\bm{v})=F^{(d)}_{p,q}(\bm{v}g{\tau}L_0h)$ for some integers $p\geq q\geq 1$ with $p+q=m$, $g\in \SO_Q^+(\R)$ and $h\in \GL_m(\R)$. Here for any such a pair $(p,q)$,
\begin{align*}
	F_{p,q}^{(d)}(\bm{w}):=\sum_{j=1}^p|w_j|^d-\sum_{j=p+1}^m|w_j|^d.
\end{align*}
We note that in view of \lemref{l:QuadClas}, $\mathscr{F}_{2,Q}$ is the space of rank $m$ quadratic forms on $\R^{n+2}$ that are indefinite on $\cV_Q^+$. We now state our counting results for linear forms and these homogeneous functions.

\begin{Thm}\label{thm:linearform}
Let  $m<n$ be a positive integer and let $\{\Omega_T\}_{T>1}\subseteq \R^m$ be a decreasing family of bounded measurable sets with $\vol(\Omega_T)=T^{-a}$ for some $0<a<n-m$. Then there is some $\nu>0$ such that for almost every linear map $L:\R^{n+2}\to \R^m$ with $Q_n|_{\ker{L}}$ indefinite we have that for all sufficiently large $T$
$$
\#\{\bm{w}\in \tilde\cL_Q: \|\bm{w}\|\leq T,\; \bm{w}L\in \Omega_T\}=c(Q,L) T^{n-m-a}(1+O(T^{-\nu})),
$$
where $c(Q,L)$ is some positive constant depending only on $Q$ and $L$.
\end{Thm}

\begin{Thm}\label{thm:homogen}
	\comm{Assume $n>m\geq 2$}. Fix a positive real $1<d<m$. Let $\{I_T\}_{T>1}\subseteq \R$ be a decreasing family of bounded measurable sets with $|I_T|=T^{-a}$ for some $0<a<n-d$. Then  there is some $\nu>0$ such that for almost every $F\in \mathscr{F}_{d,Q_n}$ and for all sufficiently large $T$ 
	\begin{align*}
		\#\{\bm{w}\in \tilde\cL_Q: \|\bm{w}\|\leq T,\ F(\bm{w})\in I_T\}=c(Q,F)T^{n-d-a}(1+O(T^{-\nu})),
	\end{align*}
	where $c(Q,F)>0$ is some constant depending only on $Q$ and $F$.
\end{Thm}

\begin{remark}\label{rmk:almostall}
	Let us now make the ``almost every'' notion in \thmref{thm:linearform} and \thmref{thm:homogen} more precise. We say that a result hods for almost every linear map $L: \R^{n+2}\to \R^m$ of rank $m$ with $Q_L$ indefinite if it holds for 
	$L=g\tau L_0h$ for all $h\in \GL_m(\R)$ and for almost every $g\in \SO^+_Q(\R)$ (with respect to a Haar measure of $\SO^+_Q(\R)$). Similarly a result holds for almost every $F\in \mathscr{F}_{d,Q}$ if it holds for $F(\bm{v})=F^{(d)}_{p,q}(\bm{v}g\tau L_0h)$ for all pairs $(p,q)$ with $p\geq q\geq 1$ and $p+q=m$, for
	all $h\in \GL_{m}(\R)$ and for almost every $g\in \SO^+_Q(\R)$.
\end{remark}
\begin{proof}[Proof of \thmref{thm:linearformintro} \comm{assuming \thmref{thm:linearform}}]
\thmref{thm:linearform} implies that for any $h\in \GL_m(\R)$, for almost every $g'\in G$ and for all sufficiently large $T$,
\begin{align*}
\#\left\{\bm{w}\in \tilde\cL_Q: \|\bm{w}\|\leq T,\; \bm{w}g'L_0h\in \Omega_T\right\}=c(Q,gL_0h) T^{n-m-a}(1+O(T^{-\nu})).
\end{align*}
Recall that $\tilde\cL_Q=\cL_Q\tau\subseteq \cV_{Q_n}^+$ with $\tau\in \GL_{n+2}(\R)$ as in \rmkref{rmk:formoftau} such that $Q(\bm{v})=Q_n(\bm{v}\tau)$. For $\bm{w}\in \tilde\cL_Q$, letting $\bm{v}=\bm{w}\tau^{-1}$ and recalling that $\|\bm{v}\|_Q=\|\bm{v}\tau\|=\|\bm{w}\|$ we have
\begin{align*}
\#\left\{\bm{v}\in \cL_Q: \|\bm{v}\|_Q\leq T,\; \bm{v}g\tau L_0h\in \Omega_T\right\}=\#\left\{\bm{w}\in \tilde\cL_Q: \|\bm{w}\|\leq T,\; \bm{w}g'L_0h\in \Omega_T\right\},
\end{align*}
where $g=\tau g'\tau^{-1}\in \SO_Q^+(\R)$. We can then finish the proof by noting that $Q_n|_{\ker(g'L_0h)}$ is indefinite if and only if $Q|_{\ker(g\tau L_0h)}$ is indefinite and that the $\tau$-conjugation map sends a full measure set in $G$ to a full measure set in $\SO_Q^+(\R)$.
\end{proof}
\thmref{thm:valuedisquadform} follows similarly from \thmref{thm:homogen} and we omit it here. For the remaining we prove \thmref{thm:linearform} and \thmref{thm:homogen} by proving necessary volume estimates.
\subsection{Volume calculation for linear \comm{maps}}
We first consider the case of a linear map and \comm{prove} the following volume estimate that is needed in order to apply \thmref{t:GenericFunctions}. In the following all balls are assumed to be centered at the origin and all implied constants may depend on the fixed parameters $n,m$ and $d$. \comm{Let $L:\R^{n+2}\to \R^m$ be a linear map of rank $m<n$ with $Q_{n,L}:=Q_n|_{\ker(L)}$ indefinite. We need to introduce a relevant measure on $\cV_{Q_{n,L}}^+:=\cV_{Q_n}^+\cap \ker(L)$. When $L=L_0$, $\ker(L_0)$ naturally identifies with $\R^{n+2-m}$ and with this identification $Q_{n,L_0}=Q_{n-m}$ and $\cV_{Q_{n,L_0}}^+=\cV_{Q_{n-m}}^+$. We denote by $m_{\cV_{Q_{n,L_0}}^+}$ the $\SO_{Q_{n,L_0}}^+(\R)$-invariant measure on $\cV_{Q_{n,L_0}}^+$ given as in \eqref{equ:lebmea} with $n-m$ in place of $n$. More generally, by \lemref{l:classification}, there exist $g_L\in G$ and $h_L\in \GL_m(\R)$ such that $L=g_LL_0h_L$. Then $\ker(L)=\ker(L_0)g_L^{-1}$ and $\cV_{Q_{n,L}}^+=\cV_{Q_{n,L_0}}^+g_L^{-1}$. We define $m_{\cV_{Q_{n,L}}^+}$ to be the pushforward of $m_{\cV_{Q_{n,L_0}}^+}$ under the right $g_L^{-1}$-multiplication map.}
\begin{Prop}\label{p:volLinOmega}
	For any linear map $L=g_LL_0h_L:\R^{n+2}\to \R^m$ of rank $m<n$ with $Q_{n,L}$ indefinite, for any measurable $\Omega\subseteq \R^m$ contained in a ball of radius $c_\Omega$ and for any $T>2c_{\Omega}\|h_L\|_{\rm op}$ we have 
	\begin{align*} 
		m_{\cV^+_{Q_n}}\left(\left\{\bm{v}\in \cV^+_{Q_n}: \|\bm{v}\|\leq T,\; \bm{v}L\in \Omega\right\}\right)&= {\frac{c_{n,m}\vol(\Omega )T^{n-m}}{|\det(h_L)|}} \left(V_L+O\left(c_\Omega^{\frac12} \|h_L\|_{\rm op}^{\frac12} \|g_L\|_{\rm op}^{n-m+\frac12}T^{-\frac12} \right)\right).
	\end{align*}
	Here 
	$c_{n,m}:=\frac{(n-m+1)\G(\frac{n+3}{2})}{(n+1)\pi^{\frac{m}{2}}\G(\frac{n-m+3}{2})}$ with $\G(s)$ the Gamma function and 
	$
	V_L:=m_{\cV^+_{Q_{n,L}}}\left(\{\bm{v}\in \cV^+_{Q_{n,L}}: \|\bm{v} \|\leq 1\}\right)
	$.
\end{Prop}
\begin{rem}
	The assumption that $Q_{n,L}$ is indefinite is necessary here. For example, for $\bm{v}L=(v_{n-m+3},\ldots, v_{n+2})$ for which $Q_{n,L}$ is definite, one can see 
	$$
	\left\{\bm{v}\in \cV^+_{Q_n}: \|\bm{v}\|\leq T,\; \bm{v}L\in \Omega\right\}\subseteq \left\{ \bm{v}\in \cV^+_{Q_n}: \|\bm{v}\|\leq \sqrt{2}c_\Omega\right\}
	$$ is uniformly bounded and independent on $T$.
\end{rem}

We first prove a smoothed version of this volume estimate.

\begin{Lem}\label{l:LinOmega}
	For any $f\in C^1_c(\R^{n+2})$ supported in a ball of radius $c_f\geq 1$ and for any $\Omega\subseteq \R^m$ contained in the unit ball we have 
	for all $T\geq 2$
	\begin{align*}
		\int_{\cV^+_{Q_n}} f\left(\tfrac{\bm{v}}{T}\right) \chi_{\Omega}(\bm{v}L_0)\, \d m_{\cV^+_{Q_n}}(\bm{v})&= c_{n,m}\vol(\Omega)\left(T^{n-m} \int_{\cV^+_{Q_{n,L_0}} }f\, \d m_{\cV^+_{Q_{n,L_0}}}+O\left( \cS_{1}(f)
		c_f^{n-m} T^{n-m-1}\right)\right), 
	\end{align*}
	{where $c_{n,m}$ is as \comm{in \propref{p:volLinOmega}} and 
	$\cS_{1}(f):=\max\left\{\|f\|_{\infty},\ \|\frac{\partial f}{\partial v_j}\|_{\infty},\ 1\leq j\leq n+2\right\}$}.
\end{Lem}
\begin{proof}
	Let 
	$\cI:=\int_{\cV^+_{Q_n}} f\left(\tfrac{\bm{v}}{T}\right) \chi_{\Omega}(\bm{v}L_0)\, \d m_{\cV^+_{Q_n}}(\bm{v})$.
	We can write any $\bm{v}\in \cV^+_{Q_n}$ as $\bm{v}=r(\bm{\alpha},1)$ with $r>0$ and $\bm{\alpha}\in S^n$ and in these coordinates, by \eqref{equ:lebmea} we have $\d m_{\cV^+_{Q_n}}(\bm{v})=r^{n-1}\, \d r\d\sigma_n(\bm{\alpha})$, where $\sigma_n$ is the rotation invariant probability measure on $S^n$ as before. We thus get
	$$\cI=\int_0^\infty\int_{S^n} f\left(r\tfrac{(\bm{\alpha},1)}{T}\right) \chi_{\Omega}(r(\bm{\alpha},1)L_0)r^{n-1}\, \d\sigma_n(\bm{\alpha}) \d r.$$
	Writing Lebesgue measure on $\R^{n+1}$ in polar coordinates $\d\bm{x}={(n+1)\nu_{n+1}}t^n\, \d t\d\sigma_n(\bm{\alpha})$ {with $\nu_{n+1}=\frac{\pi^{\frac{n+1}{2}}}{\G(\frac{n+3}{2})}$ the volume of the unit ball $B^{n+1}$ in $\R^{n+1}$}, 
	we can replace the integral over $S^n$ by an integral over $B^{n+1}$, and make a change of variables $r\mapsto rT$ to get
	$$\cI={\frac{T^n}{\nu_{n+1}}}\int_0^\infty\int_{B^{n+1}} f\left(r(\tfrac{\bm{x}}{\|\bm{x}\|},1)\right) \chi_{\Omega}\left(rT(\tfrac{\bm{x}}{\|\bm{x}\|},1)L_0\right)r^{n-1}\, \d\bm{x}\d r.$$
	Now write $\bm{x}=(\bm{y},\bm{z})$ with $\bm{y}\in \R^m$ and $\bm{z}\in \R^{n+1-m}$ and make a change of variables $\bm{y}\mapsto \|\bm{z}\|\bm{y}$ to get that 
	\begin{align*}
		\cI&=\frac{T^n}{\nu_{n+1}}\int_0^\infty\int_{\|\bm{y}\|^2+\|\bm{z}\|^2\leq 1} f\left(r(\tfrac{\bm{y}}{\sqrt{\|\bm{y}\|^2+\|\bm{z}\|^2}},\tfrac{\bm{z}}{\sqrt{\|\bm{y}\|^2+\|\bm{z}\|^2}},1)\right) \chi_{\Omega}\left( \tfrac{rT\bm{y}}{\sqrt{\|\bm{y}\|^2+\|\bm{z}\|^2}} \right)r^{n-1}\, \d\bm{y}\d\bm{z}dr\\
		&=\frac{T^n}{\nu_{n+1}}\int_0^\infty\int_{\|\bm{z}\|^2(\|\bm{y}\|^2+1)\leq 1} f\left(r(\tfrac{\bm{y}}{\sqrt{\|\bm{y}\|^2+1}},\tfrac{\bm{z}}{\|\bm{z}\| \sqrt{\|\bm{y}\|^2+1}},1)\right) \chi_{\Omega}\left( \tfrac{rT\bm{y}}{\sqrt{\|\bm{y}\|^2+1}} \right)r^{n-1} \|\bm{z}\|^m\, \d\bm{y}\d\bm{z}\d r.
	\end{align*}
	Now make another change of variables $\bm{u}=\frac{\bm{y}}{\sqrt{1+\|\bm{y}\|^2}}$ so that $\bm{y}=\frac{\bm{u}}{\sqrt{1-\|\bm{u}\|^2}}$ and note that the Jacobian of this transformation is given by $J(\bm{u})=(1-\|\bm{u}\|^2)^{-\frac{2+m}{2}}$.
	With this change of variables the region $ \|\bm{z}\|^2(\|\bm{y}\|^2+1)\leq 1$ transforms into $\|\bm{z}\|^2+\|\bm{u}\|^2\leq 1$ so that 
	\begin{align*}
		\cI&=\frac{T^n}{\nu_{n+1}}\int_0^\infty\int_{\|\bm{u}\|< 1}\int_{\|\bm{z}\|^2\leq 1-\|\bm{u}\|^2} f\left(r(\bm{u},\tfrac{\bm{z}\sqrt{1-\|\bm{u}\|^2}}{\|\bm{z}\|},1)\right) \chi_{\Omega}(rT \bm{u} )r^{n-1} \|\bm{z}\|^m J(\bm{u})\, \d\bm{z}\d\bm{u}\d r.
	\end{align*}
	Writing $\bm{z}=\zeta \bm{\omega}$ with $\zeta=\|\bm{z}\|$ and $\bm{\omega}\in S^{n-m}$ in polar coordinates and noting that $c_{n,m}=\frac{(n-m+1)\nu_{n-m+1}}{(n+1)\nu_{n+1}}$ 
	we get 
	\begin{align*}
		\cI&={c_{n,m}}T^n\int_0^\infty\int_{\|\bm{u}\|< 1}\int_{S^{n-m}}(1-\|\bm{u}\|^2)^{\frac{n-m-1}{2}} f\left(r(\bm{u},\bm{\omega} \sqrt{1-\|\bm{u}\|^2},1)\right) \chi_{\Omega}(rT \bm{u} )r^{n-1}\,\d\sigma_{n-m}(\bm{\omega}) \d\bm{u}\d r.
	\end{align*}
		Let $\cI^-$ (respectively $\cI^+$) be the above triple integral except with the integration range of $r$ given by $0\leq r\leq \frac{1}{T}$ (respectively $r>\frac{1}{T}$) so that $\cI=\cI^-+\cI^+$. We first estimate $\cI^-$ trivially by
	\begin{align*}
		\cI^-
		&\ll \|f\|_{\infty}T^n\int_0^{T^{-1}}(rT)^{-m}\vol(\Omega)r^{n-1}\d r\ll \|f\|_\infty \vol(\Omega).
		&
	\end{align*}
	Next to estimate $\cI^+$, i.e. for $r$ in the range $r> \tfrac{1}{T}$, we may assume $\|\bm{u}\|\leq \frac{1}{rT}< 1$ (since $rT\bm{u}\in \Omega$) and bound
	$$\left\|r(\bm{u},\bm{\omega} \sqrt{1-\|\bm{u}\|^2},1)-r(\bm0,\bm{\omega},1)\right\|^2\leq 2r^2\|\bm{u}\|^2 \leq 2T^{-2},$$ 
	so that 
	$$f\left(r(\bm{u},\bm{\omega} \sqrt{1-\|\bm{u}\|^2},1)\right)=f\left(r(\bm0,\bm{\omega}, 1)\right)+O\left(\tfrac{\cS_{1}(f)}{ T}\right).$$
	Using this approximation we get
	\begin{align*}
		\cI^+
		&=c_{n,m}T^n\int_{1/T}^\infty\int_{(rT)^{-1}\Omega}\int_{S^{n-m}}(1-\|\bm{u}\|^2)^{\frac{n-m-1}{2}} f\left(r(\bm0,\bm{\omega}, 1)\right) r^{n-1}\, \d\sigma_{n-m}(\bm{\omega})\d\bm{u}\d r\\
		&+O\left( \cS_{1}(f)\vol(\Omega) c_f^{n-m} T^{n-m-1}\right).
	\end{align*}
	Next approximating $(1-\|\bm{u}\|^2)^{\frac{n-m-1}{2}}=1+O((\frac{1}{rT})^2)$ uniformly for $\|\bm{u}\|\leq \frac{1}{rT}$, we can replace $(1-\|\bm{u}\|^2)^{\frac{n-m-1}{2}}$ by $1$ after adding an error of 
	$O\left(\|f\|_\infty  \vol(\Omega)T^{n-m-2}\int_{1/T}^{c_f} r^{n-m-3}\, \d r\right)$,
	which can be absorbed into the previous error term, to get that 
	%
	\begin{align*}
		\cI^+
		&=c_{n,m}T^{n-m}\vol(\Omega)\int_{1/T}^\infty\int_{S^{n-m}} f\left(r(\bm0,\bm{\omega}, 1)\right) r^{n-m-1}\, \d\sigma_{n-m}(\bm{\omega}) \d r\\
		&+O\left( \cS_{1}(f)\vol(\Omega) c_f^{n-m} T^{n-m-1}\right).
	\end{align*}
	We now add back the contribution of $0\leq r\leq \frac{1}{T}$ and note that it is bounded by $O(\|f\|_\infty \vol(\Omega))$
	to get
	\begin{align*}
		\cI
		&=c_{n,m}T^{n-m}\vol(\Omega)\int_{0}^\infty\int_{S^{n-m}} f\left(r(\bm0,\bm{\omega}, 1)\right) r^{n-m-1}\, \d\sigma_{n-m}(\bm{\omega}) \d r\\
		&+O\left( \cS_{1}(f)\vol(\Omega) c_f^{n-m} T^{n-m-1}\right).
	\end{align*}
	Finally, we conclude the proof by noting that 
		\begin{displaymath}
		\int_{0}^\infty\int_{S^{n-m}} f\left(r(\bm0,\bm{\omega}, 1)\right) r^{n-m-1}\, \d\sigma_{n-m}(\bm{\omega})\d r=\int_{\cV^+_{Q_{n,L_0}} }f\, \d m_{\cV^+_{Q_{n,L_0}}}.\qedhere
	\end{displaymath}
\end{proof}

\begin{proof}[Proof of \propref{p:volLinOmega}]
	\comm{Abbreviate $g=g_L\in G$ and $h=h_L\in \GL_m(\R)$.}
	Note that 
	$$\left\{\bm{v}\in \cV^+_{Q_n}: \|\bm{v}\|\leq T,\; \bm{v}L\in \Omega\right\}=\left\{\bm{v}\in \cV^+_{Q_n}: \|\bm{v}g^{-1}\|\leq T,\; \bm{v}L_0\in \Omega h^{-1}\right\}g^{-1}.$$
	Let $\Omega'=\Omega h^{-1}$, and since
	$m_{\cV^+_{Q_n}}$ is \comm{$G$}-invariant, it is enough to estimate 
	$$R:=m_{\cV^+_{Q_n}}\left(\left\{\bm{v}\in \cV^+_{Q_n}: \|\bm{v} g^{-1}\|\leq T,\;  \bm{v}L_0 \in \Omega'\right\}\right).$$
	By a simple rescaling argument and the estimate $c_{\Omega'}\leq c_{\Omega}\|h\|_{\rm op}$, it is enough to consider the case when $\Omega'$ is contained in the unit ball and $T\geq 2$. 
	For $0<\delta<1$ let $\psi^\pm_\delta$ approximate the indicator function of the unit ball in $\R^{n+2}$ in the sense that 
	$$\psi_\delta^+(\bm{v})=\left\lbrace\begin{array}{ll} 1& \|\bm{v}\|\leq 1,\\ 0 & \|\bm{v}\|\geq 1+\delta,\end{array}\right.$$ 
	and similarly 
	$$\psi_\delta^-(\bm{v})=\left\lbrace\begin{array}{ll} 1& \|\bm{v}\|\leq 1-\delta, \\ 0 & \|\bm{v}\|\geq 1,\end{array}\right.$$ 
	with 
	$\cS_{1}(\psi_\delta^\pm)\asymp \delta^{-1}$.  For any $g\in \comm{G}$ we define 
	$\psi_{\delta,g}^\pm(\bm{v})=\psi_{\delta}^\pm (\bm{v} g^{-1})$ and note that 
	$$ \int_{\cV^+_{Q_n}} \psi^-_{\delta,g}\left(\tfrac{\bm{v}}{T}\right)\chi_{\Omega'}(\bm{v}L_0)\, \d m_{\cV^+_{Q_n}}(\bm{v})\leq R\leq \int_{\cV^+_{Q_n}} \psi^+_{\delta,g}\left(\tfrac{\bm{v}}{T}\right)\chi_{\Omega'}(\bm{v}L_0)\, \d m_{\cV^+_{Q_n}}(\bm{v}).$$
	Applying \lemref{l:LinOmega} with $f=\psi^\pm_{\delta,g}$, noting that $\cS_{1}(f)\ll \|g\|_{\rm op}\delta^{-1}$ and $c_f\ll \|g\|_{\rm op}$ we get 
	\begin{align*}
		\int_{\cV^+_{Q_n}} \psi^\pm_{\delta,g}\left(\tfrac{\bm{v}}{T}\right)\chi_{\Omega'}(\bm{v}L)\, \d m_{\cV^+_{Q_n}}(\bm{v})&= {c_{n,m}}\vol(\Omega')T^{n-m}\left(\int_{\cV^+_{Q_{n,L_0}}}\psi^\pm_{\delta,g}\, \d m_{\cV^+_{Q_{n,L_0}}}
		+O\left( \|g\|_{\rm op}^{n-m+1} \delta^{-1} T^{-1}\right)\right).
	\end{align*}
	Further estimate (see \cite[p. 8686-8687]{KelmerYu2020})
	$$\int_{\cV^+_{Q_{n,L_0}}}\psi^\pm_{\delta,g}\, \d m_{\cV^+_{Q_{n,L_0}}}=m_{\cV^+_{Q_{n,L_0}}}\left(\left\{\bm{v}\in \cV^+_{Q_{n,L_0}}: \|\bm{v} g^{-1}\|\leq 1\right\}\right)+O\left(\|g\|_{\rm op}^{n-m}\delta\right),$$
	to get 
	\begin{align*}
		R&=c_{n,m}\vol(\Omega') T^{n-m}\left(m_{\cV^+_{Q_{n,L_0}}}\left(\{\bm{v}\in \cV^+_{Q_{n,L_0}}: \|\bm{v} g^{-1}\|\leq 1\}\right)+O\left(\|g\|_{\rm op}^{n-m}\left(\delta+ \|g\|_{\rm op}\delta^{-1}T^{-1}\right)\right)\right).
	\end{align*}
	Now setting $\delta=T^{-1/2}$ we get
	\begin{align*}
		R&= c_{n,m}\vol(\Omega') T^{n-m}\left(m_{\cV^+_{Q_{n,L_0}}}\left(\left\{\bm{v}\in \cV^+_{Q_{n,L_0}}: \|\bm{v} g^{-1}\|\leq 1\right\}\right)+O\left( \|g\|_{\rm op}^{n-m+1} T^{-1/2}\right)\right).
	\end{align*}
	Finally, \comm{recall} that $\cV^+_{Q_{n,L_0}}g^{-1}=\cV^+_{Q_{n,L}}$ and $m_{\cV^+_{Q_{n,L}}}$ is defined as the pushforward of $m_{\cV^+_{Q_{n,L_0}}}$ by $g^{-1}$. Using also the relation
	$\vol(\Omega')=\frac{\vol(\Omega)}{|\det(h)|}$ we get our result.
\end{proof}

\subsection{Volume calculation for homogenous functions}
We now state the volume estimate needed for \thmref{thm:homogen}. Before doing so, we first introduce some more notation.

Fix $1<d<m$ and given $F\in \mathscr{F}_{d,Q_n}$ with $F(\bm{v})=F_{p,q}^{(d)}(\bm{v}gL_0 h)$ for some $p\geq q\geq 1$ with $p+q=m$, $g\in \SO_{Q_n}^+(\R)$ and $h\in \GL_m(\R)$, we denote by 
$
\cV^+_{Q_n,F}=\{\bm{v}\in \cV^+_{Q_n}: F(\bm{v})=0\}.
$ 
Note that $\bm{v}\in \cV^+_{Q_n,F}$ is of the form 
$\bm{v}g=r(\bm{u}, \bm{\omega}\sqrt{1-\|\bm{u}\|^2},1)$ for some $r>0$, $\bm{\omega}\in S^{n-m}$ and $\bm{u}\in \R^m$ with $\|\bm{u}\|\leq 1$ and $F^{(d)}_{p,q}(\bm{u}h)=0$. Further decomposing $\bm{u}h=(\bm{u}_1,\bm{u}_2)$ with $\bm{u}_1\in \R^p,\bm{u}_2\in \R^q$ then the condition $F^{(d)}_{p,q}(\bm{u}h)=0$ implies that $\|\bm{u}_1\|_d=\|\bm{u}_2\|_d$. Letting $S_d^{p,q}=S_d^{p-1}\times S_d^{q-1}$ denote the product of two unit spheres with respect to the $L^d$-norm, we can write $\bm{u}h=\tfrac{\|\bm{u}h\|}{\|\tilde{\bm{\omega}}\|}\tilde{\bm{\omega}}$ with $\tilde{\bm{\omega}}\in S_d^{p,q}$. Let $t=\|\bm{u}\|\in [0,1]$ so that $\bm{u}=t\tfrac{\tilde{\bm{\omega}} h^{-1}}{\|\tilde{\bm{\omega}} h^{-1}\|}$. 
Thus we can write $\bm{v}\in \cV^+_{Q_n,F}$ of the form
$$
\bm{v}=r\left(t\tfrac{\tilde{\bm{\omega}} h^{-1}}{\|\tilde{\bm{\omega}} h^{-1}\|}, \bm{\omega}\sqrt{1-t^2},1\right)g^{-1}
$$
for some $r>0$, $\tilde{\bm{\omega}}\in S_d^{p,q}$, $\bm{w}\in S^{n-m}$ and $0\leq t\leq 1$. Under these coordinates we define the measure $m_{\cV^+_{Q_n,F}}$ on $\cV^+_{Q_n,F}$ by
\begin{align}\label{e:sigmaVQF}
	\d m_{\cV^+_{Q_n,F}}(\bm{v}):=\|\tilde{\bm{\omega}}h^{-1}\|^{d-m}(1-t^2)^{\frac{n-m-1}{2}} t^{m-d-1}r^{n-d-1}\, \d t\d r\d\sigma_{n-m}(\bm{\omega})\d\sigma_d^{p,q}(\tilde{\bm{\omega}}),
\end{align} 
where $\sigma_d^{p,q}:=\sigma_d^p\times \sigma_d^q$ is the surface area on $S_d^{p,q}$ such that $\d\bm{x}=t_1^{p-1}t_2^{q-1}\d\sigma_d^{p,q}(\tilde{\bm{\omega}})$ under the polar coordinates $\bm{x}=(t_1\bm{\omega}_1,t_2\bm{\omega}_2)\in \R^p\times \R^q$ with $t_1,t_2>0$ and $\tilde{\bm{\omega}}=(\bm{\omega}_1,\bm{\omega}_2)\in S_d^{p,q}$.

We now state our volume estimate for these homogeneous functions.
\begin{Prop}\label{p:volFinI}
Let $F\in \mathscr{F}_{d,Q_n}$ be such that $F(\bm{v})=F_{p,q}^{(d)}(\bm{v}gL_0 h)$ for some $p\geq  q\geq 1$ with $p+q=m$, $g\in G$ and $h\in \GL_m(\R)$ and let $N\geq 1$. For any measurable $I\subseteq [-N,N]$ and for any $T> 2N^{1/d}$,
	$$m_{\cV^+_{Q_n}}\left(\left\{\bm{v}\in \cV^+_{Q_n}: \|\bm{v}\|\leq T,\ F(\bm{v})\in I\right\}\right)= \frac{{c_{n,m}} |I| T^{n-d}}{|\det(h)|d}\left(V_F +O_{h,g}\left(T^{-\nu}\right) \right),$$
	where {$c_{n,m}$ is as in \propref{p:volLinOmega}}, $V_F:=m_{\cV^+_{Q_n,F}}\left(\{\bm{v}\in \cV^+_{Q_n,F}: \|\bm{v}\|\leq 1\}\right)$, $\nu=\min\{{\frac{d}{4}}, \frac{m-d}{2}\}$ and the implied constant can be bounded by some powers of $\max\{\|h\|_{\rm op}, \|h^{-1}\|_{\rm op}\}\|g\|_{\rm op}$.
\end{Prop}
As before, we first prove a smoothed version.
\begin{Lem}\label{l:VolFinI}
	Let $F\in \mathscr{F}_{d,Q_n}$ with $F(\bm{v})=F_{p,q}^{(d)}(\bm{v}L_0h)$ for some pair $(p, q)\in \N^2$ with $p+q=m$ and $h\in \GL_m(\R)$. Let $f\in C^1_c(\R^{n+2})$ be supported in a ball of radius $c_f\geq 1$. For any measurable $I\subseteq [-1,1]$ and for any $T\geq 2$,
	$$\int_{\cV^+_{Q_n}} f\left(\tfrac{\bm{v}}{T}\right) \chi_I\left(F(\bm{v})\right)\, \d m_{\cV^+_{Q_n}}(\bm{v})=\frac{ c_{n,m}|I| T^{n-d}}{|\det(h)|d}\left(\int_{\cV^+_{Q_n,F}}f\, \d m_{\cV^+_{Q_n,F}}+O_h\left( \cS_{1}(f) c_f^{n+1}  \left(  T^{d-m}+T^{-\frac{d}{2}}\log(T)\right)\right)\right),$$
	where the implied constant can be bounded by some powers of $\max\{\|h\|_{\rm op}, \|h^{-1}\|_{\rm op}\}$.
	
\end{Lem}
\begin{proof}
	Let 
$
		\cI:=\int_{\cV^+_{Q_n}} f\left(\tfrac{\bm{v}}{T}\right) \chi_I\left(F_{p,q}^{(d)}(\bm{v}L_0h)\right)\, \d m_{\cV^+_{Q_n}}(\bm{v}).
	$
	Using the same change of variables as in the proof of \lemref{l:LinOmega} we get that $\cI$ equals
	\begin{align*}
		{c_{n,m}}T^n\int_0^\infty\int_{\|\bm{u}\|<1}\int_{S^{n-m}}  f(r(\bm{u},\bm{\omega}\sqrt{1-\|\bm{u}\|^2},1))\chi_I(r^dT^dF_{p,q}^{(d)}(\bm{u}h)) (1-\|\bm{u}\|^2)^{\frac{n-m-1}{2}} r^{n-1}\, \d\sigma_{n-m}(\bm{\omega})\d\bm{u}\d r.
	\end{align*}
	After further making a change of variables $\bm{u}\mapsto \bm{u}h^{-1}$ we write the new variable $\bm{u}=\bm{u}(t_1,t_2,\tilde{\bm{\omega}})=(t_1\bm{\omega}_1,t_2\bm{\omega}_2)$ with $t_1,t_2>0$ and 
	$\tilde{\bm{\omega}}=(\bm{\omega}_1,\bm{\omega}_2)\in S_d^{p,q}$ in polar coordinates with respect to the $L^d$-norm on $\R^p$ and $\R^q$ respectively. In these coordinates $\d\bm{u}=t_1^{p-1}t_2^{q-1}\, \d\sigma_{d}^{p,q}(\tilde{\bm{\omega}})$ and $F_{p,q}^{(d)}(\bm{u})=t_1^d-t_2^d$ so that
	\begin{align*}
		\cI&= \frac{c_{n,m}T^n}{|\det(h)|}\int_0^\infty \int_{S^{p,q}_d}\int_{\{(t_1,t_2)\col \|\bm{u}(t_1,t_2,\tilde{\bm{\omega}})h^{-1}\|\leq 1\}}\int_{S^{n-m}}f\left(r(\bm{u}h^{-1},\bm{\omega}\sqrt{1-\|\bm{u}h^{-1}\|^2},1)\right)\\
		&\chi_I\left(r^dT^d(t_1^d-t_2^d)\right) (1-\|\bm{u}h^{-1}\|^2)^{\frac{n-m-1}{2}}  r^{n-1}t_1^{p-1}t_2^{q-1}\, \d\sigma_{n-m}(\bm{\omega})\d t_1\d t_2\d\sigma_d^{p,q}(\tilde{\bm{\omega}})\d r.
	\end{align*}
	Now let $I_1=I\cap [0,1]$ and $I_2=I\cap [-1,0)$ so that $\chi_I=\chi_{I_1}+\chi_{I_2}$. For $i=1,2$, let $\cI_i$ be the above integral with $\chi_I$ replaced by $\chi_{I_i}$ so that $\cI=\cI_1+\cI_2$. To compute $\cI_1$ make another change of variables 
	$s=t_1^d-t_2^d,\; t=t_2$ so that $t_1=(s+t^d)^{1/d}$ and $\d t_1\d t_2=\frac{1}{d}(s+t^d)^{\frac{1-d}{d}}\, \d s\d t$. In these new variables we get
	\begin{align*}
		\cI_1&= \frac{c_{n,m}T^n}{|\det(h)|d}\int_0^\infty \int_{S^{p,q}_d}\int_{0}^{\infty}\int_{\{t\col \|\bm{u}(s,t,\tilde{\bm{\omega}})h^{-1}\|\leq 1\}}\int_{S^{n-m}} f\left(r(\bm{u}h^{-1},\bm{\omega}\sqrt{1-\|\bm{u}h^{-1}\|^2},1)\right)\\
		&\chi_{I_1}\left(r^dT^d s\right) (1-\|\bm{u}h^{-1}\|^2)^{\frac{n-m-1}{2}}  r^{n-1}(s+t^d)^{\frac{p-d}{d}} t^{q-1}\, \d\sigma_{n-m}(\bm{\omega})\d t\d s\d\sigma_d^{p,q}(\tilde{\bm{\omega}})\d r.
	\end{align*}
	Here with a slight abuse of notation $\bm{u}(s,t,\tilde{\bm{\omega}})$ is understood as
	$\bm{u}(s,t,\tilde{\bm{\omega}})=\left((s+t^d)^{1/d}\bm{\omega}_1, t\bm{\omega}_2\right)$. 
	Note the condition $\|\left((s+t^d)^{1/d}\bm{\omega}_1, t\bm{\omega}_2\right)h^{-1}\|\leq 1$ implies that $\max\{(s+t^d)^{1/d}, t\}\ll \|h\|_{\rm op}$.
 can thus use the bounds $(s+t^d)^{\frac{p-d}{d}}\ll \|h\|_{\rm op}^{p-d}$ if $p\geq d$ and $(s+t^d)^{\frac{p-d}{d}}\leq t^{p-d}$ if $p<d$ to bound the contribution of $r\leq T^{-1}$ trivially by $O\left( \frac{\|f\|_\infty \|h\|_{\rm op}^{m-d} |I_1|}{|\det(h)|d} \right)$.
	For the contribution of $r> T^{-1}$ make the change of variables $s\mapsto (rT)^{-d} s$ and change order of integration to get that $\cI_1=\cI_1^++O\left( \frac{\|f\|_\infty \|h\|_{\rm op}^{m-d} |I_1|}{|\det(h)|d} \right)$ with
	\begin{align*}
		\cI_1^+&:= \frac{c_{n,m}T^{n-d}}{|\det(h)|d}\int_{I_1}\int_{T^{-1}}^\infty \int_{S^{p,q}_d}\int_{\{t\col \|\bm{u}(s,t,\tilde{\bm{\omega}})h^{-1}\|\leq 1\}}\int_{S^{n-m}} f\left(r(\bm{u}h^{-1},\bm{\omega} \sqrt{1-\|\bm{u}h^{-1}\|^2},1)\right)\\
		& (1-\|\bm{u}h^{-1}\|^2)^{\frac{n-m-1}{2}} r^{n-d-1}(\tfrac{s}{(rT)^d}+t^d)^{\frac{p-d}{d}} t^{q-1}\, \d\sigma_{n-m}(\bm{\omega})\d t\d\sigma_d^{p,q}(\tilde{\bm{\omega}})\d r\d s,
	\end{align*}
	where now $\bm{u}(s,t,\tilde{\bm{\omega}})=\left((\tfrac{s}{(rT)^d}+t^d)^{1/d}\bm{\omega}_1, t\bm{\omega}_2\right)$.
	Next, we bound the contribution of $t\leq(rT)^{-1}$ (conditioned $r> T^{-1}$) trivially by
	$
	O\left( \frac{\|f\|_\infty |I_1|(c_fT)^{n-m}}{|\det(h)|d}\right)
	$ and assume from now on that the integral over $t$ is restricted to $t> (rT)^{-1}$.
	When $t> (rT)^{-1}$ we can estimate $(\tfrac{s}{(rT)^d}+t^d)^{1/d}=t+O(\frac{t}{(rtT)^{d}})$ uniformly for all $s\in I_1\subseteq [0,1]$ to get that 
	\begin{align}\label{equ:apparg1}
	\|\bm{u}h^{-1}-t\tilde{\bm{\omega}}h^{-1}\|\ll \tfrac{\|h^{-1}\|_{\rm op}t}{(rtT)^d},
	\end{align}
	and assuming that $\|t\tilde{\bm{\omega}}h^{-1}\|\leq 1$ as well we can also bound
	\begin{align}\label{equ:apparg2}
	\left|\sqrt{1-\|\bm{u}h^{-1}\|^2}-\sqrt{1-\|t\tilde{\bm{\omega}}h^{-1}\|^2}\right|&\leq \sqrt{\left|\|t\tilde{\bm{\omega}}h^{-1}\|^2-\|\bm{u}h^{-1}\|^2\right|}\ll \left(\tfrac{t\|h^{-1}\|_{\rm op}}{(rtT)^{d}}\right)^{1/2}.
	\end{align}

Now for the integral over $t$, we replace the range $J_{r,T, s,\tilde{\omega}}:=\left\{t\geq (rT)^{-1}: \|\bm{u}(s,t,\tilde{\bm{\omega}})h^{-1}\|\leq 1\right\}$ by the smaller range $J'_{r, T, s,\tilde{\bm{\omega}}}:=\left\{t\geq (rT)^{-1}: \max\{\|\bm{u}(s,t,\tilde{\bm{\omega}})h^{-1}\|, \|t\tilde{\bm{\omega}}h^{-1}\|\}\leq 1\right\}$. Note that if $t\in J_{r, T, s,\tilde{\bm{\omega}}}\setminus J'_{r, T, s,\tilde{\bm{\omega}}}$, then $\|h^{-1}\|_{\rm op}^{-1}\ll t\ll \|h\|_{\rm op}$, where the first estimate follows from $\|t\tilde{\bm{\omega}}h^{-1}\|>1$ and the second estimate follows from $\|\bm{u}h^{-1}\|\leq 1$ as before. 
	Moreover, we have by \eqref{equ:apparg1}
	$$
	1<\|t\tilde{\bm{\omega}}h^{-1}\|\leq \|\bm{u}h^{-1}\|+\|\bm{u}h^{-1}-t\tilde{\bm{\omega}}h^{-1}\|\leq 1+O\left(\tfrac{ \|h^{-1}\|_{\rm op}t^{1-d}}{(rT)^d}\right).
	$$
	Since $\|h^{-1}\|_{\rm op}^{-1}\ll t\ll \|h\|_{\rm op}$, we can bound $t^{1-d}\ll \max\{\|h\|_{\rm op}, \|h^{-1}\|_{\rm op}\}^{|1-d|}$ to get that
	 \begin{align}\label{equ:ragenfort}
	 t\in \left(\tfrac{1}{\|\tilde{\bm{\omega}}h^{-1}\|}, \tfrac{1}{\|\tilde{\bm{\omega}}h^{-1}\|}+\tfrac{C_h}{(rT)^d}\right)
	 \end{align} 
	 for some constant $C_h>0$ depending only on some powers of $\max\{\|h\|_{\rm op}, \|h^{-1}\|_{\rm op}\}$.
	Now using the estimate $t^d\leq \tfrac{s}{(rT)^d}+t^d\leq 2t^d$ for all $s\in I_1\subseteq [0,1]$ and again the estimates \eqref{equ:ragenfort} and $ t\ll \|h\|_{\rm op}$ for $t\in J_{r, T, s,\tilde{\bm{\omega}}}\setminus J'_{r, T, s,\tilde{\bm{\omega}}}$ we get
	\begin{align*}
		\cI_1^+&=\frac{c_{n,m}T^{n-d}}{|\det(h)|d}\int_{I_1}\int_{T^{-1}}^\infty \int_{S^{p,q}_d}\int_{J_{r, T, s,\tilde{\bm{\omega}}}'}\int_{S^{n-m}} f\left(r(\bm{u}h^{-1},\bm{\omega} \sqrt{1-\|\bm{u}h^{-1}\|^2},1)\right)\\
		& (1-\|\bm{u}h^{-1}\|^2)^{\frac{n-m-1}{2}} r^{n-d-1}(\tfrac{s}{(rT)^d}+t^d)^{\frac{p-d}{d}} t^{q-1}\, \d\sigma_{n-m}(\bm{\omega})\d t\d\sigma_d^{p,q}(\tilde{\bm{\omega}})\d r\d s\\
		&
		+ O_h\left(d^{-1}|\det(h)|^{-1}\|f\|_{\infty}|I_1|\left( c_f^{n-m}T^{n-m}+T^{n-2d}\int_{T^{-1}}^{c_f}r^{n-2d-1}\,\d r\right)\right).
	\end{align*}
Here the implied constant is bounded by some powers of $\max\left\{\|h\|_{\rm op}, \|h^{-1}\|_{\rm op}\right\}$.
	Moreover, note that the term 
	$T^{n-2d}\int_{T^{-1}}^{c_f}r^{n-2d-1}\d r=\int_1^{c_fT}r^{n-2d-1}\d r$
	 is bounded respectively by $(c_fT)^{n-2d}$ if $n-2d>0$, by $\log(c_fT)$ if $n-2d=0$ and by $1$ if $n-2d<0$. In all cases we have 
	\begin{align*}
		\cI_1^+&=\frac{c_{n,m}T^{n-d}}{|\det(h)|d}\int_{I_1}\int_{T^{-1}}^\infty \int_{S^{p,q}_d}\int_{J_{r, T, s,\tilde{\bm{\omega}}}'}\int_{S^{n-m}} f\left(r(\bm{u}h^{-1},\bm{\omega} \sqrt{1-\|\bm{u}h^{-1}\|^2},1)\right)\\
		& (1-\|\bm{u}h^{-1}\|^2)^{\frac{n-m-1}{2}} r^{n-d-1}(\tfrac{s}{(rT)^d}+t^d)^{\frac{p-d}{d}} t^{q-1}\, \d\sigma_{n-m}(\bm{\omega})\d t\d\sigma_d^{p,q}(\tilde{\bm{\omega}})\d r\d s\\
		&
		+ O_h\left( d^{-1}|\det(h)|^{-1}\|f\|_{\infty}|I_1| \left((c_fT)^{n-m}+(c_fT)^{n-2d}\right)\right).
	\end{align*}

	Next, for any $t\in J_{r, T, s,\tilde{\bm{\omega}}}'$ we can use \eqref{equ:apparg1} and \eqref{equ:apparg2} to estimate
	$$f\left(r(\bm{u}h^{-1},\bm{\omega}\sqrt{1-\|\bm{u}h^{-1}\|^2},1)\right)=f\left(r(t\tilde{\bm{\omega}}h^{-1},\bm{\omega} \sqrt{1-\|t\tilde{\bm{\omega}}h^{-1}\|^2},1)\right)+{O_h\left(\cS_{1}(f) r\left(\tfrac{t}{(rtT)^d}+\tfrac{t^{1/2}}{(rtT)^{d/2}}\right)\right)}.$$
	Similarly, since $t\geq (rT)^{-1}$ and $s\in I\subseteq [0,1]$ we can approximate $(\tfrac{s}{(rT)^d}+t^d)^{\frac{p-d}{d}}=t^{p-d}+O(\frac{ t^{p-2d}}{(rT)^d})$ and 
	$(1-\|\bm{u}h^{-1}\|^2)^{\frac{n-m-1}{2}}=(1-\|t\tilde{\bm{\omega}} h^{-1}\|^2)^{\frac{n-m-1}{2}}+{O_h(\frac{t^{1/2}}{(rtT)^{d/2}})}$. Thus using triangle inequality and the fact that $t\leq C\|h\|_{\rm op}$ for some absolute constant $C>0$ for $t\in J_{r, T, s,\tilde{\bm{\omega}}}'$ we can estimate 
	\begin{align*}
		\cI_1^+&=\frac{c_{n,m}T^{n-d}}{|\det(h)|d}\int_{I_1}\int_{T^{-1}}^\infty \int_{S^{p,q}_d}\int_{J_{r, T, s,\tilde{\bm{\omega}}}'}\int_{S^{n-m}} f\left(r(t\tilde{\bm{\omega}}h^{-1},\bm{\omega} \sqrt{1-\|t\tilde{\bm{\omega}}h^{-1}\|^2},1)\right)\\
		&(1-\|t\tilde{\bm\omega} h^{-1}\|^2)^{\frac{n-m-1}{2}} r^{n-d-1}t^{m-d-1}\, \d\sigma_{n-m}(\bm\omega)\d t\d\sigma_d^{p,q}(\tilde{\bm{\omega}})\d r\d s\\
		&+O_h\left(d^{-1}|\det(h)|^{-1}\cS_{1}(f)|I_1|T^{n-d}\int_{T^{-1}}^{c_f}\int_{(rT)^{-1}}^{C\|h\|_{\rm op}}t^{m-d-1}r^{n-d-1}\left(\tfrac{1}{(rtT)^d}+\tfrac{rt}{(rtT)^d}+\tfrac{(r+1)t^{1/2}}{(rtT)^{d/2}}\right)\,\d t\d r\right)\\
		&+O_h\left(d^{-1} |\det(h)|^{-1}\|f\|_{\infty}|I_1| \left((c_fT)^{n-m}+(c_fT)^{n-2d}\right)\right).
	\end{align*}
	The above error terms can be bounded by 
	$$
	O_h\left( d^{-1}|\det(h)|^{-1}\cS_{1}(f)|I_1|{\left(c_f(c_fT)^{n-m}+c_f(c_fT)^{n-\frac{3d}{2}}\log(c_fT)\right)}\right),$$
	where again the implied constant is bounded by some power of $\max\{\|h\|_{\rm op},\|h^{-1}\|_{\rm op}\}$ and the term $\log(c_fT)$ is needed only when $m=\frac{3d-1}{2}$.

	Now using similar estimates as above, after adding an error term which can be absorbed by the above error terms we can replace $J_{r, T, s,\tilde{\bm{\omega}}}'$ by the larger range $J_{r, T, s,\tilde{\bm{\omega}}}'':=\left\{t\geq (rT)^{-1}: \|t\tilde{\bm{\omega}}h^{-1}\| \leq 1\right\}$ and add back the contributions of $t\leq (rT)^{-1}$ (conditioned $r> T^{-1}$) and $r\leq T^{-1}$ to get
	\begin{align*}
		\cI_1&=\frac{c_{n,m}T^{n-d}|I_1|}{|\det(h)|d}\int_{0}^\infty \int_{S^{p,q}_d}\int_{\{t\col \|t\tilde{\bm{\omega}}h^{-1}\| \leq 1\}}\int_{S^{n-m}} f\left(r(t\tilde{\bm{\omega}}h^{-1},\bm{\omega} \sqrt{1-\|t\tilde{\bm{\omega}}h^{-1}\|^2},1)\right)\\
		&(1-\|t\tilde{\bm{\omega}} h^{-1}\|^2)^{\frac{n-m-1}{2}} r^{n-d-1}t^{m-d-1}\, \d\sigma_{n-m}(\bm{\omega})\d t\d\sigma_d^{p,q}(\tilde{\bm{\omega}})\d r\\
		&+O_h\left( d^{-1}|\det(h)|^{-1}\cS_{1}(f)\left(c_f(c_fT)^{n-m}+c_f(c_fT)^{n-\frac{3d}{2}}\log(c_fT)\right)\right).
	\end{align*}
	Finally, for the main term making a change of variables $t\mapsto \frac{t}{\|\tilde{\bm{\omega}} h^{-1}\|}$ we get that 
	\begin{align*}
		\cI_1 &= \frac{T^{n-d}|I_1|}{|\det(h)|d}\left(\int_{\cV^+_{Q_n,F}}f\, \d m_{\cV^+_{Q_n,F}}+O_h\left( c_f^{n+1}\cS_{1}(f) T^{d-m}+T^{-\frac{d}{2}}\log(T)\right)\right).
	\end{align*}
	To compute $\cI_2$ we can similarly make the change of variables $s=t_2^d-t_1^d,\; t=t_1$ to get
	\begin{align*}
		\cI_2&= \frac{c_{n,m}T^n}{|\det(h)|d}\int_0^\infty \int_{S^{p,q}_d}\int_{0}^{\infty}\int_{\{t\col \|\bm{u}(s,t,\tilde{\bm{\omega}})h^{-1}\|\leq 1\}}\int_{S^{n-m}} f(r(\bm{u}h^{-1},\bm{\omega}\sqrt{1-\|\bm{u}h^{-1}\|^2},1))\\
		&\chi_{-I_2}\left(r^dT^d s\right) (1-\|\bm{u}h^{-1}\|^2)^{\frac{n-m-1}{2}}  r^{n-1}(s+t^d)^{\frac{q-d}{d}} t^{p-1}\, \d\sigma_{n-m}(\bm{\omega})\d t\d s\d\sigma_d^{p,q}(\tilde{\bm{\omega}})\d r.
	\end{align*} 
Here now $\bm{u}(s,t,\tilde{\bm{\omega}})=\left(t\bm{\omega}_1, (s+t^d)^{1/d}\bm{\omega}_2\right)$. By similar computation as above we get
\begin{align*}
\cI_2 &= \frac{c_{n,m}T^{n-d}|I_2|}{|\det(h)|d}\left(\int_{\cV^+_{Q_n,F}}f\, \d m_{\cV^+_{Q_n,F}} +O_h\left( c_f^{n+1}\cS_{1}(f) \left(T^{d-m}+T^{-\frac{d}{2}}\log(T)\right)\right)\right).
\end{align*}
Summing these two terms concludes the proof. 
\end{proof}

By a similar unsmoothing argument \comm{as in the proof of \propref{p:volLinOmega}} we \comm{now prove \propref{p:volFinI}}.

\begin{proof}
	Since $F$ is homogeneous of degree $d$, using a simple scaling trick, it suffices to prove this proposition for the case when $N=1$, i.e. $I\subseteq [-1,1]$. 
	For the remaining, let $\psi_0$ be the indicator function of the unit ball in $\R^{n+2}$ so that 
	$$
	R:=m_{\cV_{Q_n}^+}\left(\left\{\bm{v}\in \cV^+_{Q_n}: \|\bm{v}\|\leq T,\ F(\bm{v})\in I\right\}\right)= \int_{\cV^+_{Q_n}} \psi_{0}\left(\tfrac{\bm{v}}{T}\right)\chi_{I}\left(F(\bm{v})\right)\,\d m_{\cV^+_{Q_n}}(\bm{v}).
	$$
	For $0<\delta<1$ let $\psi^\pm_\delta$ approximate $\psi_0$ in the sense that 
	$$\psi_\delta^+(\bm{v})=\left\lbrace\begin{array}{ll} 1& \|\bm{v}\|\leq 1,\\ 
	0 & \|\bm{v}\|\geq 1+\delta,
	\end{array}\right.$$ 
	and similarly 
	$$
	\psi_\delta^-(\bm{v})=\left\lbrace\begin{array}{ll} 1& \|\bm{v}\|\leq 1-\delta, \\
	 0 & \|\bm{v}\|\geq 1,\end{array}\right.
	 $$ 
	with
	$\cS_{1}(\psi_\delta^\pm)\asymp \delta^{-1}$.  We define 
	$\psi_{\delta,g}^\pm(\bm{v}):=\psi_{\delta}^\pm (\bm{v} g^{-1})$ and  $F_0(\bm{v}):=F_{p,q}^{(d)}(\bm{v}L_0h)$.  Note that after making a change of variable $\bm{v}\mapsto \bm{v}g^{-1}$ we get that
	$$ \int_{\cV^+_{Q_n}} \psi^-_{\delta,g}\left(\tfrac{\bm{v}}{T}\right)\chi_{I}\left(F_0(\bm{v})\right)\, \d m_{\cV^+_{Q_n}}(\bm{v})\leq R\leq \int_{\cV^+_{Q_n}} \psi^+_{\delta,g}\left(\tfrac{\bm{v}}{T}\right)\chi_{I}\left(F_0(\bm{v})\right)\, \d m_{\cV^+_{Q_n}}(\bm{v}).$$
	Applying \lemref{l:VolFinI} with $f=\psi^\pm_{\delta,g}$ and noting that $\cS_{1}(f)\ll \|g\|_{\rm op}\delta^{-1}$ and $c_f\ll \|g\|_{\rm op}$ we get 
	\begin{align*}
		\int_{\cV^+_{Q_n}} \psi^\pm_{\delta,g}\left(\tfrac{\bm{v}}{T}\right)\chi_{I}\left(F_0(\bm{v})\right)\, \d m_{\cV^+_{Q_n}}(\bm{v})&=\frac{c_{n,m}|I| T^{n-d}}{|\det (h)|d}\left(\int_{\cV^+_{Q_n,F_0}}\psi^\pm_{\delta,g}\, \d m_{\cV^+_{Q_n,F_0}} 
		+O_{h,g}\left( \delta^{-1} \left( T^{d-m}+T^{-\frac{d}{2}}\log(T)\right)\right)\right),
	\end{align*}
where the implied constant is bounded by some powers of $\max\{\|h\|_{\rm op}, \|h^{-1}\|_{\rm op}\}\|g\|_{\rm op}$ {and the term $\log(T)$ is needed only when $m=\frac{3d-1}{2}$}.	Next we can bound 
	$$\left|\int_{\cV^+_{Q_n,F_0}}\psi^\pm_{\delta,g}\, \d m_{\cV^+_{Q_n,F_0}}-\int_{\cV^+_{Q_n,F_0}}\psi_{0,g}\, \d m_{\cV^+_{Q_n,F_0}}\right|\ll \|h\|_{\rm op}^{m-d} \|g\|_{\rm op}^{n-d}\delta $$
	to get that 
	$$\comm{R}=\frac{ c_{n,m}|I| T^{n-d}}{|\det(h)|d}\left( \int_{\cV^+_{Q_n,F_0}}\psi_{0,g}\, \d m_{\cV^+_{Q_n,F_0}}+O_{h,g}\left( \delta^{-1} \left( T^{d-m}+T^{-\frac{d}{2}}\log(T)\right) +\delta\right)\right).$$
	Let $\nu=\min\{\frac{d}{4},\frac{m-d}{2}\}$, then taking  $\delta=T^{-\nu}$ (and noting that when the term $\log(T)$ occurs, i.e. when $m=\frac{3d-1}{2}$, $\frac{m-d}{2}<\frac{d}{4}$) gives  
	\begin{align*}
		\comm{R}&= \frac{c_{n,m} |I| T^{n-d}}{|\det(h)|d}\left(\int_{\cV^+_{Q_n,F_0}}\psi_{0,g}\, \d m_{\cV^+_{Q_n,F_0}}+O_{h,g}\left(T^{-\nu}\right) \right).
	\end{align*}
	Finally, we conclude the proof by noting that $\int_{\cV^+_{Q_n,F_0}}\psi_{0,g}\, \d m_{\cV^+_{Q_n,F_0}}=\int_{\cV^+_{Q_n,F}}\psi_{0}\, \d m_{\cV^+_{Q_n,F}}$, which follows from the change of variable $\bm{v}\mapsto \bm{v}g$.
\end{proof}

\subsection{Proofs of Theorems \ref*{thm:linearform} and \ref*{thm:homogen}}
Finally we give the proofs of \thmref{thm:linearform} and \thmref{thm:homogen} which is a direct application of \thmref{t:GenericFunctions} in view of the previous volume computations. We only give the proof of \thmref{thm:linearform} and note that \thmref{thm:homogen} can be proved similarly \comm{by replacing \propref{p:volLinOmega} in the proof outline below with \propref{p:volFinI}}.

\begin{proof}[Proof of \thmref{thm:linearform}]
We first fix some notation to fit the general setting of the previous section in order to apply \thmref{t:GenericFunctions}. First let $\G\bk G$ and $\mu_{\G}$ be as fixed in the beginning of this section. We also set $\cV=\cV_{Q_n}^+$ and $m_{\cV}=\omega_Qm_{\cV^+_{Q_n}}$.
Now take $\cB=\{B_t\}_{t>1}$ with $B_t:=\left\{\bm{v}\in \R^{n+2}: \|\bm{v}\|\leq t\right\}$ and $\cO=\{\cO_{\e}\}_{0<\e<1}$ with 
$$
\cO_{\e}=\left\{g\in G: \max\left\{\|g\|_{\rm op}, \|g^{-1}\|_{\rm op}\right\}<1+\e\right\}\subseteq G.
$$
It is then not difficult to check, using the multiplicativity of operator norm, that $\cO$ satisfies condition \eqref{e:OeReg} with $d=\dim_{\R}(G/K)=n+1$, $\eta=1$ and $\delta>0$ some constant depending only on $n$. The condition \eqref{equ:stabsmper} is also immediate again in view of the multiplicativity of operator norm. Moreover, \propref{p:volLinOmega} verifies the volume condition \eqref{e:VolEst} for the family of bounded measurable sets $\{\Omega_T\}_{T>1}$ in this theorem and the family of linear forms $\{gL_0h\}_{g\in G}$ (with $h\in \GL_m(\R)$ fixed).  This theorem then follows directly by applying \thmref{t:GenericFunctions}. Finally, we also note that the normalizing factor $\omega_Q$ in the definition of $m_{\cV}$ is responsible for the dependence on $Q$ of the leading coefficient of our counting formula.
\end{proof}

\section{Quantitative intrinsic Diophantine approximation}\label{sec:intrdioapp}
In this section, we apply the general counting results obtained in section \ref{sec:meantocount} to prove Theorems \ref{thm:c1digen}--\ref{thm:khinttypegen}. We follow the notation introduced in sections \ref{sec:indioapp} and \ref{sec:conj}. \comm{More explicitly, throughout this section we consider quadratic forms $Q(\bm{x},y)=\cQ(\bm{x})-y^2$ as in \eqref{equ:qdfspfo}. For such that $Q$ we can take $\tau\in \GL_{n+2}(\R)$ of the form $\tau=\left(\begin{smallmatrix}
 \tilde{\tau} & \\
  & 1\end{smallmatrix}\right)$ with $\tilde{\tau}\in \GL_{n+1}(\R)$ such that $\cQ(\bm{x})=\|\bm{x}\tilde{\tau}\|^2$. We fix $G=\SO_{Q_n}^+(\R)$, $\G=\tau^{-1}\SO_Q^+(\Z)\tau< G$ and $\tilde\cL_Q=\cL_Q\tau\subseteq \cV_{Q_n}^+$ as in section \ref{sec:ranform}. For any finite-measure $B\subset \cV_{Q_n}^+$ we abbreviate $\omega_Qm_{\cV_{Q_n}^+}(B)$ by $|B|$. Let $S_{\cQ}$  be the ellipsoid defined in \eqref{equ:ellipsoid} and recall the norm $\|\bm{x}\|_{\cQ}:=\sqrt{\cQ(\bm{x})}$ on $\R^{n+1}$. For any $\bm{x}\in S_{\cQ}$ and $r>0$ let $\mathfrak{D}^{\cQ}_r(\bm{x})$ be the open $r$-ball on $S_{\cQ}$ be defined in \eqref{equ:opendiscell}. When $\cQ(\bm{x})=\|\bm{x}\|^2$ (so that $S_{\cQ}=S^n$) we abbreviate $\mathfrak{D}^{\cQ}_r(\bm{x})$ by $\mathfrak{D}_r(\bm{x})$.}

\subsection{Reduction to counting points on the light cone}
In this section we \comm{reduce counting rational solutions to Diophantine inequalities to counting $\tilde{\cL}_Q$-points on the light cone
$\cV_{Q_n}^+$. 
We first introduce the sets that are relevant for our reduction. Recall from \rmkref{rmk:spherehomeo} that for any $\bm{\alpha}\in S^n$ we have fixed a rotation $k_{\bm{\alpha}}\in K$ with $K$ as in \eqref{equ:macomptgpk} 
satisfying $\bm{e}_0k_{\bm{\alpha}}=\tilde{\bm{\alpha}}$ with $\bm{e}_0=(-1,\bm{0},1)\in \cV_{Q_n}^+$ and $\tilde{\bm{\alpha}}:=(\bm{\alpha},1)\in \cV_{Q_n}^+$.} For any $T>0$, $r\in (0, 1)$ and $\bm{\alpha}\in S^n$, let 
\begin{align*}
	S_{T, r,\bm{\alpha}}:=\left\{\bm{e}_0a_yk_{\bm{\alpha}'}\in\cV_{Q_n}^+: y>T^{-1},\ \bm{\alpha}'\in \mathfrak{D}_r(\bm{\alpha})\right\}
\end{align*}
be the sector with height $T$ and spherical part given by $\mathfrak{D}_r(\bm{\alpha})$. 
We also consider another family of sets: For any continuous, decreasing function $\psi: [0,\infty)\to \R_{>0}$ and for any $T>0$ define
\begin{align}\label{equ:epsitset}
E_{\psi,T}:=\left\{\bm{v}\in \cV_{Q_n}^+: \left\|\bm{e}_0-\tfrac{\bm{v}}{v_{n+2}}\right\|<\psi(v_{n+2}),\ 0< v_{n+2}<T\right\}.
\end{align}
 The starting point of this reduction is the simple observation that, for $Q(\bm{x},y)=\cQ(\bm{x})-y^2$, $(\bm{p},q)\in \cL_Q$ if and only if $\frac{\bm{p}}{q}$ is a rational point on $S_{\cQ}$ written in lowest terms. 
 We have the following simple lemma relating the counting of rational points on $S_{\cQ}$ to counting points in corresponding $\tilde\cL_Q$-translates. 

\begin{Lem}\label{lem:relatingtoell}
Fix $\bm{\alpha}\in S^n$ and let $\bm{x}=\bm{\alpha}\tilde{\tau}^{-1}\in S_{\cQ}$. Then for any $0<r<1$ and $T>1$,
\begin{align}\label{equ:relaticount}
\cN_{\cQ}(\bm{x}, r,T)=\#\left(\tilde\cL_Q\cap S_{T,r,\bm{\alpha}}\right),
\end{align}
and for any continuous decreasing function $\psi: [0,\infty)\to \R_{>0}$ and any $T>1$,
\begin{align}\label{equ:relaticount2}
\cN_{\cQ,\psi}(\bm{x}, T)=\#\left(\tilde\cL_Q k_{\bm{\alpha}}^{-1}\cap E_{\psi,T}\right).
\end{align}
\comm{Here $\cN_{\cQ}(\bm{x}, r,T)$ and $\cN_{\cQ,\psi}(\bm{x}, T)$ are the two counting functions defined in \eqref{equ:countingfuncell} and \eqref{equ:khincountfunell} respectively.}
\end{Lem}
\begin{proof}
Let $(\bm{p},q)\in \cL_Q$. For \eqref{equ:relaticount}, note that $(\bm{p},q)\tau=(\bm{p}\tilde{\tau},q)\in \tilde\cL_Q\cap S_{T,r,\bm{\alpha}}$ if and only if $\|\tfrac{\bm{p}\tilde{\tau}}{q}-\bm{\alpha}\|<r$ and $1\leq q<T$. The condition $\|\tfrac{\bm{p}\tilde{\tau}}{q}-\bm{\alpha}\|<r$ \comm{is} equivalent to $\|\tfrac{\bm{p}}{q}-\bm{x}\|_{\cQ}<r$. This finishes the proof of \eqref{equ:relaticount}. 

Next, for \eqref{equ:relaticount2}, it suffices to show $\cN_{\cQ,\psi}(\bm{x}, T)=\#(\tilde\cL_Q\cap E_{\psi,T}k_{\bm{\alpha}})$. First note that 
\begin{align*}
E_{\psi,T}k_{\bm{\alpha}}=\left\{\bm{v}\in \cV_{Q_n}^+: \left\|\tilde{\bm{\alpha}}-\tfrac{\bm{v}}{v_{n+2}}\right\|<\psi(v_{n+2}),\ 0< v_{n+2}<T\right\}.
\end{align*}
This follows from the relation $\bm{e}_0k_{\bm{\alpha}}=\tilde{\bm{\alpha}}$ and by noting that the $(n+2)$-th coordinate of $\bm{v}\in \cV_{Q_n}^+$ remains unchanged under the right multiplication of $k_{\bm{\alpha}}$. Thus $(\bm{p},q)\tau=(\bm{p}\tilde{\tau},q)\in \tilde\cL_Q\cap E_{\psi,T}k_{\bm{\alpha}}$ if and only if $\|\tilde{\bm{\alpha}}-(\tfrac{\bm{p}\tilde{\tau}}{q},1)\|<\psi(q)$ and $1\leq q<T$. The first condition is further equivalent to $\|\bm{\alpha}-\tfrac{\bm{p}\tilde{\tau}}{q}\|<\psi(q)\ \Leftrightarrow\ \|\bm{x}-\tfrac{\bm{p}}{q}\|_{\cQ}<\psi(q)$.
This finishes the proof of \eqref{equ:relaticount2}.
\end{proof}

 \subsection{Measure estimates}\label{sec:volest}
\comm{In this section we prove various measure estimates regrading the sets $S_{T,r,\bm{\alpha}}$ and $E_{\psi, T}$, which in view of \lemref{lem:relatingtoell} and the counting results in section \ref{sec:meantocount}, will contribute the main terms in our counting formulas for $\cN_{\cQ}(\bm{x}, r,T)$ and $\cN_{\cQ,\psi}(\bm{x}, T)$.}

 We first prove a measure estimate for balls and differences of balls on spheres.
Recall that the spherical measure $\sigma_n$ has the following characterization \comm{that for} any Borel set $B\subseteq S^n$,
\begin{align}\label{equ:sphermeas}
	\sigma_n(B)=\frac{\vol(\widetilde{B})}{\vol(\widetilde{S}^n)},
\end{align}
where $\widetilde{B}:=\left\{t\bm{x}: 0\leq t\leq 1, \bm{x}\in B\right\}$. In particular, $\widetilde{S}^n$ is the closed unit ball in $\R^{n+1}$ and recall that $\vol(\widetilde{S}^n)=\frac{\pi^{\frac{n+1}{2}}}{\G(\frac{n+3}{2})}$.

\begin{Lem}\label{lem:volumecyl}
For any $\bm{\alpha}\in S^n$ and for any $r>0$,  
\begin{align}\label{equ:opendisfor}
		\sigma_n(\mathfrak{D}_r(\bm{\alpha}))=\tfrac{\G\left(\frac{n+3}{2}\right)}{\pi^{\frac12}(n+1)\G\left(\frac{n+2}{2}\right)}r^n+O_n(r^{n+2}).
	\end{align}
Moreover,	for any $0<r_2<r_1<1.1$, 
	\begin{align}\label{equ:differdiscs}
		\sigma_n\left(\mathfrak{D}_{r_1}\left(\bm{\alpha}\right)\setminus \mathfrak{D}_{r_2}\left(\bm{\alpha}\right)\right)\ll_{n}r_1^n-r_2^n.
	\end{align}                                                                                                                                       
\end{Lem}
\begin{proof}
	Up to enlarging the bounding constant, \eqref{equ:opendisfor} holds trivially for $r\geq 1$; thus we may assume $0<r<1$. Since $\sigma_n$ and the Euclidean norm (on $\R^{n+1}$) are both rotation-invariant, without loss of generality, we may assume $\bm{\alpha}=(1,\bm{0})$. Then by direct computation we have
	\begin{align*}
	\mathfrak{D}_r(\bm{\alpha})=\left\{\bm{v}\in S^n: v_1>1-\tfrac{r^2}{2}\right\}
	\end{align*}
	and 
	\begin{align*}
		\mathfrak{C}_{r_1,r_2}:=\mathfrak{D}_{r_1}(\bm{\alpha})\setminus \mathfrak{D}_{r_2}(\bm{\alpha})&=\left\{\bm{v}\in S^n: 1-\tfrac{r_1^2}{2}< v_1\leq 1-\tfrac{r_2^2}{2}\right\}.
	\end{align*}
	Writing $\bm{v}\in \R^{n+1}$ as $\bm{v}=(v_1,\bm{w})$ with $v_1\in \R$ and $\bm{w}\in \R^n$, one can verify that 
	\begin{align*}
		\left\{(v_1,\bm{w}): \|\bm{w}\|^2<a^2_rv_1^2,\ 0< v_1<1-r^2/2\right\}\subseteq		\widetilde{\mathfrak{D}}_r(\bm{\alpha})\subseteq \left\{(v_1,\bm{w}): \|\bm{w}\|^2\leq a_r^2v_1^2,\ 0\leq v_1\leq 1\right\},
	\end{align*}
and
	\begin{align*}
		\widetilde{\mathfrak{C}}_{r_1,r_2}\subseteq\left\{(v_1,\bm{w})\in \R^{n+1}: a_{r_2}^2v_1^2\leq \|\bm{w}\|^2\leq a_{r_1}^2v_1^2,\ 0\leq  v_1\leq 1-\tfrac{r_2^2}{2}
		\right\},
	\end{align*}
	where for any $t\in (0,1)$, $a_t:=\frac{t\sqrt{1-t^2/4}}{1-t^2/2}$.
	This, together with \eqref{equ:sphermeas} 
	implies that
	\begin{align*}
	\frac{\vol(\widetilde{S}^{n-1})}{\vol(\widetilde{S}^n)}\int_0^{1-\frac{r^2}{2}}a_{r}^nv_1^n \, \d v_1\leq \sigma_n(\mathfrak{D}_r(\bm{\alpha}))&\leq \frac{\vol(\widetilde{S}^{n-1})}{\vol(\widetilde{S}^n)}\int_0^{1}a_{r}^nv_1^n \, \d v_1,
	\end{align*}
	and
	\begin{align*}
		\sigma_n(\mathfrak{C}_{r_1,r_2})
		&\ll_n\int_0^{1-r_2^2/2}\left(a_{r_1}^n-a_{r_2}^n\right)v_1^n \, \d v_1
		\asymp_{n} a_{r_1}^n-a_{r_2}^n.
	\end{align*}
	Now \eqref{equ:opendisfor} follows from the estimate that $a_t=t+O(t^3)$ for all $t\in (0,1)$ together with the relation $\frac{\vol(\widetilde{S}^{n-1})}{(n+1)\vol(\widetilde{S}^n)}=\frac{\G\left(\frac{n+3}{2}\right)}{\pi^{\frac12}(n+1)\G\left(\frac{n+2}{2}\right)}$. 
	
	For \eqref{equ:differdiscs}, it suffices to show $a_{r_1}^n-a_{r_2}^n\ll_n r_1^2-r_2^n$. If $r_1>2r_2$ then
\begin{align*}
	a_{r_1}^n-a_{r_2}^n<a_{r_1}^n\asymp r_1^n\asymp r_1^n-r_2^n.
\end{align*}
	If $r_1\leq 2r_2$, let  $h(t):=a_t^n$ so that $h'(t)=\frac{8na_t^{n-1}}{\sqrt{4-t^2}(2-t^2)^2}$. Then by the mean value theorem we can estimate
	\begin{align*}
		a_{r_1}^n-a_{r_2}^n&= (r_1-r_2)h'(t) 
	\end{align*} 
for some $t\in [r_2,r_1]$. Now again using the estimate $a_{t}=t+O(t^3)$ for all $t\in (0,1)$, we see that $|h'(t)|\ll_n r_2^{n-1}$ for all $t\in [r_2, r_1]\subseteq [r_2, 2r_2]$, implying that in this case we also have
\begin{align*}
	a_{r_1}^n-a_{r_2}^n&\ll_n (r_1-r_2)r_2^{n-1}<r_1^n-r_2^n,
\end{align*}
as needed.
\end{proof}

This lemma  implies the following measure estimate for the sectors considered above: For any $T>0$, $r\in (0, 1)$ and $\bm{\alpha}\in S^n$, 
\begin{align}\label{equ:sectvolest}
m_{\cV_{Q_n}^+} (S_{T,r,\bm{\alpha}})=\tfrac{1}{n}T^n\sigma_n(\mathfrak{D}_r(\bm{\alpha}))=\tfrac{\G\left(\frac{n+3}{2}\right)}{\pi^{\frac12}n(n+1)\G\left(\frac{n+2}{2}\right)}T^n(r^n+O_n(r^{n+2})).
\end{align}
Next, fix $\psi: [0,\infty)\to \R_{>0}$ a continuous decreasing function, we compute the measure of the set $E_{\psi,T}$. 
\begin{Lem}\label{lem:volcompqdioapp} 
For any $T>0$ we have
\begin{align}\label{equ:volume1}
m_{\cV_{Q_n}^+} (E_{\psi,T})=\tfrac{\G\left(\frac{n+3}{2}\right)}{\pi^{\frac12}(n+1)\G\left(\frac{n+2}{2}\right)}\cJ_{\psi}(T)+O_{n}(\cI_{\psi}(T)),
\end{align}
where $\cJ_{\psi}(T):=\int_0^Tt^{n-1}\psi(t)^n\, \d t$ and $\cI_{\psi}(T):=\int_0^Tt^{n-1}\psi(t)^{n+2}\, \d t$.
\end{Lem}

\begin{proof}
Fix $T>0$ and let $\bm{\alpha}_0=(-1,\bm{0})\in S^n$ so that $(\bm{\alpha}_0,1)=\bm{e}_0$. It is then easy to see that under the polar coordinates $(y, \bm{\alpha})\in \R_{>0}\times S^n$ the set $E_{\psi, T}$ can be described as  following 
\begin{align*}
E_{\psi,T}&=\left\{\bm{e}_0a_yk_{\bm{\alpha}}\in \cV_{Q_n}^+ : \|\bm{\alpha}-\bm{\alpha}_0\|<\psi(y^{-1}),\ y> T^{-1}\right\},
\end{align*}
which then together with \eqref{equ:opendisfor} implies that
\begin{align*}
	m_{\cV_{Q_n}^+} (E_{\psi,T})&=\int_{T^{-1}}^{\infty}\sigma_n\left(\mathfrak{D}_{\psi(y^{-1})}(\bm{\alpha}_0) \right)y^{-(n+1)}\, \d y\\
	&=\tfrac{\G\left(\frac{n+3}{2}\right)}{\pi^{\frac12}(n+1)\G\left(\frac{n+2}{2}\right)}\int_{T^{-1}}^{\infty}y^{-(n+1)}\left(\psi(y^{-1})^n+O_n\left(\psi(y^{-1})^{n+2}\right)\right)\, \d y\\
	&=\tfrac{\G\left(\frac{n+3}{2}\right)}{\pi^{\frac12}(n+1)\G\left(\frac{n+2}{2}\right)}\int_{0}^Tt^{n-1}\psi(t)^n \, \d t+O_{n}\left(\int_{0}^Tt^{n-1}\psi(t)^{n+2}\, \d t\right).
\end{align*}
This finishes the proof.
\end{proof}

Recall that in \thmref{thm:khinttypegen} the approximating function $\psi$ is only defined on the set of positive integers, and the main term is given in term of the discrete sum  $J_{\psi}(T)=\sum_{1\leq q<T}q^{n-1}\psi(q)^n$. On the other hand here for our measure estimate $\psi$ is assumed to be continuous and the main term is given in terms of the integral $\cJ_{\psi}(T)=\int_0^Tt^{n-1}{\psi}(t)^n\, \d t$.
The following lemma shows that these two terms are the same up to a smaller error term.
\begin{Lem}\label{lem:esthold}
	\comm{Assume further} $\int_0^{\infty}t^{n-1}\psi(t)^n\, \d t=\infty$. Then for any $T>1$ we have
	\begin{align}\label{equ:estseint}
	\cJ_{\psi}(T)=J_{\psi}(T)+O_{\psi}\left(J_{\psi}(T)^{\frac{n+3}{n+4}}\right).
	\end{align}
\end{Lem}
\begin{proof}
	Fix $T>1$. For any integer $1\leq q< T$, since $\psi$ is decreasing, we can estimate
	\begin{align*}
	\int_q^{q+1}t^{n-1}\psi(t)^n\,dt\leq (q+1)^{n-1}\psi(q)^n=q^{n-1}\psi(q)^n+O_n(q^{n-2}\psi(q)^n).
	\end{align*}
	Similarly,
	\begin{align*}
	\int_q^{q+1}t^{n-1}\psi(t)^n\,dt\geq q^{n-1}\psi(q+1)^n=(q+1)^{n-1}\psi(q+1)^n+O_n\left((q+1)^{n-2}\psi(q+1)^n\right).
	\end{align*}
	From these two estimates one easily sees that
	\begin{align*}
	\cJ_{\psi}(T)=J_{\psi}(T)+O_{\psi}\left(R_{\psi}(T)+1\right),
	\end{align*}
	where $R_{\psi}(T):=\sum_{1\leq q\leq T}q^{n-2}\psi(q)^n$.
	 Since $\psi$ is continuous, decreasing and $\int_0^{\infty}t^{n-1}\psi(t)^n\, \d t=\infty$, we have $J_{\psi}(T)\to\infty$ as $T\to\infty$. Thus the error $O_{\psi}(1)$ can be absorbed into the desired error $O_{\psi}\left(J_{\psi}(T)^{\frac{n+3}{n+4}}\right)$, and it remains to show $R_{\psi}(T)\ll_{n,\psi} J_{\psi}(T)^{\frac{n+3}{n+4}}$ for all $T>1$. For this by taking 
	\begin{align*}
	(p,p')=(\tfrac{n+4}{n+3}, n+4)\quad \text{and}\quad (a_q, b_q)=(q^{\frac{(n-1)(n+3)}{n+4}}\psi(q)^{\frac{n(n+3)}{n+4}}, q^{-\frac{5}{n+4}}\psi(q)^{\frac{n}{n+4}})
\end{align*}
(so that $a_qb_q=q^{n-2}\psi(q)^n$ and $a_q^p=q^{n-1}\psi(q)^n$) and applying H\"{o}lder's inequality we get
	\begin{align*}
	R_{\psi}(T)=\sum_{1\leq q\leq T}a_qb_q\leq \left(\sum_{1\leq q\leq T}a_q^p\right)^{1/p}\left(\sum_{1\leq q\leq T}b_q^{p'}\right)^{1/p'}\ll_{n, \psi} J_{\psi}(T)^{\frac{n+3}{n+4}},
	\end{align*}
	where for the second estimate we used the estimate that 
	\begin{displaymath}
	\sum_{1\leq q\leq T}b_q^{p'}=\sum_{1\leq q\leq T}q^{-5}\psi(q)^n\leq \psi(1)^n\sum_{q\in\N} q^{-5}\ll_{\psi} 1.
	\end{displaymath}
	This finishes the proof.
	\end{proof}

In view of these measure estimates the proof of  Theorems \ref{thm:c1digen}-\ref{thm:khinttypegen} will follow from the following theorems respectively, applied to the form $Q(\bm{x},y)=\cQ(\bm{x})-y^2$. \comm{Recall the constant $\beta_Q$ given in \eqref{equ:defbq}.}
\begin{Thm}\label{thm:equidis}
	For any $r,T>0$ satisfying 
\comm{$T^{-\frac{2-\beta_Q}{3-\beta_Q}}\ll_Q r< 1$}, 
for any $\bm{\alpha}\in S^n$
	  and for any $g\in \comm{G}$,
	\begin{align}\label{equ:counsector}
	\#\left(\tilde\cL_Q g\cap S_{T,r,\bm{\alpha}}\right)=|S_{T,r,\bm{\alpha}}|+ O_{Q,g}\left( r^{-\frac{n}{2n+3}}|S_{T,r,\bm{\alpha}}|^{1-\frac{2-\beta_Q}{2n+3}}\right).
			\end{align}
	\end{Thm}

\begin{Thm}\label{thm:subtranspara}
Let $\bm{\alpha}_0=(-1,\bm0)\in S^n$. 
\begin{enumerate}
\item Let $r_T=T^{-\lambda}$ for some $0<\lambda<1$. Then for any $g\in \comm{G}$, for $\sigma_{n}$-a.e. $\bm{\alpha}\in S^n$ and for all $T$ sufficiently large
\begin{align*}
\#\left(\tilde\cL_Q k_{\bm{\alpha}}^{-1}g^{-1}\cap S_{T, r_T,\bm{\alpha}_0}\right)=\left|S_{T, r_T,\bm{\alpha}_0}\right|+O_{Q,g}\left(\left|S_{T, r_T,\bm{\alpha}_0}\right|^{1-\frac{2-\beta_Q}{n+4}}\log (\left|S_{T, r_T,\bm{\alpha}_0}\right|)\right).
\end{align*}
\item Let $r_T=\frac{(\log T)^{\lambda}}{T}$ for some $\lambda>\frac{1}{n(2-\beta_Q)}$. Then for any $g\in \comm{G}$, for $\sigma_n$-a.e. $\bm{\alpha}\in S^n$ and for all $T$ sufficiently large
\begin{align*}
\#\left(\tilde\cL_Q k_{\bm{\alpha}}^{-1}g^{-1}\cap S_{T, r_T,\bm{\alpha}_0}\right)=\left|S_{T, r_T,\bm{\alpha}_0}\right|+O_{Q,g}\left(\left|S_{T, r_T,\bm{\alpha}_0}\right|^{1-\frac{\lambda(2-\beta_Q)-1}{\lambda(d+3)}}\log \left(\left|S_{T, r_T,\bm{\alpha}_0}\right|\right)\right).
\end{align*}
\end{enumerate}
\end{Thm}

\begin{Thm}\label{thm:countmg}
Assume furhter $c_Q=0$.  Let $\psi: [0,\infty)\to (0,\frac12)$ be a continuous decreasing function satisfying $\lim\limits_{t\to\infty}\psi(t)=0$ and $\int_0^{\infty}t^{n-1}\psi(t)^n\, \d t=\infty$.
Then for any $g\in \comm{G}$, for $\sigma_n$-a.e. $\bm{\alpha}\in S^n$ and for all sufficiently large $T$,
\begin{align*}
\#\left(\tilde{\cL}_Q k_{\bm{\alpha}}^{-1}g^{-1}\cap E_{\psi, T}\right)=|E_{\psi, T}|+ O_{Q,g}\left(|E_{\psi,T}|^{\frac{n+3}{n+4}}\log(|E_{\psi, T}|)\right).
\end{align*}
\end{Thm}

The remaining of this section will be devoted to proving Theorems \ref{thm:equidis}-\ref{thm:countmg}.
We will prove these counting results by applying the general counting results developed in section \ref{sec:meantocount}. \comm{Indeed as mentioned in the end of section \ref{sec:meantocount}}, \thmref{thm:equidis}, \thmref{thm:subtranspara} and \thmref{thm:countmg} will follow from \thmref{thm:cfamd}, \thmref{thm:countnoninc} and \thmref{thm:countforall} respectively. In order to apply these results, one key step is to choose suitable identity neighborhoods and to verify well-roundedness of the relevant families of sets with respect to these identity neighborhoods.

\subsection{Choosing identity neighborhoods}
In this subsection we choose appropriate families of identity neighborhoods with respect to which we verify well-roundedness. 
First, to prove \thmref{thm:equidis} we choose the following family of identity neighborhoods in $G$. 
\begin{Lem}\label{lem:volest2}
	 For any $\e,r\in (0,1)$ and $\bm{\alpha}\in S^n$, define
\begin{align}\label{equ:neighbor2}
	G_{\e,r}(\bm{\alpha}):=\left\{g\in G: \|g\|_{\rm op}<1+\e,\
	 \max\left\{\left\| \tfrac{\tilde{\bm{\alpha}}g}{\| \tilde{\bm{\alpha}}g \|}-\tfrac{\tilde{\bm{\alpha}}}{\sqrt{2}}\right\|, \left\| \tfrac{\tilde{\bm{\alpha}}g^{-1}}{\| \tilde{\bm{\alpha}}g^{-1} \|}-\tfrac{\tilde{\bm{\alpha}}}{\sqrt{2}}\right\|\right\}<r\e\right\}, 
\end{align}
where $\tilde{\bm{\alpha}}:=(\bm{\alpha},1)\in \cV_{Q_n}^+$ is as before. 
Then $G_{\e, r}(\bm{\alpha})$ is an inversion invariant identity neighborhood with 
	\begin{align*}
	\mu_G(G_{\e,r}(\bm{\alpha}))\asymp_n r^n\e^{2n+1}\quad \text{and}\quad G_{\e_1,r}(\bm{\alpha})G_{\e_2,r}(\bm{\alpha})\subseteq G_{3(\e_1+\e_2),r}(\bm{\alpha}),\quad \forall\ \e,r\in (0,1), \e_1,\e_2\in (0,\tfrac13).
	\end{align*}
\end{Lem}

\begin{proof}
To show invariance under inversion it suffices to show $\|g\|_{\rm op}=\|g^{-1}\|_{\rm op}$ for all $g\in G$.
Indeed using the $KA^+K$-decomposition we can write $g=k_1a_yk_2$ with $k_1,k_2\in K$ and $a_y\in A$ with $y\geq 1$.  Since $\|\cdot\|_{\rm op}$ is bi-$K$-invariant, we get $\|g\|_{\rm op}=\|a_y\|_{\rm op}$ and similarly $\|g^{-1}\|_{\rm op}=\|a_y^{-1}\|_{\rm op}$. A simple calculation then shows that $\|a_y\|_{\rm op}=\|a_y^{-1}\|_{\rm op}=\max\{y, y^{-1}\}$, implying that $\|g\|_{\rm op}=\|g^{-1}\|_{\rm op}$. 

Next for the volume estimate, by the $K$-invariance of the Haar measure $\mu_G$ and of the Euclidean norm, we may assume $\bm{\alpha}=\bm{\alpha}_0\in S^n$.
We claim that 
\begin{equation}\label{e:GG'}
G_{\e/4,r}'\subseteq G_{\e,r}(\bm{\alpha}_0)\subseteq G'_{4\e,r},\quad \forall\ \e, r\in (0,1),
\end{equation}
where
\begin{align*}
G'_{\e,r}:=\left\{u_{\bm{x}}a_ymk_{\bm{\alpha}'}\in G : \|\bm{x}\|<\e,\ |y-1|<\e,\ m\in M,\ \|\bm{\alpha}'-\bm{\alpha}_0\|<r\e\right\}.
\end{align*}
To show the first inclusion, fix $r, \e\in(0,1)$ and take $g=u_{\bm{x}}a_ymk_{\bm{\alpha}'}\in G'_{\e/4,r}$. 
Then $$\|g\|_{\rm op}=\|u_{\bm{x}}a_y\|_{\rm op}\leq  \|u_{\bm{x}}\|_{\rm op}\|a_y\|_{\rm op}.$$
We have $\|a_y\|_{\rm op}=\max\{y,y^{-1}\}< (1+\frac{\e}{3})$ and a simple calculation shows that 
$
\|u_{\bm{x}}\|_{\rm op}=1+\frac{\|\bm{x}\|^2+\|\bm{x}\|\sqrt{\|\bm{x}\|^2+4}}{2}< 1+\frac{\e}{2},
$
so that $\|g\|_{\rm op}< 1+\e$. Next note that 
$\tilde{\bm{\alpha}}_0g=\bm{e}_0u_{\bm{x}}a_ymk_{\bm{\alpha}'}=y^{-1}\tilde{\bm{\alpha}}'$ and hence $ \tfrac{\tilde{\bm{\alpha}}_0g}{\| \tilde{\bm{\alpha}}_0g\|}=\tfrac{\tilde{\bm{\alpha}}'}{\sqrt{2}}$ and 
$$\left\| \tfrac{\tilde{\bm{\alpha}}_0g}{\| \tilde{\bm{\alpha}}_0g\|}-\tfrac{\tilde{\bm{\alpha}}_0}{\sqrt{2}}\right\|=\tfrac{1}{\sqrt{2}}\|\bm{\alpha}'-\bm{\alpha}_0\|\leq \tfrac{r\e}{4\sqrt{2}}<r\e.$$
Next, we can use again that $\|g\|_{\rm op} < 1+\e$ to get that 
$$
\left\| \tfrac{\tilde{\bm{\alpha}}_0}{\| \tilde{\bm{\alpha}}_0g\|}-\tfrac{\tilde{\bm{\alpha}}_0g^{-1}}{\sqrt{2}}\right\|=\left\| \left(\tfrac{\tilde{\bm{\alpha}}_0g}{\| \tilde{\bm{\alpha}}_0g\|}-\tfrac{\tilde{\bm{\alpha}}_0}{\sqrt{2}}\right)g^{-1}\right\|< (1+\e)\tfrac{r\e}{4\sqrt{2}}.$$
We can use the following simple geometric observation,
\begin{align}\label{equ:unitvecest}
	\left\|\bm{u}_1-\bm{u}_2\right\|^2\leq\frac{\|\bm{u}_1-\lambda\bm{u}_2\|^2}{\lambda},
	\quad \text{$\forall\ \bm{u}_1,\bm{u}_2\in S^n$ and $ \lambda>0$},
\end{align} 
with $\bm{u}_1=\frac{\tilde{\bm{\alpha}}_0}{\sqrt{2}}$,  $\bm{u}_2=\frac{\tilde{\bm{\alpha}}_0g^{-1}}{\|\tilde{\bm\alpha}_0g^{-1}\|}$ and $\lambda=\frac{\|\tilde{\bm{\alpha}}_0 g\| \tilde{\bm{\alpha}}_0 g^{-1}\|}{2}$ to get that 
\begin{align*}
\left\|\tfrac{\tilde{\bm{\alpha}}_0}{\sqrt{2}}-\tfrac{\tilde{\bm{\alpha}}_0g^{-1}}{\|\tilde{\bm\alpha}_0g^{-1}\|}\right\|
 &\leq
   \left(\tfrac{\|\tilde{\bm{\alpha}}_0 g\| }{\|\tilde{\bm{\alpha}}_0 g^{-1}\|}\right)^{1/2} \left\|\tfrac{\tilde{\bm{\alpha}}_0}{\|\tilde{\bm{\alpha}}_0 g\|}- \tfrac{\tilde{\bm{\alpha}}_0g^{-1}}{\sqrt{2}}\right\|< (1+\e)^2\tfrac{r\e}{4\sqrt{2}}< r\e.
   \end{align*}
   This shows that $G_{\e/4,r}'\subseteq G_{\e,r}(\bm{\alpha}_0)$. The arguments for the other inclusion are similar (and simpler) and we omit the details.

Now using  the Haar measure description \eqref{equ:Haar1} and the measure decomposition $\d \mu_K(mk_{\bm{\alpha}})=\d\mu_{M}(m)\d\sigma_n(\bm{\alpha})$ and the volume estimate \eqref{equ:opendisfor}, we can compute for any $\e, r\in (0,1)$, 
\begin{align*}
\mu_G\left(G'_{\e,r}\right)=\int_{\{\bm{x}\in \R^n\col \|\bm{x}\|<\e\}}\, \d\bm{x}\int_{1-\e}^{1+\e}\, y^n\, \d y \sigma_n\left(\mathfrak{D}_{r\e}(\bm{\alpha}_0)\right)\asymp_n r^n\e^{2n+1}.\end{align*}
This estimate together with \eqref{e:GG'} implies that 
\begin{align*}
r^n\e^{2n+1}\ll _n\mu_G\left(G'_{\e/4,r}\right)\leq \mu_G\left(G_{\e,r}(\bm{\alpha}_0)\right)\leq \mu_G\left(G'_{4\e,r}\right)\ll_n r^n\e^{2n+1},\end{align*}
thus proving this volume estimate.

Finally, we prove the translation relation. For any $\e_1,\e_2\in (0, \tfrac13)$, take $g_1\in G_{\e_1,r}(\bm{\alpha})$ and $g_2\in G_{\e_2,r}(\bm{\alpha})$ we need to show $g_1g_2\in G_{3(\e_1+\e_2),r}(\bm{\alpha})$. That is, $\|g_1g_2\|_{\rm op}\leq 1+3(\e_1+\e_2)$ and 
$$ \max\left\{\left\| \tfrac{\tilde{\bm{\alpha}}g_1g_2}{\| \tilde{\bm{\alpha}}g_1g_2 \|}-\tfrac{\tilde{\bm{\alpha}}}{\sqrt{2}}\right\|, \left\| \tfrac{\tilde{\bm{\alpha}}g_2^{-1}g_1^{-1}}{\| \tilde{\bm{\alpha}}g_2^{-1}g_1^{-1} \|}-\tfrac{\tilde{\bm{\alpha}}}{\sqrt{2}}\right\|\right\}<3r(\e_1+\e_2).
$$
The first relation is easy since 
$$
\|g_1g_2\|_{\rm op}\leq \|g_1\|_{\rm op}\|g_2\|_{\rm op}<(1+\e_1)(1+\e_2)<1+3(\e_1+\e_2).
$$
For the second relation, we only prove $\left\| \tfrac{\tilde{\bm{\alpha}}g_1g_2}{\| \tilde{\bm{\alpha}}g_1g_2 \|}-\tfrac{\tilde{\bm{\alpha}}}{\sqrt{2}}\right\|<3r(\e_1+\e_2)$ and the other inequality follows similarly. For this, first observe
\begin{align*}
\left\| \tfrac{\tilde{\bm{\alpha}}g_1g_2}{\| \tilde{\bm{\alpha}}g_1 \|}-\tfrac{\tilde{\bm{\alpha}}g_2}{\sqrt{2}}\right\|\leq \left\| \tfrac{\tilde{\bm{\alpha}}g_1}{\| \tilde{\bm{\alpha}}g_1 \|}-\tfrac{\tilde{\bm{\alpha}}}{\sqrt{2}}\right\|\|g_2\|_{\rm op}<r\e_1(1+\e_2).
\end{align*}
Now applying \eqref{equ:unitvecest} with $(\bm{u}_1=\tfrac{\tilde{\bm{\alpha}}g_1g_2}{\| \tilde{\bm{\alpha}}g_1g_2 \|},\bm{u}_2=\tfrac{\tilde{\bm{\alpha}}g_2}{\| \tilde{\bm{\alpha}}g_2 \|}$ and $\lambda=\tfrac{\|\tilde{\bm{\alpha}}g_1\|\|\tilde{\bm{\alpha}}g_2\|}{\sqrt{2}\|\tilde{\bm{\alpha}}g_1g_2\|}$ we get
\begin{align*}
\left\| \tfrac{\tilde{\bm{\alpha}}g_1g_2}{\| \tilde{\bm{\alpha}}g_1g_2 \|}-\tfrac{\tilde{\bm{\alpha}}g_2}{\| \tilde{\bm{\alpha}}g_2 \|}\right\|&
\leq\left(\tfrac{\sqrt{2}\| \tilde{\bm{\alpha}}g_1 \|}{\| \tilde{\bm{\alpha}}g_2 \|\| \tilde{\bm{\alpha}}g_1g_2 \|}\right)^{\frac12}\left\| \tfrac{\tilde{\bm{\alpha}}g_1g_2}{\| \tilde{\bm{\alpha}}g_1 \|}-\tfrac{\tilde{\bm{\alpha}}g_2}{\sqrt{2}}\right\|<{r\e_1(1+\e_2)^{2}}<3r\e_1.
\end{align*}
Hence by the triangle inequality we get
\begin{align*}
\left\| \tfrac{\tilde{\bm{\alpha}}g_1g_2}{\| \tilde{\bm{\alpha}}g_1g_2 \|}-\tfrac{\tilde{\bm{\alpha}}}{\sqrt{2}}\right\|\leq \left\| \tfrac{\tilde{\bm{\alpha}}g_1g_2}{\| \tilde{\bm{\alpha}}g_1g_2 \|}-\tfrac{\tilde{\bm{\alpha}}g_2}{\| \tilde{\bm{\alpha}}g_2 \|}\right\|+\left\| \tfrac{\tilde{\bm{\alpha}}g_2}{\| \tilde{\bm{\alpha}}g_2 \|}-\tfrac{\tilde{\bm{\alpha}}}{\sqrt{2}}\right\|<3r\e_1+r\e_2<3r(\e_1+\e_2).
\end{align*}
This finishes the proof.
\end{proof}

Next, to prove \thmref{thm:subtranspara} we choose the following family of identity neighborhoods of the parabolic subgroup $P$. 
\begin{Lem}\label{lem:idnbsubgp}
For any $\e\in (0,1)$ let $P_{\e}:=G_{\e,r}(\bm{\alpha}_0)\cap P$ for some $r\in (0,1)$. Then $P_{\e}$ is independent of the choice of $r$ and it satisfies
\begin{align*}
\mu_P(P_{\e})\asymp \e^{n+1},\ \forall\ \e\in (0,1)\qquad \text{and}\qquad P_{\e_1}P_{\e_2}\subseteq P_{3(\e_1+\e_2)},\ \forall\ \e_1,\e_2\in (0, \tfrac13),
\end{align*}
where for $p=u_{\bm{x}}a_ym\in P=UAM$, $\d\mu_P(p)=y^{-1}\, \d\bm{x}\d y\d\mu_M(m)$ is the right Haar measure of $P$ as in \rmkref{rmk:parbcase}.
\end{Lem}
\begin{proof}
To show that $P_{\e}$ is independent of the choice of $r$, we claim that
\begin{align*}
P_{\e}=\left\{p\in P: \|p\|_{\rm op}<1+\e\right\}.
\end{align*}
For this, note that $\tilde{\bm{\alpha}}_0=\bm{e}_0$ and for $p=u_{\bm{x}}a_ym\in P$, $\bm{e}_0p=y^{-1}\bm{e}_0$, implying that $\frac{\tilde{\bm{\alpha}}_0p}{\|\tilde{\bm{\alpha}}_0p\|}=\frac{\tilde{\bm{\alpha}}_0}{\sqrt{2}}$. This implies that the second defining condition of $G_{\e,r}(\bm{\alpha}_0)$ is always satisfied for any $r\in (0,1)$. Hence the claim holds. From the analysis in the previous proof, we see that $P_{\e/4}'\subseteq P_{\e}\subseteq P_{4\e}'$, where
\begin{align*}
P_{\e}':=G_{\e,r}'\cap P=\left\{u_{\bm{x}}a_ym\in P: \|\bm{x}\|<\e,\ |y-1|<\e,\ m\in M\right\}.
\end{align*}
The volume estimate then follows easily from this inclusion relation. The translation relation also follows easily from the translation relation in \lemref{lem:volest2} and the definition $P_{\e}=G_{\e,r}(\bm{\alpha}_0)\cap P$.
\end{proof}

Finally, for \thmref{thm:countmg} we choose another family of identity neighborhoods of $P$. 
It turns out it will be more convenient to write elements $p\in P$ in the form $p=ma_yu_{\bm{x}}$ with $m\in M$, $a_y\in A$ and $u_{\bm{x}}\in U$. (This can be done by expressing $p^{-1}$ using the Langlands decomposition $P=UAM$.) We note that under these coordinates 
\begin{align}\label{equ:righthaarP}
\d\mu_P(ma_yu_{\bm{x}})=y^{n-1}\, \d\mu_M(m)\d y\d\bm{x}.
\end{align}
We now introduce the family of identity neighborhoods we will be working with.
\begin{Lem}\label{lem:meaestp}
For any $\e\in (0,1)$ define $\tilde{P}_{\e}:=P'_{\e}\cap (P')_{\e}^{-1}$, where
\begin{align*}
P_{\e}':=\left\{ma_yu_{\bm{x}}\in P: m\in M,\ \|a_y\|_{\rm op}<1+\tfrac{\e}{4},\ \|\bm{x}\|<\tfrac{\e}{4}\right\}.
\end{align*}
Then we have $\mu_P(\tilde{P}_{\e})\asymp_n \e^{n+1}$ for any $0<\e<1$.
\end{Lem}
\begin{proof}
Using \eqref{equ:righthaarP} and recalling that $\|a_y\|_{\rm op}=\max\{y, y^{-1}\}$ it is easy to see that 
$$
\mu_P(P'_{\e})\asymp \e^{n+1},\quad \forall\ 0<\e<1.
$$
It thus suffices to show $P'_{\e/2}\subseteq \tilde{P}_{\e}$ for all $\e\in (0,1)$. Fix $\e\in (0,1)$. Since $\tilde{P}_{\e}=P'_{\e}\cap (P')_{\e}^{-1}$ and  $P'_{\e/2}\subseteq P'_{\e}$, we only need to show $P'_{\e/2}\subseteq (P')_{\e}^{-1}$. Take $p=ma_yu_{\bm{x}}\in P'_{\e/2}$, i.e. $m\in M$, $\|a_y\|_{\rm op}<1+\frac{\e}{8}$ and $\|\bm{x}\|<\frac{\e}{8}$, we would like to show $p^{-1}\in P'_{\e}$. Recall also that $m\in M$ is of the form $m=\left(\begin{smallmatrix}
1 & & \\
 & \tilde{m} & \\
 & & 1\end{smallmatrix}\right)$ for some $\tilde{m}\in \SO_n(\R)$. By direct computation we have 
\begin{align*}
p^{-1}=u_{-\bm{x}}a_{1/y}m^{-1}=m^{-1}a_{1/y}(ma_yu_{-\bm{x}}a_{1/y}m^{-1})=m'a_{y'}u_{\bm{x}'},
\end{align*}
with $m'=m^{-1}$, $y'=y^{-1}$ and $\bm{x}'=-y(\bm{x}\tilde{m}^{-1})$. Thus $\|a_{y'}\|_{\rm op}=\|a_y\|_{\rm op}<1+\frac{\e}{4}$ and 
$$
\|\bm{x}'\|=\|-y(\bm{x}\tilde{m}^{-1})\|=y\|\bm{x}\|<\left(1+\tfrac{\e}{8}\right)\times\tfrac{\e}{8}<\tfrac{\e}{4}.
$$
This implies $p^{-1}\in P'_{\e}$, finishing the proof. 
\end{proof}

\subsection{Verifying well-roundedness}
In this subsection, we verify well-roundedness for the family $\{S_{T,r,\bm{\alpha}}\}_{T>1}$ (with $r\in (0,1)$ and $\bm{\alpha}\in S^n$ fixed), the two-parameter family $\{S_{T,r,\bm{\alpha}_0}\}_{T>1,0<r<1}$ and the family $\{E_{\psi,T}\}_{T>1}$ with respect to the three corresponding families of identity neighborhoods introduced above. 
\comm{Recall from section \ref{sec:meantocount}, for any $\beta\in [1,2), c>0$}, $\cA_{\beta, c}$ is defined to be the family of Borel sets $B\subseteq \cV_{Q_n}^+$ which satisfy \eqref{e:DB} (with $\G\bk G=\tau^{-1}\SO_Q^+(\Z)\tau\bk \SO_{Q_n}^+(\R)$, $\cL=\tilde\cL_Q$, $\cV=\cV_{Q_n}^+$, $\mu_{\G}=\tilde\mu_Q$ and $m_{\cV}=\omega_Q m_{\cV_{Q_n}^+}$). 
As discussed in \rmkref{rmk:disbd}, there exists some constant $c>0$ depending only on $Q$ such that $\cA_{\beta_Q, c}$ contains all the generalized sectors, and all the differences of two \comm{nested} generalized sectors. Moreover, if $c_Q=0$ (so that $\beta_Q=1$), $\cA_{\beta_Q,c}$ contains all the finite-measure Borel sets in $\cV_{Q_n}^+$. We fix this constant $c$ and let $\beta=\beta_Q$ for the remaining of this section.

We first verify well-roundedness for the first two families.
\begin{Lem}\label{lem:wellround}
For any $\bm{\alpha}\in S^n$, for any $T>0$ and for any $\e, r\in (0, 1)$ we have
\begin{align}\label{equ:admisspro}
S_{(1-\e)T, (1-\e)r,\bm{\alpha}}\subseteq S_{T,r,\bm{\alpha}}h\subseteq S_{(1+\e)T,(1+\e)r,\bm{\alpha}},\quad \forall\ h\in G_{\e/6,r}(\bm{\alpha}).
\end{align}
\end{Lem}
\begin{remark}\label{rmk:webdsubgp}
When $\bm{\alpha}=\bm{\alpha}_0=(-1,\bm0)\in S^n$, since $P_{\e}\subseteq G_{\e,r}(\bm{\alpha}_0)$, the relations in \eqref{equ:admisspro} also holds for any $h\in P_{\e}$.
\end{remark}
\begin{proof}
First note that since for any $\bm{v}=\bm{e}_0a_yk_{\bm{\alpha}'}=y^{-1}\tilde{\bm{\alpha}}'\in \cV_{Q_n}^+$, $\frac{\bm{v}}{\|\bm{v}\|}=\frac{\tilde{\bm{\alpha}}'}{\|\tilde{\bm{\alpha}}'\|}=\frac{\tilde{\bm{\alpha}}'}{\sqrt{2}}$ we have
\begin{align*}
	S_{T, r, \bm{\alpha}}=\left\{\bm{v}\in \cV_{Q_n}^+: 0<v_{n+2}<T,\ \left\| \tfrac{\bm{v}}{\|\bm{v}\|} -\tfrac{\tilde{\bm{\alpha}}}{\sqrt{2}}\right\|<\tfrac{r}{\sqrt{2}}\right\}.
\end{align*}
Here $\tilde{\bm{\alpha}}'=(\bm{\alpha}',1)\in \cV_{Q_n}^+$. We now prove \eqref{equ:admisspro}. We first prove the second inclusion relation in \eqref{equ:admisspro}. Take any $\bm{v}\in S_{T,r, \bm{\alpha}}$ (i.e. $0<v_{n+2}<T$ and $\left\| \tfrac{\bm{v}}{\|\bm{v}\|} -\tfrac{\tilde{\bm{\alpha}}}{\sqrt{2}}\right\|<\frac{r}{\sqrt{2}}$) and $h\in G_{\e/6,r}(\bm{\alpha})$ (i.e. $\|h\|_{\rm op}< 1+\tfrac{\e}{6}$ and $\max\left\{\left\| \tfrac{\tilde{\bm{\alpha}}h}{\| \tilde{\bm{\alpha}}h \|}-\tfrac{\tilde{\bm{\alpha}}}{\sqrt{2}}\right\|, \left\| \tfrac{\tilde{\bm{\alpha}}h^{-1}}{\| \tilde{\bm{\alpha}}h^{-1} \|}-\tfrac{\tilde{\bm{\alpha}}}{\sqrt{2}}\right\|\right\}<\frac{r\e}{6}$), we need to show $\bm{v}h\in S_{(1+\e)T,r(1+\e), \bm{\alpha}}$.
Let $(\bm{v}h)_{n+2}$ be the $(n+2)$-th coordinate of $\bm{v}h$. Since $\|h\|_{\rm op}< 1+\frac{\e}{6}$ we get 
$$
0<(\bm{v}h)_{n+2}=\tfrac{\|\bm{v}h\|}{\sqrt{2}}=\left(1+\tfrac{\e}{6}\right)\tfrac{\|\bm{v}\|}{\sqrt{2}}<(1+\e)v_{n+2}<(1+\e)T.
$$
Moreover, we have $\left\| \tfrac{\bm{v}h}{\|\bm{v}\|} -\tfrac{\tilde{\bm{\alpha}}h}{\sqrt{2}}\right\|< \frac{(1+\tfrac{\e}{6})r}{\sqrt{2}}$. 
Now apply the geometric inequality  \eqref{equ:unitvecest} with $\bm{u}_1=\tfrac{\tilde{\bm{\alpha}}h}{\|\tilde{\bm{\alpha}}h\|}$, $\bm{u}_2=\tfrac{\bm{v}h}{ \|\bm{v}h\|}$ and 
$\lambda=\tfrac{\sqrt{2}\|\bm{v}h\|}{\|\tilde{\bm{\alpha}}h\| \|\bm{v}\|}$ 
to bound
\begin{align*}
\left\| \tfrac{\bm{v}h}{ \|\bm{v}h\|} -\tfrac{\tilde{\bm{\alpha}}h}{\|\tilde{\bm{\alpha}}h\|}\right\| \leq   &
\left(\tfrac{ \|\bm{v}\|\sqrt{2}}{\|\bm{v}h\| \|\tilde{\bm{\alpha}}h\|}\right)^{1/2} \left\| \tfrac{\tilde{\bm{\alpha}}h}{\sqrt{2}}-\tfrac{ \bm{v}h}{ \|\bm{v}\|}\right\| \leq  \tfrac{(1+\tfrac{\e}{6})^2r}{\sqrt{2}}.
\end{align*}
Next, using that $\left\|\tfrac{\tilde{\bm{\alpha}}h}{\|\tilde{\bm{\alpha}}h\|}-\tfrac{\tilde{\bm{\alpha}}}{\sqrt{2}}\right\|< \frac{r\e}{6}$ we can bound 
$$\left\| \tfrac{\bm{v}h}{ \|\bm{v}h\|} -\tfrac{\tilde{\bm{\alpha}}}{\sqrt{2}}\right\|\leq  \left((1+\tfrac{\e}{6})^2+\tfrac{\sqrt{2}\e}{6}\right)\tfrac{r}{\sqrt{2}}<\tfrac{(1+\e)r}{\sqrt{2}}.$$
This shows that $\bm{v}h\in S_{(1+\e)T,(1+\e)r, \bm{\alpha}}$, thus proving the second inclusion relation in \eqref{equ:admisspro}.
The first inclusion relation in \eqref{equ:admisspro} follows from the above analysis since for any $h\in G_{\e/6}(\bm{\alpha})$ (thus also $h^{-1}\in G_{\e/6}(\bm{\alpha})$)
\begin{align*}
S_{(1-\e)T,(1-\e)r, \bm{\alpha}}h^{-1}\subseteq S_{(1-\e^2)T,(1-\e^2)r, \bm{\alpha}}\subseteq S_{T,r, \bm{\alpha}},
\end{align*}
which is equivalent to the first inclusion relation in \eqref{equ:admisspro}. 
\end{proof}

\begin{Prop}\label{prop:wellrd}
For any $\bm{\alpha}\in S^n$ and for any $r\in (0,1)$, the family $\{S_{T,r,\bm{\alpha}}\}_{T>1}$ is strongly $\cA_{\beta,c}$-well rounded with respect to the family $\left\{G_{\e/6,r}(\bm{\alpha})\right\}_{0<\e<1}\subseteq G$. Similarly, when $\bm{\alpha}=\bm{\alpha}_0$, the family $\left\{S_{T,r,\bm{\alpha}_0}\right\}_{T>1,0<r<1}$ is $\cA_{\beta,c}$-well rounded with respect to the family $\left\{P_{\e/6}\right\}_{0<\e<1}$.
\end{Prop}
\begin{proof}
For the first half of this proposition, fix $\bm{\alpha}\in S^n$ and $r\in (0,1)$ and 
denote by $\cO_{\e}=G_{\e/6,r}(\bm{\alpha})$. First note that clearly for any $T>1$, $S_{T,r, \bm{\alpha}}$ is a generalized sector. Thus $\{S_{T,r, \bm{\alpha}}\}_{T>1}\subseteq \cA_{\beta,c}$. Next, by definition, we need to show that for any $B=S_{T,r, \bm{\alpha}}$ with $T>1$ and for any $\e\in (0, 0.1)$, there exists $\overline{B}_{\e}, \underline{B}_{\e}\in\cA_{\beta,c}$ satisfying that 
\begin{align*}
\underline{B}_{\e}\subseteq \bigcap_{h\in \cO_{\e}}Bh\subseteq \bigcup_{h\in \cO_{\e}}Bh\subseteq \overline{B}_{\e},\quad |\overline{B}_{\e}\setminus\underline{B}_{\e}|\ll_n \e|B|\quad\text{and}\quad \overline{B}_{\e}\setminus \underline{B}_{\e}\in \cA_{\beta,c}.
\end{align*}
For each $\e\in (0,0.1)$, we take $\underline{B}_{\e}=S_{(1-\e)T, (1-\e)r, \bm{\alpha}}$ and $\overline{B}_{\e}=S_{(1+\e)T, (1+\e)r, \bm{\alpha}}$. Then the above first relation follows from \lemref{lem:wellround}. Note also that both $\underline{B}_{\e}$ and $\overline{B}_{\e}$ are generalized sectors and $\overline{B}_{\e}\subseteq \underline{B}_{\e}$, thus $\overline{B}_{\e}\setminus \underline{B}_{\e}\in \cA_{\beta,c}$. It thus remains to show the above measure bound. For this note that
$
\overline{B}_\e\setminus \underline{B}_\e\subseteq   R_1\bigcup R_2
$
with 
\begin{align*}
	R_1:=\left\{\bm{e}_0a_yk_{\bm{\alpha}'}\in \cV_{Q_n}^+: 0<y^{-1}\leq (1-\e)T,\ \bm{\alpha}'\in \mathfrak{D}_{r(1+\e)}(\bm{\alpha})\setminus \mathfrak{D}_{r(1-\e)}(\bm{\alpha})\right\},
\end{align*}
and 
\begin{align*}
	R_2:=\left\{\bm{e}_0a_yk_{\bm{\alpha}'}\in\cV_{Q_n}^+: (1-\e)T\leq y^{-1}\leq (1+ \e)T,\ \bm{\alpha}'\in \mathfrak{D}_{r(1+\e)}(\bm{\alpha})\right\}.
\end{align*}
We can then apply \eqref{equ:differdiscs} and \eqref{equ:opendisfor} respectively to estimate that 
\begin{align*}
	|R_1|=\frac{\omega_Q}{n}(1-\e)^nT^n\sigma_n\left(\mathfrak{D}_{r(1+\e)}(\bm{\alpha})\setminus \mathfrak{D}_{r(1-\e)}(\bm{\alpha})\right)\ll_{Q,n}  \e T^nr^n\asymp_{Q,n} \e |B|,
\end{align*}
and 
\begin{align*}
	|R_2|=\frac{\omega_Q}{n}\left((1+\e)^n-(1-\e)^n\right)T^n\sigma_n\left(\mathfrak{D}_{r(1+\e)}(\bm{\alpha})\right)\ll_{Q,n} \e T^n r^{n}\asymp_{Q,n} \epsilon |B|.
\end{align*}
This proves the above bound, and hence also the first half of this proposition.

For the second half of this proposition, note that using similar arguments with \rmkref{rmk:webdsubgp} in place of \lemref{lem:wellround} we see that for any fixed $r\in (0,1)$, the family $\left\{S_{T,r,\bm{\alpha}_0}\right\}_{T>1}$ is $\cA_{\beta,c}$-well rounded with respect to $\left\{P_{\e/6}\right\}_{0<\e<1}$ (with the parameter $\e_0=\frac13$). Since the family $\{P_{\e/6}\}_{0<\e<1}$ is independent of $r$, we see that that larger family $\left\{S_{T,r,\bm{\alpha}_0}\right\}_{T>1,0<r<1}$ is $\cA_{\beta,c}$-well rounded with respect to $\{P_{\e/6}\}_{0<\e<1}$.
\end{proof}

Finally, we verify well-roundedness for the family $\{E_{\psi,T}\}_{T>1}$. 
\begin{Lem}\label{lem:stableperb}
Let $\{E_{\psi, T}\}_{T>1}$ be \comm{as in \eqref{equ:epsitset}  with $\psi: [0,\infty)\to (0,\infty)$ decreasing and continuous and let $\{\tilde{P}_{\e}\}_{0<\e<1}$ be as in \lemref{lem:meaestp}}. 
Then we have
\begin{align}\label{equ:inclrela}
E_{\psi^{-}_{\e},(1+\e)^{-1}T}\subseteq E_{\psi,T}h\subseteq E_{\psi^+_{\e},(1+\e)T},\quad \forall\ T>0,\ h\in \tilde{P}_{\e},\ \e\in (0, \tfrac12),
\end{align}
where 
$\psi^{\pm}_{\e}(t):=(1\pm\e)\psi((1+\e)^{\mp 1}t)$ for any $t\geq 0$.
\end{Lem}
\begin{proof}
Throughout the proof, for any $\bm{v}\in\cV_{Q_n}^+$, we write it in the form $\bm{v}=(v_1,\bm{w},v_{n+2})\in \R\times \R^n\times \R_{>0}$. 
For any $g\in G$, we denote by $(v_1(g),\bm{w}(g), v_{n+2}(g)):=\bm{v}g$. 
We first note that $E_{\psi,T}$ has the following alternative description that for any $T>1$ 
\begin{align*}
	E_{\psi,T}
	&=\left\{\bm{v}\in \cV_{Q_n}^+:  2(v_1+v_{n+2})<v_{n+2}\psi^2(v_{n+2}),\ 0< v_{n+2}<T\right\}.
\end{align*}

We now prove the second inclusion relation in \eqref{equ:inclrela}. Fix $T>0$ and $\e\in (0, \tfrac12)$. Take $\bm{v}\in E_{\psi, T}$ and $h\in \tilde{P}_{\e}$, we would like to show $\bm{v}h\in E_{\psi^+_{\e},(1+\e)T}$, that is
\begin{align}\label{equ:firstgoal}
0< v_{n+2}(h)<(1+\e)T,
\end{align} 
and 
\begin{align}\label{equ:secondgoal}
2(v_{n+2}(h)+v_{n+1}(h))<(1+\e)^2v_{n+2}(h)\psi^2\left((1+\e)^{-1}v_{n+2}(h)\right).
\end{align} 
By direct computation we have
\begin{align*}
v_{1}(h)=-\bm{w}\cdot \bm{x}+\tfrac{y}{2}(1-\|\bm{x}\|^2)(v_{n+2}+v_1)-\tfrac{1}{2y}(v_{n+2}-v_1), 
\end{align*}
and
\begin{align*}
 v_{n+2}(h)=\bm{w}\cdot \bm{x}+\tfrac{y}{2}(1+\|\bm{x}\|^2)(v_{n+2}+v_1)+\tfrac{1}{2y}(v_{n+2}-v_1).
\end{align*}
Since $h\in \tilde{P}_{\e}\subseteq P'_{\e}$, we can write $h=ma_yu_{\bm{x}}$ with $m\in M$, $(1+\e/4)^{-1}<y<1+\e/4$ and $\|\bm{x}\|<\e/4$. 
For \eqref{equ:firstgoal} note that $v_{n+2}(h)>0$ is automatically satisfied since $\bm{v}h\in \cV_{Q_n}^+$; we thus only need to prove $v_{n+2}(h)<(1+\e)T$. For this, by rewriting $v_{n+2}-v_{1}=2v_{n+2}-(v_{n+2}+v_{1})$ we get
\begin{align*}
v_{n+2}(h)&=\tfrac{v_{n+2}}{y}+\bm{w}\cdot \bm{x}+\left(\tfrac{y}{2}(1+\|\bm{x}\|^2)-\tfrac{1}{2y}\right)(v_{n+2}+v_{1}),
\end{align*}
implying that (using the estimates $(1+\e/4)^{-1}<y<1+\e/4$, $\|\bm{x}\|<\e/4$ and $0<\e<\tfrac12<1$) 
\begin{align*}
\left|v_{n+2}(h)-\tfrac{v_{n+2}}{y}\right|<\tfrac{\|\bm{w}\|\e}{4}+\tfrac{5\e}{16}(v_{n+2}+v_1)=\tfrac{\e}{4}(\|\bm{w}\|+v_1)+\tfrac{\e}{16}v_1+\tfrac{5\e}{16}v_{n+2}.
\end{align*}
Now using the estimates $(\|\bm{w}\|+v_1)^2\leq 2(\|\bm{w}\|^2+v_1^2)=2v_{n+2}^2$ and $v_1\leq v_{n+2}$ we can estimate
\begin{align}
	\left|v_{n+2}(h)-\tfrac{v_{n+2}}{y}\right|\leq \left(\tfrac{\sqrt{2}}{4}+\tfrac{1}{16}+\tfrac{5}{16}\right)\e v_{n+2}<\tfrac{3\e}{4}v_{n+2},
\end{align}
which then further implies (using also $1-\frac{\e}{4}<(1+\frac{\e}{4})^{-1}<y<1+\frac{\e}{4}$)
\begin{align}\label{equ:keyest}
 \left(1-\e\right)v_{n+2}<v_{n+2}(h)<\left(1+\e\right)v_{n+2}.
\end{align}
This then immediately implies $v_{n+2}(h)<(1+\e)v_{n+2}<(1+\e)T$ as desired.

For \eqref{equ:secondgoal} note
\begin{align*}
v_{n+2}(h)+v_{1}(h)=y(v_{n+2}+v_{1}).
\end{align*}
Thus using the assumptions that $\bm{v}\in E_{\psi, T}$ (so that $2(v_{n+2}+v_1)<v_{n+2}\psi^2(v_{n+2})$) and that $\psi(\cdot)$ is decreasing and the estimates \eqref{equ:keyest} and $(1+\e/4)(1-\e)^{-1}<(1+\e)^2$ (this is true since $\e\in (0, \tfrac12)$) 
we have
\begin{align*}
2(v_{n+2}(h)+v_{1}(h))&<\left(1+\tfrac{\e}{4}\right)v_{n+2}\psi^2(v_{n+2})<(1+\e)^2v_{n+2}(h)\psi^2\left((1+\e)^{-1}v_{n+2}(h)\right),
\end{align*}
as desired. This finishes the proof of the second inclusion relation in \eqref{equ:inclrela}.

For the other inclusion relation, take any $h\in \tilde{P}_{\e}$, we need to show $E_{\psi^{-}_{\e},(1+\e)^{-1}T}\subseteq E_{\psi,T}h$, or equivalently, $E_{\psi^{-}_{\e},(1+\e)^{-1}T}h^{-1}\subseteq E_{\psi,T}$. This then follows from similar analysis and noting that $h^{-1}\in P'_{\e}$ (since $\tilde{P}_{\e}=P'_{\e}\cap (P')_{\e}^{-1}$). 
\end{proof}
\begin{Prop}\label{prop:wellbdsugp2}
Let $Q$ and $\psi: [0,\infty)\to (0,\frac12)$ be as in \thmref{thm:countmg}. Then the family $\{E_{\psi, T}\}_{T>1}$ is $\cA_{1,c}$-well rounded with respect to $\{\tilde{P}_{\e}\}_{0<\e<1}$.
\end{Prop}
\begin{proof}
We show that for any $B=E_{\psi, T}$ with $T>1$ and for any $\e\in (0, \tfrac12)$ there exist $\overline{B}_{\e},\underline{B}_{\e}\in \cA_{1,c}$ satisfying 
\begin{align*}
\underline{B}_{\e}\subseteq \bigcap_{h\in \tilde{P}_{\e}}Bh\subseteq \bigcup_{h\in \tilde{P}_{\e}}Bh\subseteq \overline{B}_{\e}\quad\text{and}\quad |\overline{B}_{\e}\setminus\underline{B}_{\e}|\ll_{Q,n} \e|B|.
\end{align*}
For this, we take $\underline{B}_{\e}=E_{\psi^{-}_{\e},(1+\e)^{-1}T}$ and $\overline{B}_{\e}=E_{\psi^{+}_{\e},(1+\e)T}$, where $\psi^{\pm}_{\e}$ are as in \lemref{lem:stableperb}. 
Clearly, both $\underline{B}_{\e}$ and $\overline{B}_{\e}$ are Borel sets with finite measure, and thus are contained in $ \cA_{1,c}$. The above inclusion relations follow from \lemref{lem:stableperb}. For the measure bound, we note that $\overline{B}_{\e}\setminus\underline{B}_{\e}\subseteq R_{1}\bigcup R_{2}$ with
\begin{align*}
R_1:=\left\{\bm{e}_0a_yk_{\bm{\alpha}}\in \cV_{Q_n}^+: \psi^{-}_{\e}(y^{-1})\leq \|\bm{\alpha}_0-\bm{\alpha}\|\leq \psi^+_{\e}(y^{-1}),\ 0< y^{-1}\leq (1+\e)^{-1}T\right\},
\end{align*}
and 
\begin{align*}
R_2:=\left\{\bm{e}_0a_yk_{\bm{\alpha}}\in \cV_{Q_n}^+: \|\bm{\alpha}_0-\bm{\alpha}\|\leq \psi^+_{\e}(y^{-1}),\ (1+\e)^{-1}T\leq y^{-1}\leq (1+\e)T\right\}.
\end{align*}
By \lemref{lem:volumecyl} we have (noting also that by our assumptions $\psi_{\e}^+(t)<1$ for all $t\in [0,\infty)$) 
\begin{align*}
|R_1|\ll_{Q,n}\int_{0}^{(1+\e)^{-1}T}t^{n-1}\left(\psi^+_{\e}(t)^n-\psi^-_{\e}(t)^n\right)\, \d t \quad \text{and}\quad |R_2|\ll_{Q,n} \int_{(1+\e)^{-1}T}^{(1+\e)T} t^{n-1}\psi^+_{\e}(t)^n\, \d t.
\end{align*}
This, together with the measure estimate $|B|\asymp_{Q,n} \cJ_{\psi}(T)$ (cf. \lemref{lem:volcompqdioapp}) implies that
\begin{align*}
|\overline{B}_{\e}\setminus \underline{B}_{\e}|&\ll_{Q,n}\int_{0}^{(1+\e)T}t^{n-1}\psi^+_{\e}(t)\, \d t- \int_{0}^{(1+\e)^{-1}T}t^{n-1}\psi^-_{\e}(t)\, \d t\\
&=\left((1+\e)^{2n}-(1-\e)^{2n}\right)\int_0^Tt^{n-1}\psi(t)^{n}\, \d t
\ll_{Q,n}\e|B|.
\end{align*}
This finishes the proof.
\end{proof}

\subsection{Proofs of counting results}
\comm{We now combine all the results from previous sections to prove Theorems \ref{thm:equidis}-\ref{thm:countmg} which then imply Theorems \ref{thm:c1digen}-\ref{thm:khinttypegen} respectively.}

In view of \lemref{lem:volest2} and \propref{prop:wellrd}, 
\thmref{thm:equidis} now easily follows from \thmref{thm:cfamd} as follows.
\begin{proof}[Proof of \thmref{thm:equidis}]
Fix $\bm{\alpha}\in S^n$ and $r\in (0,1)$. 
Let $\cB=\left\{S_{T,r,\bm{\alpha}}\right\}_{T>1}$ and $\cO=\left\{\cO_{\e}\right\}_{0<\e<1}$ with $\cO_{\e}=G_{\e/6,r}(\bm{\alpha})$. Then by \propref{prop:wellrd}, $\cB$ is strongly $\cA_{\beta,c}$-well rounded with respect to $\cO$, and by \lemref{lem:volest2}, $\cO$ satisfies the assumption \eqref{e:OeReg} with $\delta=c_n'r^n$ 
(for some $0<c'_n<1$ depending only on $n$), $d=2n+1$ and $\eta=3$. In particular, we have $0<\delta<1$. Next, note that by \eqref{equ:sectvolest} we have
 $|S_{T,r,\bm{\alpha}}|\gg_{Q,n} T^nr^n$. 
Thus the condition $|S_{T,r,\bm{\alpha}}|> \delta^{-\frac{1}{2-\beta}}$ is satisfied whenever 
$T^{-\frac{2-\beta}{3-\beta}}\ll_{Q,n} r$.
Hence for 
$T^{-\frac{2-\beta}{3-\beta}}\ll_{Q,n} r<1$
we can apply \thmref{thm:cfamd} to get for any $g\in \comm{G}$,
\begin{align*}
\#\left(\tilde\cL_Q g\cap S_{T,r,\bm{\alpha}}\right)=|S_{T,r,\bm{\alpha}}|+ O_{Q,g}\left( r^{-\frac{n}{2n+3}}|S_{T,r,\bm{\alpha}}|^{1-\frac{2-\beta}{2n+3}}\right),
\end{align*}
finishing the proof.
\end{proof}
\begin{proof}[Proof of \thmref{thm:c1digen}]
\comm{Apply \thmref{thm:equidis} with $g=\text{id}$ and the relation \eqref{equ:relaticount} we get for any $r,T>0$ satisfying $T^{-\frac{2-\beta}{3-\beta}}\ll_{Q,n} r<1$ and for any $\bm{x}\in S_{\cQ}$,
\begin{align*}
\cN_{\cQ}(\bm{x}, r,T)&=|S_{T,r,\bm{\alpha}}|+ O_{Q}\left( r^{-\frac{n}{2n+3}}|S_{T,r,\bm{\alpha}}|^{1-\frac{2-\beta}{2n+3}}\right),
\end{align*}
where $\bm{\alpha}=\bm{x}\tilde{\tau} \in S^n$.
Then this theorem follows from the first measure estimate in \eqref{equ:sectvolest} for $S_{T,r,\bm{\alpha}}$ and the relations $|S_{T,r,\bm{\alpha}}|=\omega_Qm_{\cV_{Q_n}^+}(S_{T,r,\bm{\alpha}})$ and $\sigma_n(\mathfrak{D}_r(\bm{\alpha}))=\sigma_{\cQ}(\mathfrak{D}^{\cQ}_r(\bm{x}))$.}
\end{proof}

Similarly, we can prove \thmref{thm:subtranspara} by applying \thmref{thm:countnoninc}.
\begin{proof}[Proof of \thmref{thm:subtranspara}]
Let $\cB=\left\{B_{t,r}\right\}_{t>1, 0<r<1}$, where $B_{t,r}:=S_{t^{1/n}, r^{1/n},\bm{\alpha}_0}$. By the second half of \propref{prop:wellrd} we know that $\cB$ is $\cA_{\beta,c}$-well rounded with respect to $\cO=\{P_{\e/6}\}_{0<\e<1}\subseteq P$. Next, we verify that $B_{t,r}$ satisfies the measure assumption in \thmref{thm:nonestgen}. First, by \eqref{equ:sectvolest} we have 
\begin{align}\label{equ:volcomps}
|B_{t,r}|=|S_{t^{1/n}, r^{1/n},\bm{\alpha}_0}|=\varkappa_{\cQ}tr+O_{Q,n}(tr^{\frac{n+2}{n}})\asymp_{Q,n} tr.
\end{align}
Moreover, by similar analysis as in the proof of \propref{prop:wellrd}, we see that for any $t_1>t_2>1$ and $0<r_2<r_1<1$,
\begin{align*}
\left|B_{t_1,r_1}\setminus B_{t_2, r_2}\right|
&\ll_n t_2\sigma_n\left(\mathfrak{D}_{r_1^{1/n}}(\bm{\alpha}_0)\setminus \mathfrak{D}_{r_2^{1/n}}(\bm{\alpha}_0) \right)+(t_1-t_2)\sigma_n\left(\mathfrak{D}_{r_1^{1/n}}(\bm{\alpha}_0)\right)\\
&\ll_n t_1(r_1-r_2)+(t_1-t_2)r_1.
\end{align*}
This verifies condition \eqref{r:mestimate}. Thus for case (1) we can apply case (1) of \thmref{thm:countnoninc} with $H= P$, $\iota_P=\mathfrak{s}$, $\mu_{P\bk G}=\sigma_n$ (cf. \rmkref{rmk:parbcase}) and $d=n+1$ (cf. \lemref{lem:idnbsubgp}) to get for any $p\in P$, for $\sigma_n$-a.e. $\bm{\alpha}\in S^n$ and for $t$ sufficiently large
\begin{align*}
\#\left(\tilde\cL_Q k_{\bm{\alpha}}^{-1}p^{-1}\cap B_{t, t^{-\lambda}}\right)=|B_{t, t^{-\lambda}}|+O\left(|B_{t, t^{-\lambda}}|^{1-\frac{2-\beta}{n+4}}\log \left(|B_{t,t^{-\lambda}}|\right)\right).
\end{align*}
The counting formula in case (1) then follows after making a change of variable $T=t^{1/n} $.
Case (2) follows similarly using case (2) of \thmref{thm:countnoninc}.
\end{proof}
\begin{proof}[Proof of  \thmref{thm:countforgenericgen}]
Case (1) and (2) of \thmref{thm:countforgenericgen} follow by applying case (1) and (2) of \thmref{thm:subtranspara} respectively (with $p=\text{id}$), together with the relation \eqref{equ:relaticount} and the fact that the right $\tilde{\tau}^{-1}$-multiplication map from $S^n$ to $S_{\cQ}$ sends a full $\sigma_n$-measure set in $S^n$ to a full $\sigma_{\cQ}$-measure set in $S_{\cQ}$.
\end{proof}
\begin{proof}[Proof of \thmref{thm:countmg}]
Let $\cB=\left\{E_{\psi, T}\right\}_{T>1}$ \comm{and note that} $\cB$ is an increasing family. Moreover, the function $T\mapsto |E_{\psi, T}|$ is continuous \comm{and unbounded} in view of \lemref{lem:volcompqdioapp} and the \comm{assumptions that $\psi$ is continuous and} $\int_0^{\infty} t^{n-1}\psi(t)^n\, \d t=\infty$. Thus $\cB$ satisfies $\left\{|B|: B\in \cB\right\}\supseteq (V,\infty)$ for some $V>0$. Moreover, by \propref{prop:wellbdsugp2}, $\cB$ is $\cA_{1,c}$-well rounded with respect to $\{\tilde{P}_{\e}\}_{0<\e<1}$. Thus we can apply \thmref{thm:countforall} (with $\beta=1$ and $d=n+1$ (cf. \lemref{lem:meaestp})) to get for any $p\in P$, for $\sigma_n$-a.e. $\bm{\alpha}\in S^n$ and for all sufficiently large $T$,
\begin{align*}
\#\left(\tilde\cL_Q k_{\bm{\alpha}}^{-1}p^{-1}\cap E_{\psi, T}\right)=|E_{\psi, T}|+O_Q\left(|E_{\psi,T}|^{\frac{n+3}{n+4}}\log\left(|E_{\psi, T}|\right)\right).
\end{align*}
This finishes the proof.
\end{proof}

\begin{proof}[Proof of \thmref{thm:khinttypegen}]
We first assume $\psi(\N)\subseteq (0,\frac12)$. With slight abuse of notation, we still denote by ${\psi}: [0,\infty)\to \R_{>0}$ a decreasing continuous extension of $\psi: \N\to \R_{>0}$ satisfying that $\psi(t)=\psi(1)$ for all $0\leq t<1$. {In particular, we have $\psi([0,\infty))\subseteq (0,\frac12)$.} Moreover, since this continuous extension is decreasing, the assumptions $\lim_{q\to\infty}\psi(q)=0$ and $\sum_{q\in \N}q^{n-1}\psi(q)^n=\infty$ imply that $\lim\limits_{t\to\infty}{\psi}(t)=0$ and $\int_0^{\infty}t^{n-1}{\psi}(t)^n\, \d t=\infty$. We can thus apply \thmref{thm:countmg} (with $p=\text{id}$), the relation \eqref{equ:relaticount2} and the fact that the right $\tilde{\tau}^{-1}$-multiplication map from $S^n$ to $S_{\cQ}$ sends a full $\sigma_n$-measure set in $S^n$ to a full $\sigma_{\cQ}$-measure set in $S_{\cQ}$ to get for $\sigma_{\cQ}$-a.e. $\bm{x}\in S_{\cQ}$ and for all sufficiently large $T$,
\begin{align*}
\cN_{\cQ,\psi}(\bm{x}, T)=|E_{\psi, T}|+ O_Q\left(|E_{\psi,T}|^{\frac{n+3}{n+4}}\log(|E_{\psi, T}|)\right).
\end{align*}
Next, applying the measure estimate \eqref{equ:volume1} for $E_{{\psi},T}$, together with \lemref{lem:esthold} and the simple estimates $\cJ_{\psi}(T)\asymp J_{\psi}(T)$ and $\cI_{{\psi}}(T)\asymp I_{\psi}(T)$ we get
\begin{align*}
\cN_{\cQ,\psi}(\bm{x}, T)&=n\varkappa_{\cQ}\cJ_{\psi}(T)+O_{\cQ}(\cJ_{\psi}(T)^{\frac{n+3}{n+4}}\log(\cJ_{\psi}(T))+\cI_{\psi}(T))\\
&=n\varkappa_{\cQ}J_{\psi}(T)+O_{\cQ,\psi}\left(J_{\psi}(T)^{\frac{n+3}{n+4}}\right)+ O_{\cQ}(J_{\psi}(T)^{\frac{n+3}{n+4}}\log(J_{\psi}(T))+I_{\psi}(T)).
\end{align*}
We can then finish the proof by absorbing the first error term into the second. 

{Finally, for a general $\psi: \N\to (0,\infty)$, since $\lim\limits_{q\to\infty}\psi(q)=0$, there exists some $q_0\in\N$ such that $\psi(q)<\frac12$ for all $q\geq q_0$. Define $\psi_0: \N\to (0,\frac12)$ by $\psi_0(q):=\psi(q_0)$ if $q<q_0$ and $\psi_0(q)=\psi(q)$ if $q\geq q_0$. Applying the previous result for $\psi_0$ we get for $\sigma_{\cQ}$-a.e. $\bm{x}\in S_{\cQ}$ and for all sufficiently large $T$,
\begin{align*}
\cN_{\cQ,\psi_0}(\bm{x}, T)&=n\varkappa_{\cQ}J_{\psi_0}(T)+O_{\cQ,n,\psi}(J_{\psi_0}(T)^{\frac{n+3}{n+4}}\log(J_{\psi_0}(T))+I_{\psi_0}(T)).
\end{align*}
We can then conclude the proof by noting that $\cN_{\cQ,\psi}(\bm{x}, T)=\cN_{\cQ,\psi_0}(\bm{x}, T)+O_{\psi}(1)$, $J_{\psi}(T)=J_{\psi_0}(T)+O_{\psi}(1)$ and $I_{\psi}(T)=I_{\psi_0}(T)+O_{\psi}(1)$ for any $T>1$ and $\bm{x}\in S_{\cQ}$.}
\end{proof}

\bibliographystyle{alpha}

\end{document}